\documentclass[11pt]{amsart}
 
\input{xy} 
\xyoption{all}
\usepackage{graphicx}
\usepackage{color}

\newtheorem{theorem}{Theorem} [section]
\newtheorem{prop}[theorem]{Proposition} 
\newtheorem{lemma}[theorem]{Lemma}
\newtheorem{cor}[theorem]{Corollary}

\numberwithin{equation}{section} 
\numberwithin{figure}{section}

\newcommand\C{{\mathbb C}}

\newcommand\R{{\mathbb R}}
\newcommand\Z{{\mathbb Z}}

\newcommand\N{{\mathbb N}}
\newcommand\D{{\mathbb D}}
\newcommand\Hyp{{\mathbb H}}
\newcommand\cM{\mathcal{M}}
\newcommand\cC{\mathcal{C}}
\newcommand\cD{\mathcal{D}}
\newcommand{\cT}{{\mathcal T}}
\newcommand{\cF}{{\mathcal F}}

\newcommand{\cR}{\mathcal{R}}
\newcommand{\cL}{\mathcal{L}}
\newcommand{\cP}{\mathcal{P}}
\newcommand{\cS}{\mathcal{S}}
\newcommand{\cX}{\mathcal{X}}
\newcommand{\cB}{\mathcal{B}}
\newcommand\eps{\varepsilon}

\renewcommand\phi{\varphi}
\newcommand\Aut{\mathrm{Aut}}
\newcommand\Stab{\mathrm{Stab}}

\newcommand\iso{\simeq} %isomorphism -~

\newcommand\modspace {\mod} %\mod with space before
\renewcommand\mod{\operatorname{mod}}  %\mod without extra space
\renewcommand\gcd {\operatorname{gcd}} %gcd
\newcommand\lcm {\operatorname{lcm}}  %lcm

\newcommand\del{\partial} 

\newcommand\<{\langle} %innerproduct brackets 
\renewcommand\>{\rangle} %%

\renewcommand\Im {\operatorname{Im}} 
\newcommand\Top {\mathrm{Top}}
\newcommand \Tw {\operatorname{Twist}}
\newcommand\TP {\mathrm{TP}}

\newcommand\Sol {\mathrm{Sol}}

\newcommand{\IS}{\mathbb{S}}

\begin{document}

\title{The classification of polynomial basins of infinity}

\author{Laura DeMarco and Kevin Pilgrim}

\begin{abstract}  
We consider the problem of classifying the dynamics of complex polynomials $f: \C \to \C$ restricted to their basins of infinity.  We synthesize existing combinatorial tools --- tableaux, trees, and laminations --- into a new invariant of basin dynamics we call the pictograph.  For polynomials with all critical points escaping to infinity, we obtain a complete description of the set of topological conjugacy classes.  We give an algorithm for constructing abstract pictographs, and we provide an inductive algorithm for counting topological conjugacy classes with a given pictograph.
\end{abstract}

\date{\today}

\thanks{We thank Jan Kiwi and Curt McMullen for helpful discussions.  Our research was supported by the National Science Foundation and the Sloan Foundation.  MSC2010:  37F10, 37F20.}

\maketitle

%\tableofcontents

\thispagestyle{empty}

%%%%%%
%%%%%%
\section{Introduction}

This article continues a study of the moduli space of complex polynomials $f:\C\to\C$, in each degree $d\geq 2$, in terms of the dynamics of polynomials on their basins of infinity.  Our main goal is the development of combinatorial tools to classify the topological conjugacy classes of a polynomial $f$ restricted to its basin
	$$X(f) = \{z\in\C: f^n(z)\to\infty\}.$$
The basin $X(f)$ is an open, connected subset of $\C$.  Further, when all critical points of $f$ lie in $X(f)$, the basin is a rigid Riemann surface, admitting a unique embedding to $\C$ (up to affine transformations).  In this case, the restriction $f: X(f)\to X(f)$ uniquely determines the conformal  conjugacy class of $f: \C\to\C$.  Thus, our results provide a combinatorial classification of topological conjugacy classes of polynomials in the shift locus.  	

We introduce the {\em pictograph} as a combinatorial diagram, possibly infinite, to encode the dynamical system $(f, X(f))$.  It combines discrete data with analytic information.  Our main result can be stated as follows:

\begin{theorem}  \label{maintheorem1}
The pictograph is a topological-conjugacy invariant.  For any given pictograph $\cD$, the number of topological conjugacy classes $\Top(\cD)$ of basins $(f, X(f))$ with pictograph $\cD$ is algorithmically computable from the discrete data of $\cD$.
\end{theorem}

\noindent
We provide the ingredients for an inductive, algorithmic computation of $\Top(\cD)$, though we give the full details only in degree 3 (see Theorems \ref{deg3twist} and \ref{deg3infinite}) where the computations are more tractable.  The computation is achieved by an analysis of the quasiconformal twist deformations on the basin of infinity, as introduced in \cite{Branner:Hubbard:1} and \cite{McS:QCIII}.  We provide examples both with $\Top(\cD) = 1$ and with $\Top(\cD)>1$.  Indeed, in every degree $d>2$, there exist examples with $\Top(\cD)>1$.   Even though the pictograph is not a complete invariant, it ``knows" its failure.

As we shall see, a pictograph can be defined and constructed abstractly.  We prove that every abstract pictograph arises from some polynomial.  Moreover, there is a natural moduli space $\cB_d$ parameterizing the restrictions of polynomials to their basins of infinity \cite{DP:basins}, and once the pictograph and critical escape rates are fixed, the locus in $\cB_d$ with this data admits the following description (see \S\ref{sec:conjugacy}). 

\begin{theorem} \label{maintheorem2}
Fix a pictograph $\cD$ and a list of $N$ compatible critical escape rates.  Then  the locus in $\cB_d$ of basin dynamical systems with this data is a compact locally trivial fiber bundle over the torus $(S^1)^N$ with totally disconnected fibers.   The total space is foliated by $N$-manifolds, and the leaves are in bijective correspondence with topological conjugacy classes of basin dynamics within the total space.  Over the shift locus, the fibers are finite. 
\end{theorem}  

\noindent
In a sequel to this paper, we will show that this bundle has a natural group structure, and that the count of topological conjugacy classes has a purely algebraic interpretation.

\bigskip\noindent{\bf Motivation and context.}
This article continues a series of articles (\cite{DP:basins}, \cite{DP:heights}) in which we study polynomial dynamics by concentrating on the basin of infinity.  Our goal has been to understand the structure and organization of topological conjugacy classes within the moduli space $\cM_d$ of conformal conjugacy classes.  In degree $d=2$, there are only two topological conjugacy classes of basins $(f, X(f))$, distinguished by the Julia set being connected or disconnected \cite[Theorem 10.1]{McS:QCIII}.  In every degree $>2$, there are infinitely many topological conjugacy classes of basins, even among the structurally stable polynomials in the shift locus.  

Our methods and perspective are inspired by the two foundational articles of Branner and Hubbard on polynomial dynamics which lay the groundwork and treat the case of cubic polynomials in detail.  In \cite{Branner:Hubbard:1}, critical escape rates and the quasiconformal deformation of the basin of infinity play an important role, while \cite{Branner:Hubbard:2} introduces combinatorial and geometric invariants.  In this article, we build upon the methods of \cite{Branner:Hubbard:2} and \cite{DM:trees} (compare also \cite{Blanchard:Devaney:Keen}, \cite{Perez:frames}, \cite{Kiwi:combinatorial}, \cite{Milnor:cubics}) to classify conjugacy classes by studying the recurrent behavior of critical points.   While the Branner-Hubbard tableau (or equivalently, the Yoccoz $\tau$-sequence) records the first-return map along the critical nest, the {\em pictograph} we define records the first return to a ``decorated" critical nest.  The definition works in all degrees.  Viewed as an invariant of a polynomial dynamical system, the pictograph synthesizes the pattern and tableau of Branner and Hubbard, the metric tree equipped with dynamics of DeMarco and McMullen, and the laminations of Thurston. 
 
 In \cite{DP:heights}, we defined a decomposition of the moduli space of polynomials $\cM_d$ in terms of the critical escape rates.  Passing to the corresponding quotient $\cM_d \to \cT_d^*$, the image $\cT_d^*$ has, on a dense open subset corresponding to the shift locus, the structure of a cone over a locally-finite simplicial complex.  The top-dimensional simplices are in one-to-one correspondence with topological conjugacy classes of structurally stable polynomials in the shift locus.  The combinatorics defined here can be used to encode the simplices of the complex $\cT_d^*$.  The algorithms we provide in this article are implemented in degree $d=3$ in \cite{DS:count} and \cite{DS:euler}.  
 
In the theory of dynamical systems, the study of a system like $(f, X(f))$ is somewhat nonstandard.  On the one hand, since all points tend to $\infty$ under iteration, the system is transient.  On the other hand, the structure of $(f, X(f))$, with an induced dynamical system on its Cantor set of ends, carries enough information to recover the full entropy of the polynomial $(f, \C)$ (see \cite[Theorem 1.1]{DM:trees}).  Further, it has been shown recently that the conformal class of $(f, X(f))$ uniquely determines the conformal class of $(f, \C)$ unless there is a critical point in a periodic end of $X(f)$ \cite{Yin:Zhai}.

\bigskip\noindent{\bf Structure of the article.}
The article is divided into four parts:

\medskip\noindent
{\bf Part I.  Local structure.}  We begin in Section \ref{sec:laminations} by studying local models for polynomial branched covers.  We introduced the local model surface and local model map in \cite{DP:basins}.  Here, we show that the conformal structure of a local model surface is recorded by a {\em finite lamination}.  The local model map, up to symmetries, can be recovered from the domain lamination and its degree.  

\medskip\noindent
{\bf Part II.  The tree of local models.}  Throughout Part II (Sections \ref{sec:trees}--\ref{sec:symmetries}), we fix one polynomial $f: \C\to\C$ of degree $d\geq 2$, and we study the restriction to its basin of infinity $X(f)$.  We may assume that at least one critical point of $f$ lies in $X(f)$; otherwise, the Julia set of $f$ is connected, and the restriction $(f, X(f))$ is conformally conjugate to $(z^d, |z|>1)$.  We begin with a review of the polynomial tree $(F, T(f))$, defined in \cite{DM:trees}, and we introduce the {\em tree of local models} $(\cF, \cX(f))$ associated to $f$.  There are semiconjugacies
	$$(f, X(f)) \to (F, T(f))$$
induced by a natural quotient map, and 
	$$(\cF, \cX(f)) \to (f, X(f))$$
induced by a natural  {\em gluing quotient map}.  We define the {\em spine} of the tree of local models and show (Proposition \ref{tlm spine}) that the dynamical system $(\cF, \cX(f))$ is determined by the first-return map on its spine.  

Both trees and trees of local models can be defined abstractly, from a list of axioms.  Following the proof of the realization theorem for trees \cite[Theorem 1.2]{DM:trees}, we prove the analogous realization theorem (Theorem \ref{tlm realization}):  every abstract tree of local models comes from a polynomial.    

\medskip\noindent
{\bf Part III.  The moduli space and topological conjugacy.}
In Part III (Section \ref{sec:conjugacy}), we study our dynamical systems in families.  We begin by recalling facts about the quasiconformal deformation theory of polynomials from \cite{McS:QCIII}, specifically the twisting and stretching deformations on the basin of infinity.  These quasiconformal conjugacies parametrize the topological conjugacy classes of basins. We show (Theorem \ref{tlm invariance}) that the tree of local models is a twist-conjugacy invariant.

The moduli space $\cM_d$ of conformal conjugacy classes of polynomials has the structure of a complex orbifold, viewed as the space of complex coefficients modulo the conjugation action by affine transformations. The moduli space $\cB_d$ of polynomial basins $(f,X(f))$ has a natural Gromov-Hausdorff topology coming from a flat conformal metric on the basin. By the main result of \cite{DP:basins}, the projection $\cM_d \to \cB_d$ is continuous and proper, with connected fibers, and a homeomorphism on the shift locus.  We examine the structure of the subset $\cB_d(\cF, \cX)$ of basins in $\cB_d$ with a given tree of local model models, and we provide ingredients for the proof of Theorem \ref{maintheorem2}.  We show $\cB_d(\cF, \cX)$ is a compact, locally trivial fiber bundle over a torus with totally disconnected fibers.  The twisting deformation induces the local holonomy of the fiber bundle.  Indeed, the orbits of the twisting deformation are in one-to-one correspondence with the topological conjugacy classes of polynomial basins with the given tree of local models $(\cF, \cX)$.  The counting of these topological conjugacy classes (to be discussed in detail in Part IV) is done via an analysis of the monodromy of the twisting action on $\cB_d(\cF, \cX)$.

\medskip\noindent
{\bf Part IV. Combinatorics and counting.}
In Part IV (Sections \ref{sec:spine}--\ref{sec:counting}), we introduce our main object, the pictograph.  Given a tree of local models, its pictograph is the collection of lamination diagrams along its spine, labelled by the orbits of the critical points.  The tree of local models can be reconstructed from its pictograph and the list of critical escape rates (Proposition \ref{spine}).  We prove that the pictograph is a topological-conjugacy invariant (Theorem \ref{spine invariance}).  It is finer than the Branner-Hubbard tableau, the $\tau$-sequence, and the DeMarco-McMullen tree. 

We describe the abstract construction of these pictographs in Section 8, and we prove Theorem \ref{maintheorem1} in Section 11, providing the inductive arguments for counting topological conjugacy classes associated to each pictograph.  We treat the case of cubic polynomials first and in greatest detail, in Sections 9--10.  We describe how from the pictograph one can recover the tree code (introduced in \cite{DM:trees}), the Branner-Hubbard tableau of \cite{Branner:Hubbard:2}, and the Yoccoz $\tau$-function.   We use these arguments to show (\S 10.1) the existence of cubic polynomials with 
\begin{itemize}
\item	the same tableau but different trees;
\item	the same tree but different pictographs; and
\item	the same pictograph but not topologically conjugate.
\end{itemize}
The examples we provide are structurally stable and in the shift locus of $\cM_3$, though they can be easily produced for cubic polynomials with only one escaping critical point.  In \cite{Branner:cubics}, Branner showed that there is only one Fibonacci solenoid in the moduli space $\cM_3$ of cubic polynomials; that is, fixing the critical escape rate, the union of all cubic polynomials with the Fibonacci marked grid (defined in \cite{Branner:Hubbard:2}) forms a connected solenoid.   Using our combinatorial techniques, we give (\S 10.2) an example of a pictograph for which the tableau (and $\tau$-sequence) are critically recurrent and for which the corresponding locus in $\cM_3$ consists of two solenoidal components.

%%%%%%
%%%%%%

\bigskip\bigskip

\begin{center}\textsc{{\bf I. Local structure}}\end{center} 

\nopagebreak
\section{Local models and finite laminations}  \label{sec:laminations}

From \cite[\S 4]{DP:basins}, a {\em local model map} is a branched cover 
	$$(Y,\eta) \to (X,\omega)$$
between Riemann surfaces equipped with holomorphic 1-forms, of a particular topological type.  These branched covers model restrictions of a polynomial branched cover $f: \C\to\C$.   In this section, we introduce finite laminations and their branched covers, as combinatorial representations of the local model maps.  We prove:

\begin{theorem} \label{branched cover}
\begin{enumerate}
\item A local model surface $(X,\omega)$ is uniquely determined, up to isomorphism, by its associated lamination and the heights of its inner and outer annuli.
\item A branched covering $(Y,\eta) \to (X, \omega)$ is uniquely determined, up to post-composition by an isomorphism of $(X,\omega)$, by the data consisting of the lamination associated to $(Y,\eta)$, the heights of its inner and outer annuli of $(Y, \eta)$, and the degree.
\end{enumerate}
\end{theorem}

\subsection{Local models} \label{local model section} 
A {\em local model surface} is a pair $(X,\omega)$ consisting of a planar Riemann surface $X$ and holomorphic 1-form $\omega$ on $X$ obtained in the following manner.  Begin with a slit rectangle in the plane
	$$R = \{x+iy: 0 < x < 2\pi, \; h_{min} < y < h_{max}\} \setminus \Sigma,$$
where $\Sigma$ is a (possibly empty) finite union of vertical slits of the form 
	$$\Sigma_j = \{x+iy: x = x_j, h_{min} < y \leq c_0\}$$ 
for a distinguished value of $c_0 \in (h_{min}, h_{max})$.  The surface $X$ is a closure of $R$, obtained 
by identifying pairs of vertical edges via horizontal translations.  The 1-form $\omega$ is the defined by $dz$ in the coordinates on $R$.  The $y$-coordinate in the rectangular representation induces a {\em height function} $h: X \to \R$.  The {\em central leaf} $L_X$ of the local model surface is the level set $\{z: h(z) = c_0\}$ containing all (if any) zeros of $\omega$; these are the images of the topmost points of the slits.   The complement $X\setminus L_X$ is a disjoint union of the {\em outer annulus} given by the quotient of $\{ c_0 < y < h_{max}\}$ and finitely many {\em inner annuli} given by the quotient of $\{h_{min} < y < c_0\}$.  
For convenience, we often suppress mention of the 1-form and write simply $X$ for $(X,\omega)$. 

Note that the rectangular coordinates can be recovered from the pair $(X,\omega)$ via integration $\phi(z) = \int_{z_0}^z \omega$, and the height function is given by $h(z) = \Im \phi(z)$, up to the addition of a real constant.

A {\em local model map} (or simply a {\em local model}) is a branched cover between local model surfaces 
	$$f:  (Y,\eta) \to (X,\omega)$$
such that 
	$$\eta = \frac{1}{\deg f} \; f^*\omega$$
and $f$ sends the central leaf of $Y$ to the central leaf of $X$.  In \cite[Lemma 4.2]{DP:basins}, we observed that every local model arises as the restriction of a polynomial branched cover $f: (\C, \eta) \to (\C, \omega)$ for a meromorphic 1-form $\omega$ having purely imaginary residue at each of its poles.  

An {\em isomorphism} of local model surfaces is a degree 1 local model map.

\subsection{Finite laminations}
Let $C$ be an oriented Riemannian 1-manifold, isometric to $\R/2\pi\Z$ with the standard metric and affine structure.  A {\em finite lamination} is an equivalence relation $L$ on $C$ such that 
\begin{enumerate}
\item	each equivalence class is finite, 
\item	all but finitely many classes are trivial (consist only of one element), and 
\item	classes are {\em unlinked.}
\end{enumerate}
The third condition means that if two pairs of equivalent points $\{a,b\}$ and $\{c,d\}$ lie in distinct equivalence classes, then $a$ and $b$ are in the same connected component of $C\setminus\{c,d\}$.  We will deal exclusively with finite laminations, so we henceforth drop the adjective ``finite''.  More general types of laminations play a crucial role in the classification of the dynamics of polynomials; cf. \cite{thurston:iterated}, \cite{kiwi:char_laminations}.  

A lamination is conveniently represented by a planar {\em lamination diagram}, defined as follows.  Given a lamination on the circle $C$, choose an orientation-preserving, isometric identification of $C$ with $S^1 = \{z\in\C: |z|=1\}$.  
\begin{figure}
\includegraphics[width=1.5in]{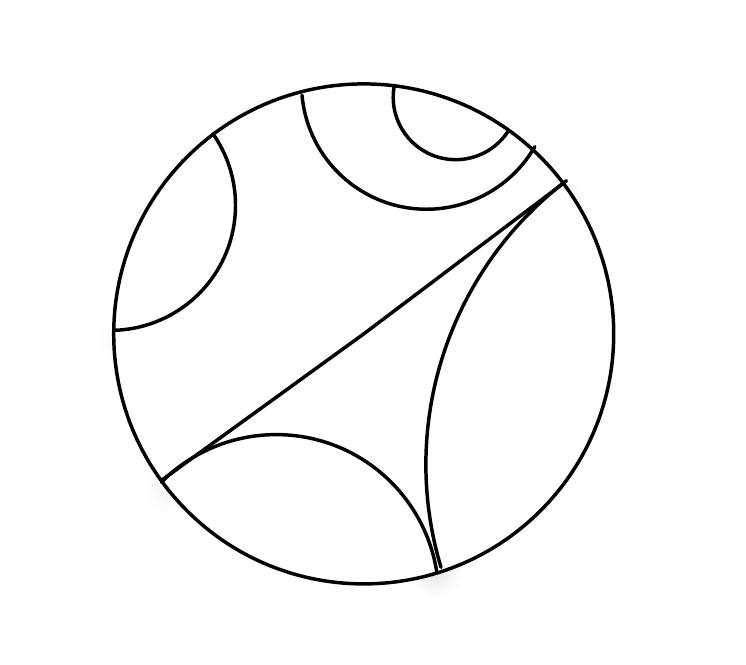}
\caption{A finite lamination with four non-trivial equivalence classes.}  \label{lamination}
\end{figure}
For each nontrivial equivalence class, join adjacent points in $C$ by the hyperbolic geodesic ending at those points, as in Figure \ref{lamination}.  Condition (3) that classes are unlinked guarantees that the hyperbolic geodesics do not intersect.

Two laminations are {\em equivalent} if there exists an orientation-preserving isometry (i.e. rotation) of their underlying circles taking one to the other.  Thus, there is no distinguished marking by angles on the circle $C$.

\subsection{Laminations and local model surfaces}  
Let $(X,\omega)$ be a local model surface, and let $L_X$ be its central leaf.  The 1-form $\omega$ induces an orientation and length-function on $L_X$, giving it the structure of the quotient of a circle by a finite lamination.  Therefore, there is a uniquely determined finite lamination associated to the local model surface $(X,\omega)$.

\begin{lemma} \label{lamination core}
A finite lamination $L$ determines a local model surface $(X,\omega)$, up to the heights of its inner and outer annuli.  
\end{lemma}

\proof
As we defined in \S\ref{local model section}, a local model surface is determined by its rectangular representation.  For any values $-\infty \leq h_{min} < c_0 < h_{max} \leq \infty$, we can construct a surface $(X,\omega)$ from the rectangle $\{0 < x < 2\pi, \; h_{min} < y < h_{max}\}$ with central leaf determining lamination $L$.  Indeed, choose any point on the circle $C$ to represent the edges $\{x=0 = 2\pi\}$.  For each point on $C$ in a non-trivial equivalence class, place a vertical slit from the bottom to height $y=c_0$.  Vertical edges leading to points in an equivalence class are paired by horizontal translation if they are joined by a hyperbolic geodesic in the diagram for $L$.  The unlinked condition in the definition of the finite lamination guarantees that $X$ is a planar Riemann surface.  The 1-form $dz$ on the rectangle glues up to define the form $\omega$.  It is immediate to see that the local model surface $(X,\omega)$ is determined up to isomorphism, once the values of $h_{min}$, $c_0$, and $h_{max}$ have been chosen.  %The surface is determined uniquely if we also specify the identification between its central leaf and the lamination $L$.
\qed

\subsection{Branched covers of laminations}
\begin{figure}
\includegraphics[width= 3in]{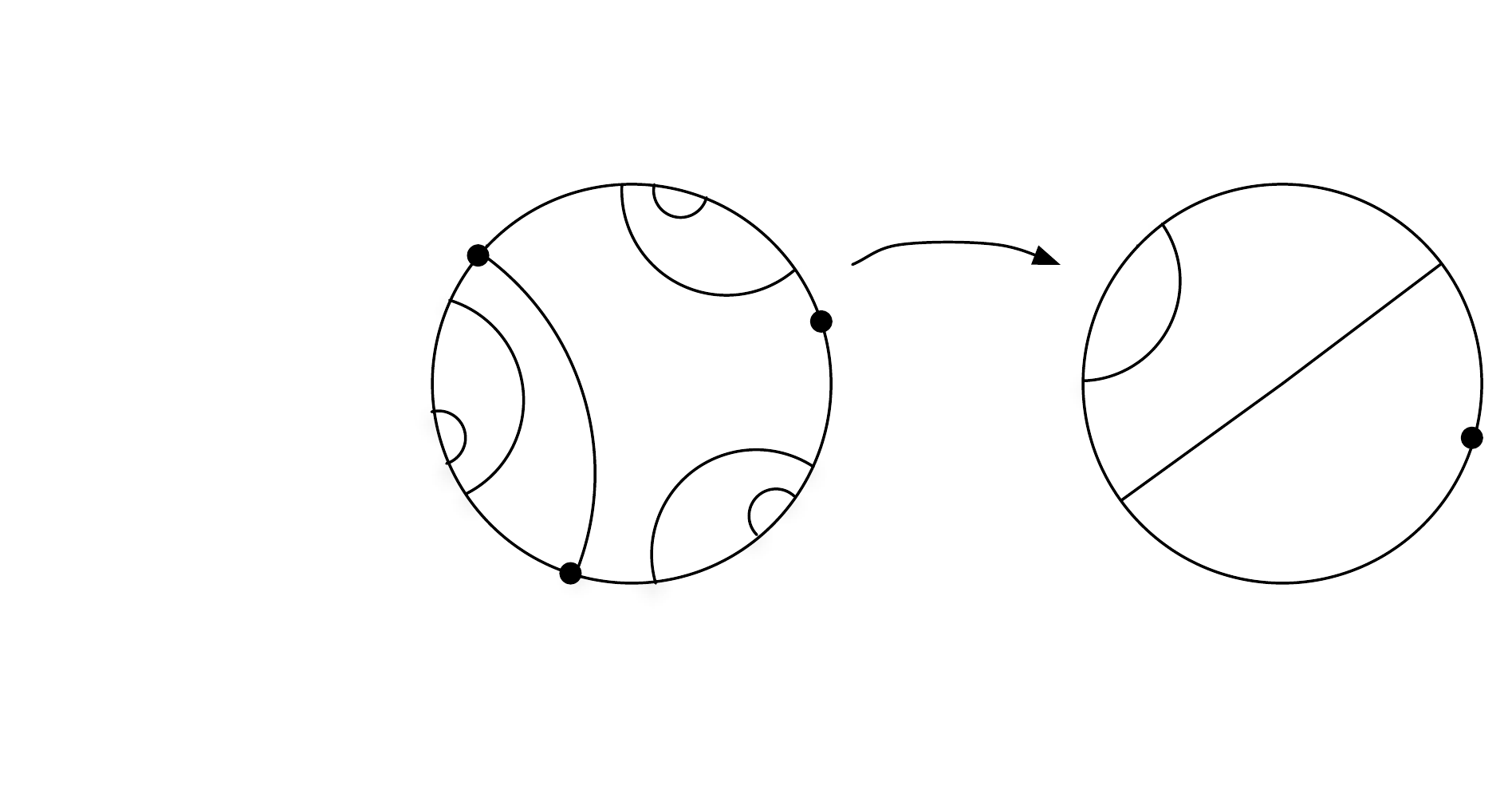}
\caption{A degree 3 branched cover of laminations.  The three marked points on the left are sent to the marked point on the right.} \label{lamination branched cover}
\end{figure}

If $L_1$ and $L_2$ are finite laminations, a {\em branched covering} $L_1 \to L_2$ is an orientation-preserving covering map $\delta: C_1 \to C_2$ on their underlying circles such that 
\begin{enumerate}
\item	$\delta$ is affine; i.e. $\delta(t) = ((\deg \delta) \, t + c) \mod 2\pi$ where each $C_i\iso \R/2\pi\Z$;
\item	for each equivalence class $A$ of $L_1$, the image $\delta(A)$ is equal to an (entire) equivalence class of $L_2$; and
\item	$\delta$ is {\em consecutive-preserving}.
\end{enumerate}
Consecutive-preserving means that for each equivalence class $A$ of $L_1$, either the image class $\delta(A)$ is trivial,  or consecutive points $x, y\in A$ (with respect to the cyclic ordering on $A$) are sent to consecutive points $\delta(x), \delta(y)$ in $\delta(A)$.  See Figure \ref{lamination branched cover}.

\begin{lemma} \label{equivalence}
A branched cover of laminations is determined by its domain and degree, up to symmetries.  More precisely, given branched covers $\delta: L_1 \to L_2$ and $\rho:M_1 \to M_2$ of the same degree, and given an isometry $s_1:  L_1 \to M_1$, there exists a unique isometry $s_2 : L_2\to M_2$ such that the diagram
$$\xymatrix{ L_1 \ar[r]^{s_1}    \ar[d]_\delta & M_1\ar[d]^{\rho} \\
			L_2 \ar[r]_{s_2}  &  M_2   }$$
commutes.  
\end{lemma}

\proof
The lamination diagram of $L_2$ is determined by $L_1$ and the degree; indeed, the rules for a branched covering determine the equivalence classes of $L_2$, as the images of those of $L_1$.  By hypothesis, there exists an isometry $s_1: L_1 \to M_1$ taking equivalence classes to equivalence classes, and $\delta$ and $\rho$ have the same degree.  Therefore, there exists an isometry $r_2: L_2 \to M_2$.  

Note that any branched cover $L_1 \to L_2$ of degree $d$ is determined by the image of a single point; this is because, in suitable coordinates, the covering is given by $t\mapsto d \, t$.  Now, suppose there exists an equivalence class $x$ in $L_1$ such that $r_2\circ \delta(x) \not= \rho\circ s_1(x)$ in $M_2$.  Combining the above facts, there is a uniquely determined symmetry $s: M_2\to M_2$ sending $r_2\circ\delta(x)$ to $\rho\circ s_1 (x)$.  We conclude that $s\circ r_2 \circ \delta = \rho\circ s_1$.  Set $s_2 = s\circ r_2$.
\qed

\medskip
In the proof of Theorem \ref{maintheorem1}, we will need to compute orders of rotation symmetry of certain lamination diagrams and record how these symmetry orders transform under branched covers.  

\begin{lemma}  \label{symmetry orders}
Let $\delta: L_1 \to L_2$ be a branched cover of laminations of degree $d$.  If $L_1$ has a rotational symmetry of order $k$, then $L_2$ has a rotational symmetry of order $k/\gcd(k,d)$.  
\end{lemma}

\proof
Suppose $s_1: L_1\to L_1$ is a symmetry of order $k$, so it rotates the circle underlying $L_1$ by $2\pi/k$.  By Lemma \ref{equivalence}, there exists a unique rotational symmetry $s_2: L_2\to L_2$ so the diagram 
$$\xymatrix{ L_1 \ar[r]^{s_1}    \ar[d]_\delta & L_1\ar[d]^{\delta} \\
			L_2 \ar[r]_{s_2}  &  L_2   }$$
commutes.  Fix coordinates on $L_1$ and $L_2$ so $\delta(t) = d\, t \mod 2\pi$.  As $s_1$ is a translation by $2\pi/k$, it follows that $s_2$ is a translation by $2\pi d/k \mod 2\pi$.  Therefore, $s_2$ is a symmetry of order $k/\gcd(k,d)$.  
\qed

\subsection{Gaps and local degrees}
In this work, a {\em gap} $G$ in a finite lamination is a connected component in the unit disk of the complement of the lamination, such that this component meets the boundary circle in a collection of arcs of positive length.  In other words, we consider the complementary components of the hyperbolic convex hulls of the equivalence classes of the lamination.  Or equivalently, a gap is a maximal open subset of the circle such that any pair of points in $G$ is unlinked with any pair of equivalent points.  

We remark that our terminology conflicts with that of \cite{thurston:iterated}; when an equivalence class consists of three or more points, we do not consider the ideal polygon it bounds in the disk as a gap.  The following lemmas are immediate from the definitions.

\begin{lemma} \label{closed gap}
A branched cover $\delta$ takes the closure of a gap surjectively to the closure of a gap.
\end{lemma}

Note that a gap itself is not necessarily mapped onto a gap by a branched cover, as seen by the example in Figure \ref{lamination branched cover}.  For the gap on the right containing the marked point, its preimage contains two gaps, one of which fails to map surjectively (as it misses the marked point).   

\begin{lemma} [and definition]  \label{gap degree}
The local degree of a gap $G$,
	$$\deg(\delta, G) = \frac{k|G|}{|\delta G|},$$
where $|G|$ is the length of $G$, is a positive integer which coincides with the topological degree of $\delta | G$. 
\end{lemma}

\begin{lemma}[and definition]  \label{eq class degree}
The {\em local degree} of $\delta$ at an equivalence class $A$,  
	$$\deg(\delta,A) = \frac{\#A}{\#\delta(A)},$$
is a positive integer and coincides with the topological degree of $\delta|A$.  
\end{lemma}

\subsection{Critical points of a lamination branched cover}  Suppose $\delta$ is a branched cover of laminations. 
An equivalence class $A$ is {\em critical} if $\deg(\delta, A) >1$, and a gap $G$ is {\em critical} if $\deg(\delta, G) >1$.  Abusing terminology, we refer to critical equivalence classes and critical gaps as {\em critical points} of $\delta$. 

\begin{lemma} \label{label count}
The total number of critical points of a branched cover $\delta$, computed by
	$$\sum_A (\deg(\delta, A) -1) + \sum_G (\deg(\delta, G) - 1),$$
is equal to $\deg\delta - 1$.
\end{lemma}

\proof
By collapsing equivalence classes to points, a lamination determines, and is determined by, a planar, tree-like 1-complex with a length metric of total length $2\pi$.  Tree-like means that it is the boundary of the unique unbounded component of its complement.  A branched covering determines, and is determined up to equivalence by, a locally isometric branched covering map between such complexes which extends to a planar branched covering in a neighborhood.  This lemma therefore follows from the usual Riemann-Hurwitz formula.
\qed

\subsection{Branched covers of laminations and local models}
Let $f: (Y,\eta) \to (X,\omega)$ be a local model map.  Let $L_X$ and $L_Y$ denote the finite laminations associated to the central leaves of $X$ and $Y$.  It is immediate to see that $f$ induces a branched cover of laminations $L_X\to L_Y$.

Conversely, we have:

\begin{lemma}  \label{loc}
A branched cover of finite laminations $\delta: L_1\to L_2$ determines a local model map, up to the heights of the inner and outer annuli of the local model surfaces.
\end{lemma}

\proof
Let $\delta: L_1\to L_2$ be a branched cover of laminations of degree $k$.  By Lemma \ref{lamination core}, we may construct, for each $i= 1, 2$, a local model surface $(X_i,\omega_i)$ so that its central leaf is identified with the lamination $L_i$; we may choose the heights of the inner and outer annuli to be $h_{min} = -\infty$ and $h_{max} = +\infty$ for $i=1,2$.  Because the height function is fixed, up to an additive constant, the surface $(X_i, \omega_i)$ is uniquely determined.  By choosing both $h_{min}$ and $h_{max}$ to be infinite, each inner annulus $(X_i, \omega_i)$ is isomorphic to the punctured disk $\{0 < |z| < 1\}$ equipped with the 1-form $ri \, dz/z$, where $r>0$ is the length of the corresponding gap in $L_i$, while the outer annulus of $(X_i, \omega_i)$ is isomorphic to the punctured disk $\{1 < |z| < \infty\}$ equipped with the 1-form $i \, dz/z$.

In these punctured-disk local coordinates, we extend $\delta$ by $z^k$, sending the outer annulus of $X_1$ to that of $X_2$.  For each gap $G$ of $L_1$, we extend $\delta$ by $z^{\deg(\delta,G)}$ in its punctured-disk coordinates. The local degree $\deg(\delta, G)$ is well-defined by Lemma \ref{gap degree}, and the extension is well-defined by Lemma \ref{closed gap}.  By construction, we obtain a branched cover $f: (X_1,\omega_1)\to (X_2, \omega_2)$ of degree $k$ such that $f^*\omega_2 = k \, \omega_1$ and $f$ induces the lamination branched cover $\delta: L_1\to L_2$.  

Note that finite choices of heights $h_{min}$ and $h_{max}$ give rise to local model maps that are restrictions of the constructed $f$.  In this case, there is a compatibility condition on the heights:  if an annulus has finite modulus $m$, then any degree $k$ cover is an annulus of modulus $m/k$.  If the domain surface $(X_1, \omega_1)$ has central leaf at height $h_0 \in \R$, and $h_{min} = h_{-1}$ and $h_{max} = h_1$, then the image surface $(X_2, \omega_2)$ will have outer annulus of height $k(h_1 - h_0)$ and inner annuli of height $k(h_0 - h_{-1})$.
\qed

\subsection{Proof of Theorem \ref{branched cover}}  The first conclusion is the content of Lemma \ref{equivalence}.  More functorially, an isometry between local model surfaces is determined by its restriction to the corresponding central leaves, so in particular the group of isometric symmetries of a local model surface is faithfully represented by the group of symmetries of its lamination.  More generally, since a branch cover $(Y, \eta) \to (X, \omega)$ in local Euclidean coordinates has differential a multiple of the identity, it is also determined by its restriction to the associated central leaves.  The second conclusion then follows from Lemma \ref{loc}.
\qed

%%%%%%
%%%%%%

\bigskip\bigskip

\begin{center}\textsc{{\bf II. The tree of local models}}\end{center} 

In Sections 3--5, we focus on the dynamics of a single polynomial $f: \C \to \C$ with $K(f)$ disconnected.  Naturally associated to $f$ will be two other dynamical systems, depending only on the basin dynamics $(f, X(f))$: the polynomial tree $(F, T(f))$, discussed in \S 3, and the {\em tree of local models}, in \S 4.  \S 5 is devoted to the important topic of symmetries of trees of local models.  

\nopagebreak

\section{Polynomial trees}
\label{sec:trees}

In this section, we recall some important definitions and facts from \cite{DM:trees} about polynomial trees.  We prove:

\begin{prop} \label{tree spine}
A polynomial tree $(F, T)$ is uniquely determined by the first-return map $(R, S_1(T))$ on a unit simplicial neighborhood of its spine.
\end{prop}

\noindent
The case of degree 3 polynomial trees is covered by \cite[Theorem 11.3]{DM:trees}.

\subsection{The metrized polynomial tree}  \label{tree}
Fix a polynomial $f: \C\to \C$ of degree $d\geq 2$.  Assume that at least one critical point of $f$ lies in the basin of infinity $X(f)$, so that its filled Julia set $K(f) = \C\setminus X(f)$ is not connected.  The escape-rate function is defined by
\begin{equation} \label{escape rate}
	G_f(z) = \lim_{n\to\infty} \frac{1}{d^n} \log^+ |f^n(z)|;
\end{equation}
it is positive and harmonic on the basin $X(f)$.  The {\em tree} $T(f)$ is the quotient of $X(f)$ obtained by collapsing each connected component of a level set of $G_f$ to a point.  

There is a unique locally-finite simplicial structure on $T(f)$ such that the vertices coincide with the grand orbits of the critical points.  The polynomial induces a simplicial branched covering
	$$F: T(f) \to T(f)$$
of degree $d$.  

The function $G_f$ descends to the {\em height function} 
	$$h_f : T(f) \to \R_+$$
which induces a metric on $T(f)$.  Adjacent vertices $v$ and $w$ have distance $|h_f(v) - h_f(w)|$.  

Let $E(f)$ be the set of edges in $T(f)$.  Note that the preimage in $X(f)$ of each edge in $T(f)$ is an annulus.  It maps by $f$ as a covering map to an image annulus.  The degree of these restrictions defines the {\em degree function} 
	$$\deg_f:  E(f) \to \N$$
of the tree $(F, T(f))$.  

\subsection{Fundamental edges and vertices}  \label{fundamental}
In this paragraph, we introduce some terminology and notation that will be employed throughout the paper. Let $v_0$ denote the highest vertex of valence $>2$, and let $v_1, v_2, \ldots, v_{N}:=F(v_0)$ be the consecutive vertices above $v_0$ in increasing height; we refer to $v_0, \ldots, v_{N-1}$ as the {\em fundamental vertices} and the edges $e_i$ joining $v_{i-1}$ and $v_i$, $i=1, \ldots, N$, as {\em fundamental edges}.  They will play an important role later.   Note that the union of the fundamental edges is a fundamental domain for the action of $F$.  

\subsection{The Julia set and weights} 
The {\em Julia set} $J(F)$ of the tree $(F, T(f))$ is the set of ends at height 0; 
we let $\overline{T}(f)=T(f)\cup J(F)$.  The quotient map $X(f) \to T(f)$ extends continuously to $\C\to \overline{T}(f)$, collapsing each connected component of $K(f)$ to a point of $J(F)$.  In \cite{DM:trees}, a probability measure on the Julia set of $F$ is constructed which coincides with the pushforward of the measure of maximal entropy for a polynomial under the natural projection $K(f) \to J(F)$.  

As in \S\ref{fundamental}, let $v_0$ be the highest branching vertex in $T(f)$.  For any vertex $v$ below $v_0$, its {\em level} $l(v)$ is the least integer $l\geq 0$ so that $h(F^l(v)) \geq h(v_0)$; this implies that $F^{l(v)}=v_j$ is a fundamental vertex for some $j \in \{0, \ldots, N-1\}$.  Denote the Julia set below $v$ by $J_v(F)$.  The measure $\mu_F$ is constructed by setting 
\begin{equation} \label{weight}
	\mu_F(J_v(F)) = \frac{\deg(v)\deg(F(v))\cdots\deg(F^{l(v)-1}(v))}{d^{l(v)}}.
\end{equation}
We refer to this quantity $\mu_F(J_v(F))$ as the {\em weight} of the vertex $v$; it will be used in \S\ref{sec:tlm} in the construction of the tree of local models.   

The numerator in (\ref{weight}) admits the following interpretation which will be used later.  For any edge $e$ below $v_0$, let $A_e$ be the annulus in $X(f)$ over $e$.  If $e$ is the edge above and adjacent to vertex $v$, then the ratio of moduli
	$$\mod(A_e)/\mod(A_{e_j})$$ 
is the reciprocal of the numerator in (\ref{weight}), where $e_j$ is the unique fundamental edge in the orbit of $e$.  This ratio will be called the {\em relative modulus} of the annulus $A_e$.

\subsection{Example:  degree 2}
\begin{figure}
\includegraphics[width=1.2in]{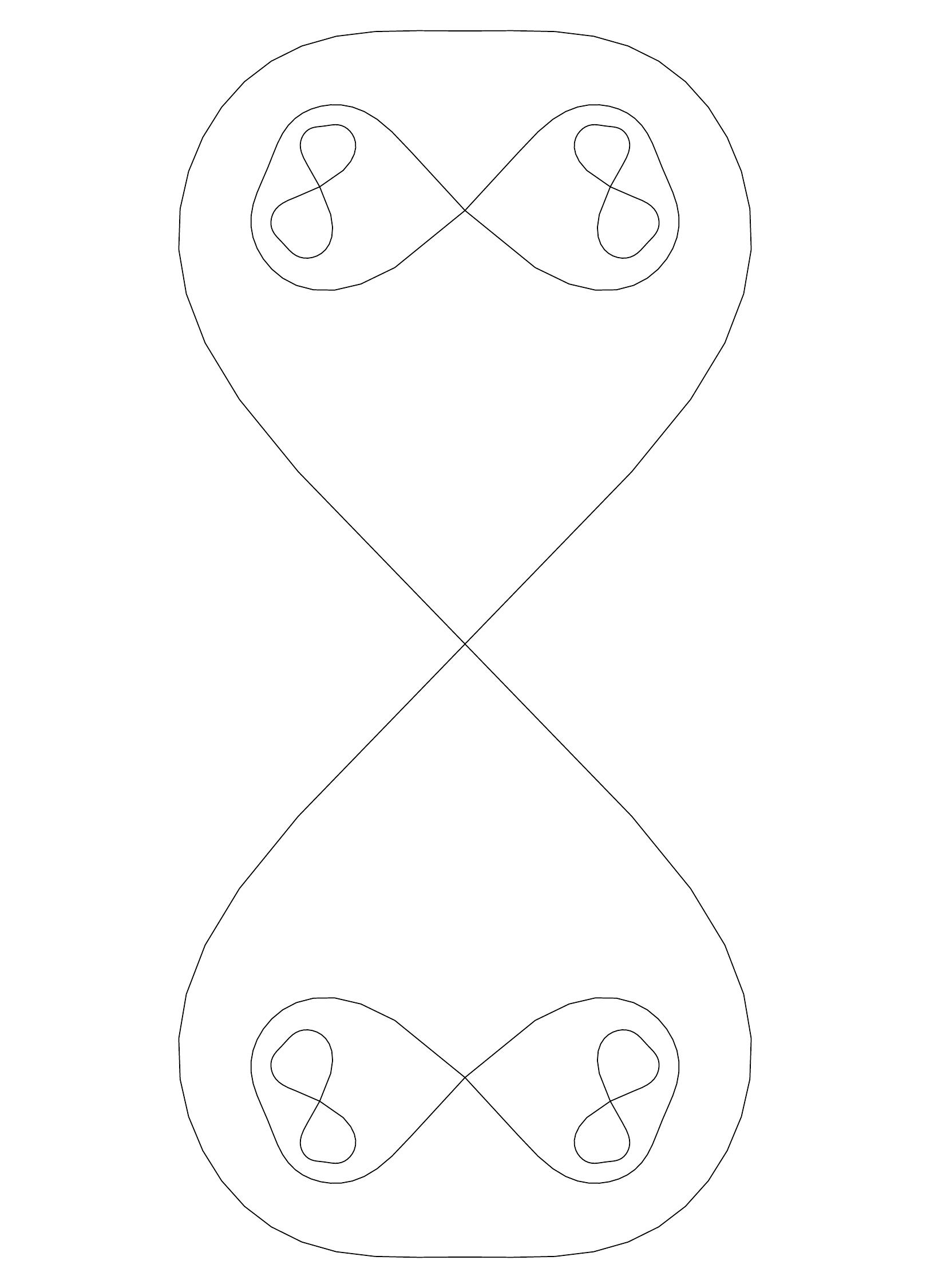} \qquad
\includegraphics[width=1.8in]{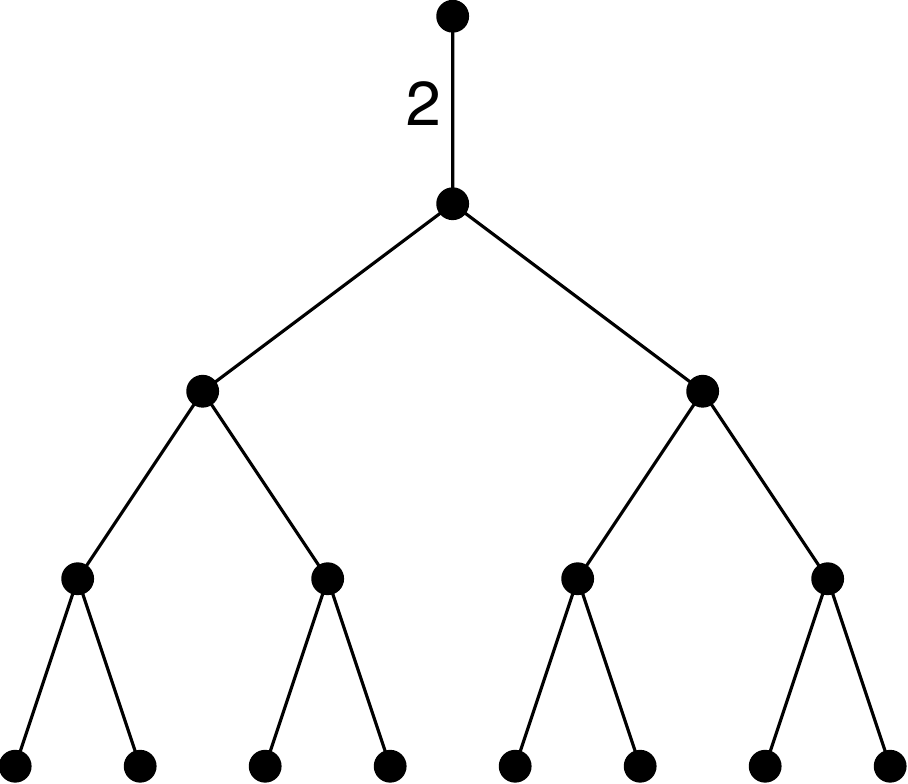} 
\caption{Critical level sets of $G$ and part of the tree associated to a quadratic polynomial with disconnected Julia set (from \cite[Figure 1]{DM:trees}).}
\label{quad tree}
\end{figure}

Trees in degree 2 are very simple to describe; up to scaling of the height metric, there is only one possibility.  Let $f_c(z) = z^2 + c$, and assume that $c$ is not in the Mandelbrot set, so the Julia set $J(f_c)$ is a Cantor set.  The level sets of the escape-rate function $G_c$ break the plane into a dyadic tree.  That is, for each $h> G_c(0)$, the level curve $\{G_c = h\}$ is a smooth topological circle, mapping by $f_c$ as a degree 2 covering to its image curve $\{G_c = 2h\}$; the level set $\{G_c = G_c(0)\}$ is a figure 8, with the crossing point at 0.  Each bounded complementary component of the figure 8 maps homeomorphically by $f_c$ to its image; there are thus copies of the figure 8 nested in each bounded component.  Consequently, level curves $\{G_c = G_c(0)/2^n\}$ are unions of figure 8s for all positive integers $n$; all other connected components of level curves in $X(f_c)$ are topological circles.   See Figure \ref{quad tree}.

The tree $T(f_c)$ has a unique highest branch point $v_0$, at the height $h(v_0) = G_c(0)$, and all vertices below zero have valence 3.  The action of $F: T(f_c) \to T(f_c)$ is uniquely determined, up to conjugacy, by the condition that $h(F(v)) = 2 h(v)$ for every vertex $v$ and that $F$ takes open sets to open sets.  Thus, the pair $(F, T(f_c))$ is completely determined by the height of the highest branch point, $G_c(0)$.

\subsection{The polynomial tree, abstractly defined}  \label{abstract tree}
In \cite{DM:trees}, it is established that these polynomial tree systems $(F,T(f))$ are characterized by a certain collection of axioms.  

We state the axioms here for convenience.  By a tree $T$, we mean a locally finite, connected, 1-dimensional simplicial complex without cycles.  Denote the set of edges by $E$ and the vertices by $V$.  For a given vertex $v\in V$, let $E_v$ denote the set of edges adjacent to $v$.  A simplicial map $F: T\to T$ is of {\em polynomial type} if 
\begin{enumerate}
\item	$T$ has no endpoints (vertices of valence 1);
\item	$T$ has a unique isolated end;
\item	$F$ is proper, open, and continuous;
\item	the grand orbit of any vertex includes a vertex of valence $\geq 3$, where $x,y\in T$ lie in the same {\em grand orbit} if $F^m(x) = F^n(y)$ for some positive integers $m, n$; and
\item	there exists a {\em local degree function} $\deg: E\cup V \to \N$ for $F$, satisfying 
	$$\deg(v) = \sum_{e\in E_v, \; F(e)=F(e')} \deg(e)$$
at every vertex $v$, for any given edge $e'$ adjacent to $v$,  and
\begin{equation} \label{critical point}
	2 \deg(v) - 2 \geq \sum_{e\in E_v} (\deg(e)-1) 
\end{equation}
at every vertex $v$.  
\end{enumerate}
It follows from the axioms that the topological degree of $F$ is well defined and finite, and it satisfies
	$$\deg F = \max_{v\in V} \; \deg v.$$
Further, it is proved in \cite{DM:trees} that the degree function for $(F, T)$, if it exists, is unique.  

A vertex $v$ is a {\em critical point} of $(F, T)$ if we have strict inequality in the relation (\ref{critical point}).  There are at most $d-1$ critical points, counted with multiplicity, and there is at least one critical point in the grand orbit of every vertex.

Theorem 7.1 of \cite{DM:trees} states that every pair $(F,T)$ of polynomial type is in fact the quotient of a polynomial of degree $\deg F$, with its simplicial structure uniquely determined by condition (4).  The critical vertices are the images of the critical points of the polynomial.  We sketch the proof of this realization theorem below in \S\ref{tree realization}.

Any tree of polynomial type $(F,T)$ can be endowed with a {\em height metric} $d_h$, which is linear on edges and the length of any edge $e$ satisfies $|e| = |F(e)|/ (\deg F)$.  There is a finite-dimensional space of possible height metrics, and each induces a {\em height function} 
	$$h: T\to \R_+$$ 
where $h(x)$ is the distance from $x$ to set of non-isolated ends (the {\em Julia set} of $F$).  The distance function can be recovered from $h$ by $d_h(v,w) = |h(v) - h(w)|$ on adjacent vertices.  When equipped with a height metric, we refer to the triple $(F,T,h)$ as a {\em metrized polynomial tree}.  The realization theorem of \cite{DM:trees} states further that every metrized polynomial tree $(F,T,h)$ arises from a polynomial $f: \C\to\C$ of degree $\deg F$ with $h$ as the height function $h_f$ induced by $G_f$.  

\subsection{The spine of the tree}
The {\em spine} $S(T)$ of a tree $(F, T)$ is the convex hull of its critical points and critical ends.  In other words, it is the connected subtree consisting of all edges with degree $>1$.  For example, in degree 2, $S(T)$ includes the highest branching vertex and the ray leading to infinity.  

We let $S_1(T)$ denote a unit neighborhood of the spine:  it includes all vertices at combinatorial distance $\leq 1$ from $S(T)$.  Let $(R, S_1(T))$ denote the first-return map of $(F, T)$ on $S_1(T)$.  More precisely: for a vertex or edge $u \in S_1(T)$ let $r(u)=\min\{ i > 0 : F^i(u) \in S_1(T)\}$; we then set $R(u)=F^{r(u)}$.  The map $R$ is semi-continuous: if $e$ is an edge above $v$, then $r(e) = r(v)$. 

We now prove that the full tree $(F, T)$ is determined by the renormalization $(R, S_1(T))$.  The argument is by induction on descending height.  The spirit of this argument will be used again to establish Propositions \ref{tlm spine}, Proposition \ref{gluing the spine}, Lemma \ref{spine automorphisms}, and Proposition \ref{spine}. 

\medskip\noindent{\em Proof of Proposition \ref{tree spine}.}    It suffices to show that the tree $T$ and map $F$ can be reconstructed, since the local degree function is uniquely determined \cite[Theorem 2.9]{DM:trees}.

Let $v_0$ denote the highest branching vertex.  All edges above $v_0$ have degree $d$ and are contained in the spine.  Above $v_0$, $F$ acts by translation along the ray $[v_0, \infty)$.  

The {\em star} $T_v$ of a vertex $v$ is the union of the vertex and its adjacent edges.  The unit neighborhood $S_1(T)$ of the spine of $T$ includes the star of $v_0$; the action of $R$ on this star collapses all edges below $v_0$ to a single edge of degree $d$.  Thus, we know the tree $T$ and the action of $F$ on all vertices at combinatorial distance $\leq 1$ from $v_0$.  Now suppose we have reconstructed $(F,T)$ at all vertices with combinatorial distance $\leq n$ from $v_0$.  Assume that $v$ is a vertex at combinatorial distance $n$ which is {\em not} in the spine $S(T)$.  Then $v$ has degree 1.  Consequently, its star is a copy of the star of its image $F(v)$, and the map $F$ on the star is uniquely determined.  

Now suppose $v$ is a vertex in $S(T)$ at combinatorial distance $n$ from $v_0$.  Then its star is contained in $S_1(T)$.  Let $w =  F^k(v)$ be the first return of $v$ in $S(T)$.  Note that $w$ may not coincide with $R(v)$, but $R$ allows us to reconstruct the action of $F^k$ from the star of $v$ onto the star of $w$.  For all $0 < j < k$, the star of $F^j(v)$ is mapped with degree 1 to the star of $F^{j+1}(v)$.  Therefore, the action of $F$ on the star of $v$ is exactly the action of $F^k$ on the star of $v$, where the image star (around $w$) is identified with the star of $F(v)$.  In this way, we have extended our construction of $(F, T)$ to combinatorial distance $n+1$.  
\qed

\subsection{Remark:  cubic polynomials and the spine}
The tree dynamical system $(F, T)$, is {\em not} determined by the first-return map on the spine $S(T)$ alone.  In fact, in the case of cubic polynomial trees, the data of the first return to $S(T)$ is equivalent to the data of the Branner-Hubbard tableau (or the Yoccoz $\tau$-sequence) which records the return of the critical point to its critical nest.  See \cite{DM:trees} for examples of distinct cubic trees with the same tableau; the examples are also shown in \S\ref{sec:examples}.

\subsection{Realization of trees, a review}  \label{tree realization}  In this subsection, we give an overview of the realization theorem for trees.  This construction motivated the definition of trees of local models, introduced in \S\ref{sec:tlm}.   

Begin with a metrized polynomial tree $(F,T, h)$ in the shift locus, so its critical heights are all positive.  A polynomial in the shift locus with tree $(F,T, h)$ is constructed as follows:
\begin{enumerate}
\item  {\bf Inflate the vertices.} Inductively on descending height, choose a local model surface $(X_v,\omega_v)$ over each vertex $v$ of $T$, ``modelled on" the vertex $v$ in $T$.
\item  {\bf Local realization.} 	The induction for Step (1) is done by choosing, for each vertex $v$, a local model map $$(X_v, \omega_v) \to (X_{F(v)}, \omega_{F(v)})$$ ``modelled on" $F$ at $v$.  The condition (5) on local degrees in \S\ref{abstract tree} guarantees the existence of such a local model.  The result of this step is a collection of local model maps, indexed by the vertices of $T$;  the domain $X_v$ of each is equipped with a natural projection to the star $T_v$. 

\item	 {\bf Glue.}  Over each edge of $T$, say joining $v$ to $v'$, glue the outer annulus of $X_v$ with the corresponding inner annulus of $v'$ so that the local model maps extend holomorphically (we do this more formally in \S \ref{glue} below).  After gluing all edges, we obtain a holomorphic map $f$ from a rigid, planar Riemann surface $X$ to itself.
\end{enumerate}
By uniformization, $X$ lies in the Riemann sphere, and $f$ extends uniquely to a polynomial whose basin dynamics $(f, X(f))$ is isomorphic to $(f, X)$.  By construction, the metric tree dynamics $(F, T(f), h_f)$ of the polynomial $f$  is isomorphic to that of the given metric tree $(F, T, h)$.

 A general metrized polynomial tree $(F,T,h)$---one with critical points at height zero---is realized by a compactness and continuity argument:  approximate $(F,T,h)$ by trees in the shift locus, realize each approximate tree by a polynomial, and pass to a convergent subsequence.  By continuity of metrized trees $f\mapsto (F, T(f), h_f)$, the limiting polynomial will have tree $(F,T,h)$.

%%%%%%
%%%%%%

\bigskip\bigskip
\section{The tree of local models}  \label{sec:tlm}

In this section we introduce the tree of local models $(\cF, \cX(f))$, as an intermediate dynamical object between the tree dynamics $(F, T(f))$ and the basin dynamics $(f, X(f))$.  It is intermediate in the sense that the basin dynamics $(f, X(f))$ determines the tree of local models dynamics $(\cF, \cX(f))$ which in turn determines the tree dynamics $(F, T(f))$.  Perhaps counterintuitively, however, it is not intermediate in the sense that the natural semiconjugacies are arranged via  
\[ (\cF, \cX(f)) \to (f, X(f)) \to (F, T(f));\]
the point here is that the gluing quotient map $\cX(f) \to X(f)$ depends on $f$.

Trees of local models can also be defined and constructed abstractly, and we show:

\begin{theorem} \label{tlm realization}
Every abstract tree of local models $(\cF, \cX)$  arises from a polynomial basin $(f, X(f))$.
\end{theorem}

\noindent
The proof is similar to the realization of abstract polynomial trees in \cite{DM:trees}.

As with polynomial trees (Proposition \ref{tree spine}), we consider the spine $(\cR, \cS)$ of the tree of local models $(\cF, \cX)$ and prove:

\begin{prop}  \label{tlm spine}
A tree of local models $(\cF, \cX)$ is uniquely determined by its first-return map on its spine $(\cR, \cS)$.
\end{prop}

\begin{prop}  \label{gluing the spine}
A tree of local models $(\cF, \cX)$ and a gluing along  its  spine $(\cR, \cS)$ determines uniquely a basin dynamical system $(f, X(f))$. 
\end{prop}

\noindent
While the basin dynamics $(f, X(f))$ in Proposition \ref{gluing the spine} is unique up to conformal conjugacy, the polynomial $f$ is not.

\subsection{The tree of local models, defined abstractly}  \label{abstract tlm} 
Let $(F,T, h)$ be a metrized polynomial tree.  For a vertex $v$, let $T_v$ denote the star of $v$.  A {\em tree of local models} over $(F,T,h)$ is a collection of triples $\{(f_v, (X_v, \omega_v), \pi_v): v \in V\}$, indexed by the vertices $v$ of $T$, such that  for each vertex $v$, 
\begin{enumerate}
\item	 the pair $(X_v, \omega_v)$ is a local model surface which is ``modelled on" the star $T_v$.  Specifically, there exists a marking homeomorphism 
	$$\pi_v: T(X_v,\omega_v) \to T_v$$ 
from the quotient tree of $(X_v, \omega_v)$, obtained by collapsing the leaves of the horizontal foliation of $\omega_v$ to points, to the star $T_v$.    We require further that $\pi_v$ is an isometry from the induced metric $\mu_F(J_v(F))\, |\omega_v|$ on $T(X_v,\omega_v)$ to the height metric on $T_v$, where $\mu_F(J_v(F))$ is the weight of $v$, defined in equation (\ref{weight}).    
\item the map 
	$$f_v: (X_v,\omega_v) \to (X_{F(v)}, \omega_{F(v)})$$ 
is a local model map which is ``modelled on" $F$ at $v$.  Specifically, via the marking homeomorphisms $\pi_v$ and $\pi_{F(v)}$, the restriction $F: T_v \to T_{F(v)}$ is the quotient of $f_v$, and the local degree function on $T_v$ coincides with the local degree of $f_v$ on leaves. 
\end{enumerate}
By condition (1), the heights of the inner and outer annuli in $X_v$ are controlled by the metric on $(F, T, h)$.  By condition (2), the widths of these annuli are also controlled, and therefore the moduli are determined.  In fact, the moduli coincide with the lengths of edges of $(F,T,h)$ in the {\em modulus metric} of \cite{DM:trees}.  

The data of a tree of local models can be bundled together to define a holomorphic degree $d$  branched cover 
	$$\cF: \cX \to \cX$$
of a disconnected Riemann surface $\cX$ to itself.  The 1-forms $\omega_v$ define a conformal metric $|\omega_v|$ on each $X_v$ (with singularities at the zeros of $\omega_v$).  On each local model surface $X_v$, the holomorphic map $\cF$ scales this metric by the factor $\deg(v)$.

\subsection{Equivalence of trees of local models}  \label{tree equivalence}
A tree of local models $\cF_1: \cX_1 \to \cX_1$ is {\em equivalent} to the tree of local models $\cF_2: \cX_2 \to \cX_2$ if there exists a holomorphic isometry 
	$$i: \cX_1 \to \cX_2$$ 
which induces a conjugacy 
	$$\cF_2 \circ i = i \circ \cF_1$$
while respecting the underlying tree structure.  That is, the isometry $i$ projects, via the marking homeomorphisms, to an isometry of polynomial trees,
	$$i: T_1 \to T_2$$
which conjugates the tree dynamics of $F_1$ to that of $F_2$.  

In particular, an automorphism of a tree of local models $\cF: \cX\to\cX$ is a holomorphic isometry $\cX\to \cX$ which induces an isometry of the underlying tree $T\to T$ and commutes with $\cF$.

\subsection{The tree of local models associated to a polynomial}
Let $f$ be a polynomial of degree $d$ with disconnected Julia set.  Let $G_f$ be its escape-rate function, and let 
	$$\omega = 2i \, \del G_f.$$  
Then $\omega$ is a holomorphic 1-form on $X(f)$, and each level curve $\{z\in X(f): G_f(z) = c\}$ has length $2\pi$ in the metric $|\omega|$.  The functional equation $G_f(f(z)) = d\, G_f(z)$ implies that 
	$$\frac{1}{d} f^*\omega = \omega.$$
Form the metrized polynomial tree $(F,T(f), h_f)$ as above.  Consider the projection $\pi_f: X(f) \to T(f)$ from the basin of infinity to the tree.  For each vertex $v\in T(f)$, let $X_v$ be the preimage in $X(f)$ of the star $T_v$, and set  
	$$\omega_v = \frac{1}{\mu_F(J_v(F))} \; \omega = \frac{2i \, d^{l(v)}}{\deg(f^{l(v)}|X_v)} \; \del G_f,$$
where the weight $\mu_F(J_v(F))$ was defined in equation (\ref{weight}).  Then the pair $(X_v, \omega_v)$ forms a local model surface, where each horizontal leaf of $\omega_v$ is the fiber over a point in $T_v$, and the central leaf is the fiber over $v$.  The normalization of $\omega_v$ is chosen so that the central leaf has length $2\pi$ in the metric $|\omega_v|$.  The Riemann surface $\cX$ is the disjoint union of surfaces $X_v$, indexed by the vertices $v$ in $T(f)$.

The restriction of the polynomial $f|X_v$ defines a local model
	$$f_v: (X_v, \omega_v) \to (X_{F(v)}, \omega_{F(v)}).$$
Indeed, the level $l(v)$ satisfies $l(F(v)) = l(v)-1$ whenever $l(v)>0$, so
	$$\omega_v = \frac{d^{l(v)}}{\deg(f^{l(v)}|X_v)} \; \frac{1}{d} \, f^*\omega = \frac{d^{l(F(v))}}{\deg(f|X_v) \deg(f^{l(F(v))}|X_{F(v)})} \; f^*\omega = \frac{1}{\deg(f|X_v)} \; f^*\omega_{F(v)},$$
as required for a local model map.  Therefore, the data $\{f_v, (X_v, \omega_v), \pi_v\}$ with markings induced from the projection $\pi_f$ define a tree of local models $(\cF, \cX(f))$ over $(F, T(f), h_f)$.

\subsection{The spine of the tree of local models} 
Fix a tree of local models $(\cF, \cX)$, and let $(F, T)$ be the underlying polynomial tree.  The {\em spine} $\cS$ of $\cX$ is the subset of $\cX$ lying over the spine $S(T)$ of the tree.  Like the spine of the underlying tree, there is an associated renormalization, the first-return map $\cR: \cS \to \cS$; for each vertex $v\in S(T)$, it is defined by  $\cR|_{X_v} = \cF^{s(v)}: X_v \to X_{F^{s(v)}} $ where $s(v) = \min\{i>0 : F^i(v) \in S(T)\}$.  Unlike the first-return map we consider for $(F, T)$, we do not thicken the spine to a unit neighborhood.   Note now that $\cR: \cS \to \cS$ is a holomorphic dynamical system; it particular, it is continuous and, in the natural Euclidean coordinates from the 1-form, is a homothety with scaling factor $\deg(v)$ away from singular points. 

We now prove that $(\cF, \cX)$ is uniquely determined by the first-return map $(\cR, \cS)$.  The proof proceeds along exactly the same lines as the proof of Proposition \ref{tree spine}. 

\medskip\noindent
{\em Proof of Proposition \ref{tlm spine}.}
The first observation is that the underlying metrized-tree dynamics $(F, T, h)$ can be recovered from the first-return map $(\cR, \cS)$.  Indeed, the local model surface over any vertex $v\in S(T)$ collapses to the star of $v$ (and determines the metric, locally).  Thus, the spine $\cS$ of $\cX$ determines the unit combinatorial-neighborhood $S_1(T)$ and a return map $R_1: S_1(T) \to S_1(T)$.  Strictly speaking, $R_1$ is not the {\em first} return on $S_1(T)$, but rather, the first return from the spine to itself, together with the action on stars.  Applying the {\em proof} of Proposition \ref{tree spine}, we are able to recover the full tree dynamical system $(F, T)$.  

As in the proof of Proposition \ref{tree spine}, we reconstruct $(\cF, \cX)$ from $(\cR, \cS)$ inductively on descending height.  We begin with vertices $v$ in the spine and use $\cR$ to reconstruct $\cF$ on the local model surface over $v$.  All other vertices have degree 1, so the map $\cF$ and surface $\cX$ are uniquely determined.  
\qed

\subsection{The gluing quotient map} \label{glue}
Suppose $(f, X(f))$ is a basin dynamical system.  For each vertex $v \in T(f)$, there is an inclusion $X_v(f) \hookrightarrow X(f)$.  The totality of these inclusions define a canonical semiconjugacy 
$(\cF, \cX(f)) \to (f, X(f))$ between the dynamics $\cF$ on the tree of local models $\cX(f)$ induced by $f$ and that of $f$ on its basin $X(f)$. If $v$ and $v'$ are incident, with $v'$ above $v$, the inclusions $X_v(f) \hookrightarrow X(f)$ and $X_{v'}(f) \hookrightarrow X(f)$ have the property that the image of the outer annulus of $X_v(f)$ coincides with that of the inner annulus of $X_{v'}(f)$.  The composition of the first with the inverse of the second gives a conformal isomorphism between these annuli.  

We conclude that a polynomial determines (i) a {\em gluing}: a collection $$\iota^f=\{\iota^f_e\}_{e \in E}$$ of conformal isomorphisms $\iota^f_e$ from the outer annulus of $X_v(f)$ to the corresponding inner annulus of $v'$ commuting with $\cF(f)$, one for each edge $e$ of the tree $T(f)$, and (ii) a corresponding {\em gluing quotient map} 
$$g_f: \cX(f) \to X(f).$$  
Note that an isomorphism $(f, X(f)) \to (g, X(g))$ lifts under the gluing projections to an isomorphism $(\cF(f), \cX(f)) \to (\cF(g), \cX(g))$. 

Conversely, given an abstract tree of local models $(\cF, \cX)$, one may consider abstract gluings as well.  An (abstract) {\em gluing} is a collection $\iota=\{\iota_e\}_{e \in E}$ of conformal isomorphisms $\iota_e$ from the outer annulus of $X_v$ to the corresponding inner annulus of $v'$ commuting with $\cF$, where $v$ and $v'$ are joined by an edge $e$.  An abstract  gluing $\iota$ defines a gluing quotient map $\cX \to \cX/\iota = X^\iota$ to an abstract planar Riemann surface to which the dynamics of $\cF$ descends to yield a proper degree $d$ holomorphic self-map $f^\iota: X^\iota \to X^\iota$.  In this way, a gluing $\iota$ defines a holomorphic semiconjugacy $(\cF, \cX) \to (f^\iota, X^\iota)$.   

Recall the definition of the fundamental edges and vertices, from \S\ref{fundamental}.  The choice of gluings along the $N$ fundamental edges determines the gluings at all vertices above $v_0$.  As with the tree of local models, a gluing can also be reconstructed from its action on the spine of $(\cF, \cX)$.  The proof of Proposition \ref{gluing the spine} is the same inductive argument employed now twice before. 

\medskip\noindent
{\em Proof of Proposition \ref{gluing the spine}.}
Fix a tree of local models $(\cF, \cX)=\{(f_v, (X_v, \omega_v), \pi_v): v \in V\}$ over a metric polynomial tree $(F, T, h)$ and let $(\cR, \cS)$ be the first-return map to its spine. 

Suppose we are given the data consisting of the gluings $\iota_e, e \in S(T),$ along the spine. 
Note that the the gluings at all edges above $v_N$ are determined by those at the fundamental edges $e_1, \ldots, e_N$.  
We proceed inductively on descending height.  Let $n \geq 0$ and suppose $\iota_e$ is defined along all edges joining vertices with combinatorial distance $\leq n$ from the highest branching vertex $v_0$.  Let $v$ be a vertex at distance $n+1$, joined by edge $e$ above it to vertex $v'$.  If $v$ lies in the spine, then $\iota_e$ has already been defined.  If $v$ is not in the spine, then $\deg(f_v: X_v \to X_{F(v)})=1$ and $f_{v'}$ has degree 1 on the inner annulus corresponding to $e$.  Setting $\iota_e:=f_{v'}^{-1} \circ \iota_{F(e)} \circ f_v$ gives the unique extension of the gluing along $e$ so that the needed functional equation is satisfied.  

The previous paragraph shows that gluings along the spine determine gluings on the whole tree of local models, yielding a holomorphic degree $d$ branched covering map $f$ of an abstract planar Riemann surface $X$ to itself.   The proof of the realization theorem (Theorem \ref{tlm realization} below) shows that the abstract basin dynamics $(f, X)$ is holomorphically conjugate to that of some polynomial.  
\qed

\subsection{Realization theorem} We now prove Theorem \ref{tlm realization}.  It may be useful to review the proof sketch of the tree realization theorem, given in \S\ref{tree realization}.  The final step in the proof is a continuity argument; to make the continuity argument go through in the setting of trees of local models, we rely on some observations from \cite{DP:basins}, in particular the proof of Lemma 3.2 there.  

\medskip\noindent
{\em Proof of Theorem \ref{tlm realization}.}
Let $(\cF, \cX)=\{(f_v, (X_v, \omega_v), \pi_v): v \in V\}$ be a tree of local models over the metrized tree $(F, T, h)$. When the tree $(F,T,h)$ lies in the shift locus, so that all critical points have positive height, the realization goes through as for trees.  The first two steps of tree realization (as described above in \S\ref{tree realization}) are already achieved by the given data.  We glue the local models, appeal to the uniformization theorem and extendability of the dynamics on the glued surface to all of $\C$, and conclude the existence of a polynomial in the shift locus realizing the given tree of local models.  

Now suppose the  tree $(F,T,h)$ has critical points in its Julia set (i.e. of height 0).  By \cite[Theorem 5.7]{DM:trees}, we can approximate $(F,T, h)$ by a sequence of metrized trees $(F_k, T_k, h_k)$ so that $(F_k, T_k, h_k)$ is isometrically conjugate to $(F, T, h)$ at all heights $\geq 1/k$, and further, all critical points of $(F_k, T_k, h_k)$ have height $\geq 1/k$.  We may construct trees of local models $(\cF_k, \cX_k)$ over each $(F_k, T_k, h_k)$, so that when restricted to heights above $1/k$, the dynamics of $\cF_k$ is holomorphically conjugate to that of $\cF$.   

Choose arbitrarily a gluing $\iota$ for $(\cF, \cX)$.  For each $k$, via the identification from the above conjugacies, we transport the gluing $\iota$ to a partially defined gluing for $(\cF_k, \cX_k)$, and we choose an extension arbitrarily to obtain a gluing $\iota_k$ for $(\cF_k, \cX_k)$.
By the first paragraph, these determine polynomials $f_k$ which we may assume are monic and centered.  Each of the polynomials $f_k$ has the same maximal critical escape rate, so by passing to a subsequence, we may assume the $f_k$ converge locally uniformly on $\C$ to a polynomial $f$.  

As in the proof of \cite[Lemma 3.2]{DP:basins}, the local uniform convergence $f_k \to f$ on $\C$ implies that for any $t >0$ the restrictions to $\{t \leq G_k(z) \leq 1/t\}$ converge uniformly to $f$ on $\{t \leq G_f(z) \leq 1/t\}$ and the 1-forms $\omega_k = \del G_k$ converge on this subset to $\omega = \del G_f$; indeed, the escape-rate functions are harmonic where positive, so the uniform convergence implies the derivatives also converge.  We therefore conclude that the tree of local models associated to $f$ is isomorphic to $(\cF, \cX)$.  
\qed

%%%%%%
%%%%%%

\bigskip\bigskip\section{Symmetries in the tree of local models} \label{sec:symmetries}

A tree of local models may admit many nontrivial automorphisms.   The group of such symmetries, unsurprisingly, will play an important role in the problem of counting topological conjugacy classes.

\subsection{The automorphism group} \label{automorphisms}
Let $\Aut(\cF, \cX)$ be the conformal automorphism group of the tree of local models $(\cF, \cX)$, as defined in \S\ref{tree equivalence}.  While any basin of infinity in degree $d$ has only a finite group of automorphisms, which is necessarily a subgroup of the cyclic group of order $d-1$, the group $\Aut(\cF, \cX)$ can be large and complicated.   Consider the following examples.

For any degree 2 tree of local models, we have $\Aut(\cF, \cX) \iso \Z/2\Z$.  The unique nontrivial automorphism is generated by an order-two rotation of the local model surface containing the critical point.  It acts trivially on all local models above the critical point.  The action on all vertices below the critical point is uniquely determined by the dynamics, because every such vertex is mapped with degree 1 to its image.

By contrast, consider the tree of local models for the cubic polynomial $f(z) = z^2 + \eps z^3$ for any small $\eps$.  This polynomial has one fixed critical point and one escaping critical point.  While the basin $(f, X(f))$ has no nontrivial automorphisms, the tree of local models has an automorphism of infinite order, acting by a rotation of order 2 at the vertex of the spine which is the preimage of the escaping critical point.  If the escape rate of the critical point in $X(f)$ is $M$, the rotation of order 2 at height $M/3$ induces an order $2^n$ rotation at the vertex in the spine of height $M/3^n$.  The action on the local model surface at each vertex of local degree 1 is uniquely determined; similarly for the vertices at heights greater than $M$.  In fact, for this example, $\Aut(\cF, \cX)$ is isomorphic to the profinite group $\Z_2$, the $2$-adic integers under addition; this follows from Lemma \ref{spine automorphisms} below. 

\subsection{Local symmetry at a vertex} \label{local symmetry}
Denote by $\IS^1$ the quotient group of $\R$ by the subgroup $2\pi \Z$. 
The group of orientation-preserving isometries of a Euclidean circle of circumference $2\pi$ is then canonically isomorphic to $\IS^1$ via the map which measures the displacement between a point and its image.  

Let $(\cF, \cX)$ be a tree of local models with underlying tree $(F, T)$.  Let $v$ be any vertex of $T$.  The outer annulus of $X_v$ is metrically the product of an oriented Euclidean circle $C$ of circumference $2\pi$ with an interval.  Let $\Stab_v(\cF,\cX)$ denote the stabilizer of $v$ in $\Aut(\cF, \cX)$, i.e. all $\Phi \in \Aut(\cF, \cX)$ with $\Phi(X_v) = X_v$.  Any element of this stabilizer induces a conformal automorphism of $X_v$.  Because this automorphism must preserve the outer annulus of $X_v$, it is necessarily a rotation.  Consequently, there is a well-defined homomorphism 
	$$\Stab_v(\cF, \cX) \to \Aut(X_v, \omega_v) \hookrightarrow \mbox{\rm Isom}^+(C)=\IS^1.$$

\begin{lemma}
\label{lemma:stabilizer_finite} 
For every vertex $v$, the image of $\Stab_v(\cF, \cX)$ in $\IS^1$ is a finite cyclic group $\Z/k(v)\Z$.  
\end{lemma}

\proof
Because elements of $\Aut(\cF,\cX)$ must commute with the dynamics, the points of the critical grand orbits are permuted, preserving heights; every local model surface $(X_v, \omega_v)$ contains at least one and finitely many such points.  Therefore the image of $\Stab_v(\cF,\cX)$ in the group of rotations is finite.
\qed

\medskip\noindent
The order $k(v)$ is called the {\em local symmetry} of $(\cF, \cX)$ at vertex $v$.

\subsection{Profinite structure}
Fix a height $t>0$.  Consider the automorphism group, similarly defined, of the restriction $(\cF_t, \cX_t)$ of the dynamics of $(\cF, \cX)$ to the local models over vertices with height $\geq t$.   Restriction gives a map $\Aut(\cF, \cX) \to \Aut(\cF_t, \cX_t)$; denote its  image by $\Aut_{(\cF, \cX)}(\cF_t, \cX_t)$.  If $0<s<t$ then restriction gives a compatible natural surjection 
	$$\Aut_{(\cF, \cX)}(\cF_s, \cX_s) \to \Aut_{(\cF, \cX)}(\cF_t, \cX_t).$$  
The structure of $\Aut_{(\cF, \cX)}(\cF_t, \cX_t)$ for large positive values of $t$ is easy to compute.  Recall the definition of the fundamental vertices from \S\ref{fundamental}.

\begin{lemma}  \label{upper automorphisms}
Let $(\cF, \cX)$ be a tree of local models over $(F, T)$ with $N$ fundamental vertices $v_0, \ldots, v_{N-1}$.  Fix $j \in \{0, \ldots, N-1\}$ and let $t$ be the height of $v_j$ in $T(F)$.  Then
	$$\Aut_{(\cF, \cX)}(\cF_t, \cX_t) \iso \prod_{i=j}^{j+N-1} \Z/k(v_i)\Z$$
where $k(v_i)$ is the local symmetry of $(\cF, \cX)$ at $v_i$.  
\end{lemma}

\proof
The automorphism group $\Aut(\cF, \cX)$ stabilizes all vertices $v_j$ with $j\geq 0$, and the cyclic group $\Z/k(v_i)\Z$ is the stabilizer of $v_i$.  The action of any automorphism at vertex $v_i$ uniquely determines its action at all vertices in its forward orbit, by Lemma \ref{equivalence}.  The fundamental vertices are in distinct grand orbits, so the automorphism group is a direct product.  
\qed

\begin{lemma} \label{profinite}
For any tree of local models, $\Aut(\cF,\cX)$ is a profinite group, the limit of the collection of finite groups $\Aut_{(\cF, \cX)}(\cF_t, \cX_t)$. 
\end{lemma}

\proof
It remains to show that the groups $\Aut_{(\cF, \cX)}(\cF_t, \cX_t)$ are finite for each $t>0$.  The group $\Aut_{(\cF, \cX)}(\cF_t, \cX_t)$ is finite by Lemma \ref{upper automorphisms} for all $t$ large enough.  From Lemma \ref{equivalence}, the action of an automorphism $\phi\in\Aut(\cF, \cX)$ at any vertex $v$ determines uniquely its action at the image vertex $F(v)$.  The vertices of a given height must be permuted by an automorphism, so we may apply Lemma \ref{lemma:stabilizer_finite} to conclude that $\Aut(\cF_t, \cX_t)$ is finite for any $t>0$.  The restriction maps allow us to view the full automoprhism group $\Aut(\cF, \cX)$ as an inverse limit.  
\qed

\medskip

\subsection{The spine and its automorphism group}  \label{automorphism induction}
Now suppose $(\cF, \cX)$ is a tree of local models with first-return map $(\cR, \cS)$ on its spine. Since $(\cR, \cS)$ is again a dynamical system, it too has an automorphism group $\Aut(\cR, \cS)$ which is naturally a profinite group.  It follows that $\Aut(\cR, \cS)$ is inductively computable; the base case is covered by Lemma \ref{upper automorphisms} at height $t = h(v_0)$.  Furthermore, in the shift locus, the subtree of $\cS$ below $v_0$ is finite, and $\Aut(\cR, \cS)$ is a finite group which is inductively computable in finite time.  

The restriction of any automorphism $\phi\in\Aut(\cF, \cX)$ to the spine $\cS$ is an automorphism of $(\cR, \cS)$.  Indeed, $\phi$ preserves local degree, and the spine consists of all vertices with local degree $>1$.  The following lemma then implies that $\Aut(\cF, \cX)$ is inductively computable from the data of $(\cR, \cS)$.

\begin{lemma}  \label{spine automorphisms}
The map 
	$$\Aut(\cF, \cX) \to \Aut(\cR, \cS), $$ 
which sends an automorphism to its restriction to the spine, is an isomorphism.
\end{lemma}

\proof 
Suppose $\psi \in \Aut(\cR, \cS)$ is given.  We use the usual inductive argument to show $\psi=\phi|_\cS$ for a unique $\phi \in \Aut(\cF, \cX)$.  Define $\phi$ by $\phi=\psi$ on the local model surfaces at and above the vertex $v_0$.  For the induction step, suppose $\phi$ has been constructed at all vertices with combinatorial distance at most $n$ from $v_0$, commuting with $\cF$.  Let $v'$ be a vertex at combinatorial distance exactly $n$ from $v_0$ and suppose $v$ is just below $v'$.  If $v \in \cS$ we set $\phi_v=\psi_v$ on the surface $X_v$.  If $v \not\in \cS$ then by induction $\phi$ has already been defined on $v'$ and on $w = F(v)$.  Let $\hat{w}=\phi(w)$, $\hat{v}' = \phi(v')$, and denote the image of $v$ under $\phi$, yet to be defined, by $\hat{v}$. 

The restriction $\phi_{v'}$ to $X_{v'}$ uniquely determines $\hat{v}$, because an automorphism must preserve the tree structure.   In addition, $\phi$ commutes with the local model maps at each vertex, so $f_{\hat{v}'}\circ \phi_{v'} = \phi_{F(v')} \circ f_{v'} $, and we deduce that $\hat{w} = F(\hat{v})$.  Since $v \not\in\cS$, neither is $\hat{v}$, and the local model maps $f_v: X_v \to X_w$ and $f_{\hat{v}}: X_{\hat{v}} \to X_{\hat{w}}$ are isomorphisms.  So the automorphism $\phi$ must send $X_v$ to $X_{\hat{v}}$ via the composition $(f_{\hat{w}})^{-1} \circ \phi_w \circ f_v$; this composition defines the extension $\phi_v$.   In this way, we have extended $\phi$ uniquely from combinatorial distance $n$ to combinatorial distance $n+1$, completing the proof.  
\qed

%Lemma \ref{spine automorphisms} implies that the automorphism groups of the trees of local models associated to the quadratic and cubic polynomials mentioned above have the indicated form.  In the quadratic case, the spine consists just of the vertices at and above $v_0$, and the only choice is the order two local symmetry at $v_0$.  In the cubic case, the spine consists of local model surfaces lying over vertices $\ldots v_2, v_1, v_0, v_{-1}, v_{-2}, \ldots$ lying on a bi-infinite ray in $T$.  We have $k(v_0)=1$ and for each $i < 0$, the map $f_{v_i}$ has degree two, and indeed $k(v_i)=2^{|i|}$. 

\subsection{Symmetries in degree 3}
We will use the following lemma in our computations for cubic polynomials in \S\ref{sec:degree3}.  

\begin{lemma}
\label{lemma:cubic tlm symmetries}
Suppose $f$ is a cubic polynomial with critical points $c_1, c_2$ and $(\cF, \cX)$ is its tree of local models. 
\begin{enumerate}
\item If $c_1=c_2$, then $\kappa(v_0)=3, \kappa(v_1)=1$, and $\Aut(\cF, \cX)$ is cyclic of order $3$.
\item If the heights of $c_1, c_2$ are the same and $c_1\neq c_2$, then either 
\begin{enumerate}
\item $\kappa(v_0)=\kappa(v_1)=1$ and $\Aut(\cF, \cX)$ is trivial, or 
\item $\kappa(v_0)=\kappa(v_1)=2$ and $\Aut(\cF, \cX)$ is cyclic of order $2$.
\end{enumerate}
\item In all other cases, the order of local symmetry of each fundamental vertex is equal to $1$.
\end{enumerate}
\end{lemma}

Case 1 occurs when exactly when $f(z)=z^3+c$ lies outside the connectedness locus; case 2(b) when $f$ admits an automorphism.  

\proof In case 1, the number of fundamental vertices is $N=1$, the local model map $f_{v_0}$ has a deck group of order $3$, and its range $X_{v_1}$ has a distinguished point, the unique critical value.  Case 2 is similar.  If $f$ has an automorphism, then there are symmetries of order 2 at $v_0$ and its image commuting with $f_{v_0}$; thus $k(v_0) = k(v_1) = 2$.  If $f$ fails to have an automorphism but $c_1\not= c_2$, there are no symmetries at $v_0$ and consequently no symmetries at $v_1$, so $k(v_0) = k(v_1) = 1$.  The conclusions about the automorphism groups then follow immediately from Lemma \ref{spine automorphisms}. 

To prove the last statement, suppose that the two critical points have distinct heights.  Then the local model map from $X_{v_0}$ to its image is a degree 3 branched cover with a unique critical point  (of multiplicity 1) in the surface $X_{v_0}$.  Such a branched cover has no symmetries, so $k(v_0) = k(F(v_0)) = 1$.  Further, if the two critical points are in distinct foliated equivalence classes, then the local model surface $X_{v_1}$ has a unique marked point on its central leaf (its intersection with the orbit of critical point $c_2$) that must be preserved by any automorphism; therefore, the local symmetry at $v_1$ will be 1.  
\qed

%%%%%%
%%%%%%

\bigskip\bigskip

\begin{center}\textsc{{\bf III. In the moduli space}}\end{center} 

\nopagebreak

\section{Topological conjugacy}
\label{sec:conjugacy}

In this section, we remind the reader of the quasiconformal deformation theory of polynomials, following \cite{McS:QCIII}.  We show that the tree of local models is invariant under topological conjugacies that preserve critical escape rates.  In other words:

\begin{theorem} \label{tlm invariance}
The tree of local models is a twist-conjugacy invariant.
\end{theorem}

We recall the topology on $\cB_d$, the moduli space of basins $(f, X(f))$ introduced in \cite{DP:basins}, and we study the locus $\cB_d(\cF, \cX)\subset \cB_d$ of basins with a given tree of local models $(\cF, \cX)$.  Recalling (Proposition \ref{gluing the spine}) that a gluing of $(\cF, \cX)$ determines a basin dynamical system, we refer to elements of $\cB_d(\cF, \cX)$ as {\em gluing configurations}.  

 If $\kappa$ is the reciprocal of a positive integer, we denote by $\kappa \IS^1$ the quotient group $\R/2\pi \kappa \Z$ of $\IS^1$ by the group generated by a rotation of order $\kappa^{-1}$.  Carefully accounting for symmetries in $(\cF, \cX)$, we show:

\begin{theorem} 
\label{thm:bundle_of_gluing}
Let $(\cF,\cX)$ be a tree of local models with $N$ fundamental edges.  Given a basepoint $(f, X(f)) \in \cB_d(\cF, \cX)$, there is a  continuous projection 
	$$\cB_d(\cF, \cX) \to (\kappa\IS^1)^N$$ 
for some positive integer $\kappa^{-1}$, giving $\cB_d(\cF, \cX)$ the structure of a compact, locally trivial fiber bundle over an $N$-torus whose fibers are totally disconnected.  The twisting action is the lift of the holonomy induced by rotations in each coordinate, and the orbits form the leaves of a foliation of $\cB_d(\cF, \cX)$ by $N$-dimensional manifolds.  The leaves are in bijective correspondence with topological conjugacy classes within the space $\cB_d(\cF, \cX)$.  For $(\cF, \cX)$ in the shift locus, the fibers are finite.  
\end{theorem}

\noindent
Consequently, the problem of classifying basin dynamics up to topological conjugacy amounts to computing the monodromy action of twisting in the bundle $\cB_d(\cF, \cX)$.  Leading to the proof of Theorem \ref{maintheorem1}, we observe:

\begin{cor}
\label{cor:monodromy}
Under the hypotheses of Theorem \ref{thm:bundle_of_gluing}, let $\theta \in (\kappa \IS^1)^N$ be any point in the base torus.  Then the set of topological conjugacy classes in $\cB_d(\cF, \cX)$ is in bijective correspondence with the orbits of the monodromy action of $\Z^N = \pi_1((\kappa \IS^1)^N, \theta)$ on the fiber above the basepoint $\theta$.
\end{cor}

\subsection{Fundamental subannuli}  \label{subannuli}
Fix a polynomial representative $f: \C\to\C$ of its conjugacy class, and let $G_f$ be its escape-rate function. The {\em foliated equivalence class} of a point $z$ in the basin $X(f)$ is the closure of its grand orbit 
	$$\{w \in X(f): \exists \; n, m\in \Z, f^n(w) = f^m(z)\}$$ 
in $X(f)$.  Let $N$ be the number of distinct foliated equivalence classes containing critical points of $f$.   Note that $N=0$ if and only if the Julia set of $f$ is connected, if and only if the maximal critical escape rate 
	$$M(f) = \max\{G_f(c): f'(c) =0\}$$ 
is zero.  For $N>0$, these critical foliated equivalence classes subdivide the fundamental annulus 
	$$A(f) = \{z\in X(f): M(f) < G_f(z) < d\, M(f)\}$$ into $N$ {\em fundamental subannuli} $A_1, \ldots, A_N$ linearly ordered by increasing escape rate.   

The number $N$ coincides with the number of fundamental edges or vertices of the tree $(F, T(f))$, as defined in \S\ref{fundamental}.  For each $i = 1, \ldots, N$, the annulus $A_i$ lies over the fundamental edge $e_i$.

\subsection{Quasiconformal deformations of the basin}  \label{qc}
For each conformal conjugacy class of polynomial $f \in \cM_d$, there is a canonical space of marked quasiconformal deformations of $f$ supported on the basin of infinity.  The general theory, developed in \cite{McS:QCIII}, shows that this space admits the following description; see also \cite{DP:heights}.  The wring motion of \cite{Branner:Hubbard:1} is a special case.

One can define quasiconformal stretching and twisting deformations on each of the subannuli $A_j$ independently so that the resulting deformation of the basin $X(f)$ is continuous and well-defined and an isometry on each horizontal leaf.   We will parametrize the deformations of each subannulus by $t + is$ in the upper half-plane $\Hyp$, acting by the linear transformation 
	$$\left( \begin{array}{cc} 1 & t \\ 0 & s \end{array} \right)$$
on a rectangular representative of the annulus in $\R^2$, of width 1 and height equal to the modulus $\mod A_j$, with vertical edges identified.  Extending the deformation to the full basin of infinity by the dynamics of $f$, the deformation thus defines an analytic map 
\[ \Hyp^N \to \cM_d, \]
sending point $(i, i, \ldots, i)$ to $f$.  By construction, the twisting deformations (where $s=1$ in each factor) preserve escape rates, while stretching (with $t=0$ in each factor) preserves external angles.  Both stretching and twisting send horizontal leaves isometrically to horizontal leaves.

An important idea of \cite{McS:QCIII} in this context is that any two polynomial basins $(f, X(f))$ and $(g, X(g))$ which are topologically conjugate are in fact quasiconformally conjugate by a homeomorphism of the above type:  it has a decomposition into $N$ stretching and twisting factors, each factor determined by its effect on the $N$ fundamental subannuli.  Moreover, if 
the forward orbits of two critical points meet a common level set in the closure $\{M(f) \leq G_f \leq  d\cdot M(f)\}$ of the fundamental annulus, the arc length (angular difference) between these points is preserved under any topological conjugacy.   
 (See \S 5 of \cite{McS:QCIII}.)

\subsection{Normalization of twisting} \label{twist normalization}
For the proofs of Theorems \ref{conjugacy classes}, \ref{deg3twist}, and \ref{deg3infinite}, it will be convenient to work with the following normalization for the twisting action.   Fix $f\in \cM_d$ and consider the real analytic map 
	$$\Tw_f: \R^N \to \cM_d$$
which parametrizes the twisting deformations in the $N$ fundamental subannuli of $f$, sending the origin to $f$.   We normalize the parameterization $\Tw_f$ so that the basis vector
	$${\bf e}_j = (0, \ldots, 0, 1, 0, \ldots, 0)\in \R^N$$
induces a full twist in the $j$-th fundamental subannulus.  That is, if $\mod A_j$ is the modulus of the $j$-th subannulus of $f$, then $\Tw_f(t \, e_j)$ corresponds to the action of $1 + i \, t/\mod A_j \in\mathbb{H}$ in the coordinates described in \S\ref{qc}.

\subsection{Twisting and the tree of local models}
We now prove that the tree of local models is invariant under the twisting deformation.  More precisely, a twisting deformation induces, via restriction to central leaves and extension by isometries, an isomorphism between trees of local models.

\medskip\noindent
{\em Proof of Theorem \ref{tlm invariance}.}  
Fix a tree of local models associated to a polynomial $f \in \cM_d$ and suppose a twisting deformation conjugates $(f, X(f))$ to $(g,X(g))$ by a quasiconformal map $h$.  Then $h$ induces an isomorphism $H: (F, T(f)) \to (G, T(g))$, and so for each $v \in T(f)$ the restriction of $h$ gives an affine map of local model surfaces $h_v: (X_v(f), \omega_v) \to (X_{H(v)}(g), \omega_{H(v)})$.  Since $h_v$ is an isometry on the corresponding central leaves, it extends to an isometry $\phi_v: X_v(f) \to X_{H(v)}(g)$.   The dynamics of $f$ and of $g$ is locally a constant scaling, so $\phi = \{\phi_v\}$ yields an isomorphism $(\cF(f), \cX(f)) \to (\cF(g), \cX(g))$.  
\qed

\subsection{The space of basins $\cB_d$}  
We begin by recalling results from \cite{DP:basins}.  The set $\cB_d$ of conformal conjugacy classes of basins $(f, X(f))$ inherits a natural Gromov-Hausdorff topology: two basins $(f, X(f)), (g, X(g))$ are $\epsilon$-close if there is a relation $\Gamma$ on the product $\{\epsilon < G_f < 1/\epsilon\} \times \{\epsilon < G_g< 1/\epsilon\} $ which is $\epsilon$-close to the graph of an isometric conjugacy.  The natural projection $\pi: \cM_d \to \cB_d$ is continuous, proper, and monotone (i.e. has connected fibers).  Both spaces are naturally stratified by the number $N$ of fundamental subannuli and the projection respects this stratification.  While discontinuous on $\cM_d$, twisting is continuous on each stratum $\cB_d^N$, by \cite[Lemma 5.2]{DP:heights}.

\subsection{The bundle of gluing configurations}  \label{bundle}
Fix a tree of local models $(\cF, \cX)$.  Recall from \S\ref{glue} that an abstract  gluing $\iota$ defines a quotient map $\cX \to \cX/\iota = X^\iota$ to an abstract planar Riemann surface to which the dynamics of $\cF$ descends to yield a proper degree $d$ holomorphic self-map $f^\iota: X^\iota \to X^\iota$.  In this way, a gluing $\iota$ defines a holomorphic semiconjugacy $(\cF, \cX) \to (f^\iota, X^\iota)$.   The holomorphic conjugacy class of $(f^\iota, X^\iota)$ we call the associated {\em gluing configuration}.  Given an abstract tree of local models $(\cF, \cX)$, we let $\cB_d(\cF, \cX)\subset \cB_d$ be the collection of all gluing configurations. Theorem \ref{tlm realization} implies this is nonempty.   

 We begin with an identification of $\cB_d(\cF, \cX)$ as a set. 
 
The automorphism group $\Aut(\cF, \cX)$ acts naturally on the set of gluings as follows.  Given an automorphism $\Phi \in \Aut(\cF, \cX)$ and a gluing $\iota = \{\iota_e\}_{e \in E}$, the gluing $\Phi.\iota$ is the collection of isomorphisms $\{(\Phi.\iota)_e\}_{e \in E}$ defined as follows.  Suppose edge $e\in E$ joins $v$ to the vertex $v'$ above it; set $\hat{v}=\Phi^{-1}(v)$ and $\hat{e}=\Phi^{-1}(e)$ and define 
\[ (\Phi.\iota)_e := \Phi_{\hat{v}'} \circ \iota_{\hat{e}} \circ \Phi_v^{-1}.\]
Put another way: a gluing $\iota$ defines an equivalence relation $\sim_{\iota}$, which is a subset of $\cX \times \cX$; the gluing $\Phi.\iota$ corresponds to the equivalence relation given by $(\Phi \times \Phi)(\sim_{\iota}) \subset \cX \times \cX$. 

\begin{prop}
\label{thm:bgc_as_set}
The natural map $\iota \mapsto (f^\iota, X^\iota)$ descends to a bijection between $\Aut(\cF, \cX)$-orbits of gluings and gluing configurations. 
\end{prop}

\proof In one direction, an automorphism sending one gluing to another, by definition, descends to holomorphic map conjugating the two gluing configurations.  In the other, a conjugacy between two gluing configurations lifts to an automorphism sending the first corresponding gluing to the second.
\qed

\medskip
With respect to the topology on the space of basins $\cB_d$, we now prove that the set of all gluing configurations forms a compact fiber bundle over a torus.  

The main idea in the proof of Theorem \ref{thm:bundle_of_gluing} is the following. The torus forming the base of the bundle parameterizes the gluing choices along the fundamental edges.  However, there is no canonical identification of $(\cF, \cX)$ with $(\cF, \cX(f))$.  The ambiguity is an element of $\Aut(\cF, \cX)$.  Hence we will pass to a convenient quotient of this torus which erases this ambiguity.   Twisting deformations alter these gluing choices in a continuous way.  Once a gluing  has been chosen over the fundamental edges, the remaining choices for gluing may be selected in stages, inductively on descending height.  At each stage, the set of choices is finite, so the totality of such choices is naturally either finite or a Cantor set.

\medskip\noindent
{\em Proof of Theorem \ref{thm:bundle_of_gluing}.}
Let $v_0, \ldots, v_{N-1}$ denote the fundamental vertices of the underlying polynomial tree and $v_N = F(v_0)$.  Recall the definition of the local symmetry $k(v)$ of $(\cF, \cX)$ at a vertex $v$, given in \S\ref{local symmetry}.  Let 
	$$\kappa = \frac{1}{ \lcm\{ k(v_0), k(v_1), \ldots, k(v_N) \} } \; .$$
We define a continuous projection 
	$$\cB_d(\cF, \cX) \to (\kappa\IS^1)^N$$
which will define the fiber bundle structure.  

Fix any basepoint $(f, X(f))$ in $\cB_d(\cF,\cX)$.  Choose any leaf $\gamma$ of the vertical foliation of the basin $(f, X(f))$, so $\gamma$ is an external ray (and non-singular above the height of $v_0$).  Let $(g, X(g)) \in \cB_d(\cF,\cX)$.   This means there is an isomorphism $\phi: (\cF(f), \cX(f)) \to (\cF(g), \cX(g))$ which restricts to isomorphisms
	$$\phi_i:  X_i(f) \to X_i(g)$$
between the local model surfaces over vertices $v_0, v_1, É, v_N$ for $f$ and $g$; the local isomorphisms $\phi_i$ must send points in the critical orbits of $f$ to those of $g$.  The $\phi_i$ are canonical only up to pre-composition by the restriction of an element of $\Aut(\cF(f),\cX(f))$ to $X_i(f)$.   In particular, the chosen leaf $\gamma$ determines a collection of $k(v_i)$ vertical leaves in $X_i(g)$.

For the basin $(g, X(g))$, where the local model surfaces $X_i(g)$ have all been glued, there is a well-defined angle displacement between the distinguished vertical leaves in $X_i(g)$ and those of $X_{i-1}(g)$, as a value in the circle $\kappa\IS^1$.  The value of this angle displacement, for each $i = 1, \ldots, N$, defines the projection.  

Local triviality and continuity of the projection can be seen from the twisting action.  We proved in \cite[Lemma 5.2]{DP:basins} that the twisting action of $\R^N$ is well-defined on the stratum $\cB_d^N$ of polynomial basins with $N$ fundamental subannuli.  It is continuous and locally injective.  The space $\cB_d(\cF, \cX)$ is invariant under twisting by Theorem \ref{tlm invariance}.  
The definitions of the twisting action and of the projection imply that twisting by  $t$ in the $i$th fundamental subannulus translates the $i$th coordinate of the image under projection to the base by $t \bmod 2\pi \kappa $.  
It follows that twisting defines a local holonomy map between fibers in the bundle of gluing configurations and the space of gluing configurations is foliated by $N$-manifolds whose leaves are the orbits under the twisting action. 

Recall from \S\ref{qc} that two polynomial basins are topologically conjugate and have the same critical escape rates if and only if they are equivalent by a twisting deformation.  Thus, the topological conjugacy classes within the space of gluing configurations are easily seen to be in one-to-one correspondence with the twisting orbits, i.e. leaves. 

We now show that the fibers are totally disconnected.  Recall the gluing construction used in the proof of Theorem \ref{tlm realization}.   First, fix a point $b$ in the base torus of the projection.  This corresponds to choosing one from among finitely many choices of gluings  over the fundamental edges joining $v_0, \ldots, v_N$; this determines all gluings at vertices above $v_0$.  The collection of gluing choices is now made sequentially by descending height.  At the inductive stage, we have a vertex $v$ joined up to a vertex $v'$ along an edge $e$ of degree $k$; there are $k$ choices for the gluing isomorphism over $e$.  After a choice is made at every vertex in the tree, we obtain a holomorphic self-map $f: X\to X$ which is conformally conjugate to $f: X(f) \to X(f)$ for some polynomial $f$, by Theorem \ref{tlm realization}.   All basins in $\cB_d(\cF, \cX)$ over the basepoint $b$ are obtained in this way.  By the discreteness of gluing choices at each vertex and the definition of the Gromov-Hausdorff topology on $\cB_d$, for any fixed combinatorial distance $n$ from $v_0$, the set of gluing configurations which can be produced using the continuous choices corresponding to the basepoint $b$ and to a fixed set of choices at the finite set of vertices $v$ below and at distance at most $n$ from $v_0$ is an open set in $\cB_d(\cF, \cX)$.   In this way, we see that each gluing configuration over the basepoint $b$ is in its own connected component and the fibers are finite if $(\cF, \cX)$ lies in the shift locus.  

It remains to show that the bundle of gluing configurations is compact.  By properness of the critical escape rate map $f \mapsto \{G_f(c): f'(c)=0\}$ on the space of basins $\cB_d$ \cite{DP:heights}, the bundle must lie in a compact subset of $\cB_d$.   Let $(f_n, X(f_n))$ be any sequence in the bundle converging to a basin $(f, X(f))$.  Exactly as in the proof of Theorem \ref{tlm realization}, we may deduce that $(f, X(f))$ has the same tree of local models, and is therefore in the bundle of gluing configurations; see also \cite[Lemma 3.2]{DP:basins}.
\qed

\begin{lemma} \label{Cantor fiber}
Let $(\cF, \cX)$ be a tree of local models.  If the fibers in the bundle of gluing configurations have infinite cardinality, then they are homeomorphic to Cantor sets.
\end{lemma}

\proof
The fibers are compact and totally disconnected by the previous lemma.  By Brouwer's topological characterization of the Cantor set \cite[Thm. 2-97]{hocking:young:topology}, we need only show the fiber is perfect.   From the inductive construction of basins from the tree of local models, we see that the fiber of the bundle of gluing configurations has infinite cardinality if and only if there are conformally inequivalent gluing choices at an infinite collection of heights tending to 0.  By the definition of the Gromov-Hausdorff topology on the space of basins $\cB_d$, basins are close if they are ``almost" conformally conjugate above some small height $t>0$.  Consequently, any basin in the bundle of gluing configurations can be approximated by a sequence where a different gluing choice has been made at heights $\to 0$.
\qed

\subsection{The bundle of gluing configurations, in degree 2} 
We can give a complete picture of the bundle of gluing configurations in degree two.
Let $(\cF, \cX)$ be a tree of local models in degree 2.  In the notation of the proof of Theorem \ref{thm:bundle_of_gluing}, we have $N=1$, $k(v_0) = 2$, and $k(v_1) = 1$, so $\kappa = 1/2$.  Since every edge below $v_0$ has degree one, once the basepoint $b \in (1/2)\IS^1$ corresponding to the gluing along the fundamental edge $e$ joining $v_0$ and $v_1$ has been chosen, the remaining gluings are uniquely determined.  Hence the projection map $\cB_d(\cF, \cX) \to (1/2)\IS^1$ is $1$-to-$1$  and the bundle of gluing configurations $\cB_2(\cF, \cX) \subset \cB_2$ is homeomorphic to a circle.  In more familiar language: it is the image of an equipotential curve around the Mandelbrot set in the moduli space $\cM_2$ via the homeomorphism from the shift locus in $\cM_2$ to that of $\cB_2$.  A full loop around the Mandelbrot set corresponds to an external angle displacement running from $0$ to $2\pi/2$.  In fact, this is the same as the loop in Blanchard-Devaney-Keen inducing the generating automorphism of the shift \cite{Blanchard:Devaney:Keen}; the two lobes of the central leaf at $v_0$ are interchanged under the monodromy generator.  

\subsection{The bundle of gluing configurations in degree 3}  \label{degree 3 bundles}
In degree three, we can give a complete succinct picture of the bundle of gluing configurations in a few special cases.  The remaining ones are handled by Theorems \ref{deg3twist} and \ref{deg3infinite}. 

Suppose $f$ is a cubic polynomial with an automorphism and disconnected Julia set.  Then both critical points escape at the same rate, the automorphism has order $2$, and it interchanges the two critical points and their distinct critical values.  It is easy to see that there is a unique branched cover of laminations of degree $3$ with this symmetry.  It follows that, for a given critical escape rate, there is a unique tree of local models $(\cF, \cX)$ with this configuration.  By Lemma \ref{lemma:cubic tlm symmetries}, $k(v_0)=k(v_1)=2$, so $\kappa=1/2$.
Like in the quadratic case, a basin of infinity is uniquely determined by the gluing of the local models along the fundamental edge, because all edges below $v_0$ have local degree $1$.  But it now takes {\em two} turns around the base (= one full twist in the fundamental annulus) to return to a given basin, because the angle displacement between a critical point and its critical value is an invariant of conformal conjugacy.  Thus the projection $\cB_3(\cF, \cX) \to (1/2)\IS^1$ is $2$-to-$1$ and the bundle of gluing configurations $\cB_3(\cF, \cX)$ is homeomorphic to a circle.  
 
Suppose $f$ is a cubic polynomial where the two critical points coincide and escape to infinity, so it has a monic and centered representation as $f(z) = z^3 + c$, with $c$ not in the connectedness locus.  Let $(\cF, \cX)$ be its tree of local models.  By Lemma \ref{lemma:cubic tlm symmetries}, $k(v_0)=3$, $k(v_1)=1$, and $\Aut(\cF, \cX)$ is cyclic of order $3$. The bundle $\cB_3(\cF, \cX)$ therefore projects to the circle $(1/3)\IS^1$.  Again, since all edges below $v_0$ map by degree $1$, the gluing at the fundamental edge determines the basin. 
Going around this base circle of length $1/3$ forms a closed loop in $\cB_3(\cF, \cX)$, because a basin of infinity is uniquely determined by the gluing along the fundamental edge; the bundle is homeomorphic to a circle.  Note that lifting this path to the family $\{z^3 + c: c\in\C\}$ induces only a half-loop around the connectedness locus, since $z \mapsto -z$ conjugates $z^3+c$ to $z^3-c$.

Suppose $(\cF, \cX)$ is a tree of local models for a cubic polynomial in the shift locus with $N=2$ fundamental edges.  By Lemma \ref{lemma:cubic tlm symmetries}, there are no symmetries over $v_0, v_1, v_2$, so $\kappa = 1$.   The base of the fiber bundle is $\IS^1\times \IS^1$.   For polynomials in the shift locus, there exists a height $t>0$ such that all vertices below height $t$ have local degree 1, so all fibers of the fiber bundle must be finite.  In fact, the bundle of gluing configurations $\cB_3(\cF, \cX)$ is homeomorphic to a finite union of smooth 2-tori; compare \cite[Theorem 1.2]{DP:heights}.  The question of how many tori comprise this finite union is answered by Theorem \ref{deg3twist}; each torus corresponds to a distinct topological conjugacy class of polynomials.  

Finally, let $(f, X(f))$ be any other basin in the space $\cB_3$, so it has $N=1$ fundamental edge and there are are no symmetries at the fundamental vertices by Lemma \ref{lemma:cubic tlm symmetries}.  Therefore, $\kappa =1$ and the base of the fiber bundle is the circle $\IS^1$.  The fibers are necessarily finite if $f$ lies in the shift locus, but the fibers can be finite or infinite in the case where one critical point lies in the filled Julia set.  The number of connected components in the bundle $\cB_3(\cF, \cX)$ (and their topological structure) is given in Theorems \ref{deg3twist} and \ref{deg3infinite}.

%%%%%%
%%%%%%

\bigskip\bigskip

\begin{center}\textsc{{\bf IV. Combinatorics and algorithms}}\end{center} 

\nopagebreak

\section{The pictograph}
\label{sec:spine}

In this section, we define the pictograph $\cD(f)$ of a polynomial $f$, with degree $d\geq 2$ and disconnected Julia set. Using Theorem \ref{tlm invariance}, we will first show 

\begin{theorem} \label{spine invariance}
The pictograph is a topological-conjugacy invariant.
\end{theorem} 

Formally speaking, the pictograph is a static object; there is no map.  We next show nevertheless that 

\begin{prop} \label{spine}
The tree of local models $(\cF, \cX(f))$ is determined up to holomorphic conjugacy by the pictograph $\cD(f)$ and either 
\begin{itemize}
\item[(i)] the critical escape rates $\{G_f(c): f'(c)=0\}$, or equivalently, 
\item[(ii)] the lengths of the fundamental edges $e_1, \ldots, e_N$.
\end{itemize}
\end{prop}

\noindent
It follows that the pictograph encodes the symmetry group $\Aut(\cF, \cX(f))$ of the tree of local models.
Proposition \ref{spine} and Theorem \ref{thm:bundle_of_gluing} immediately imply Theorem \ref{maintheorem2}.

\subsection{Pictographs}  \label{dec spines}
Fix a polynomial $f$ of degree $d$ with disconnected Julia set and $N$ fundamental subannuli (defined in \S\ref{subannuli}).  Let $(\cF, \cX)$ be the tree of local models associated to $f$, and let $(F, T)=(F, T(f))$ be its simplicial polynomial tree.  Recall that the spine $S(T)$ of the tree is the convex hull of its critical points and critical ends.   The {\em pictograph} is a collection of lamination diagrams, one for each vertex in $S(T)$ at and below vertex $v_N$, each diagram labelled by its intersection with the critical orbits.  It is defined as follows.

Let as usual $v_0$ be the highest branching vertex of $T(f)$ and set $v_N = F(v_0)$.  We consider the vertices $v\in S(T)$ which are at and below the height of $v_N$.  For each such vertex $v$, record the lamination diagram for the central leaf of the local model surface $X_v(f)$.  We join lamination diagrams by an edge if the corresponding vertices are joined by an edge in $S(T)$.  This forms a spine of lamination diagrams.  

To define the labels, we first choose an indexing of the critical points $\{c_1, c_2, \ldots, c_{d-1}\}$ of $f$.  It is convenient to index them in nonincreasing height order, so $G_f(c_i)\geq G_f(c_j)$ if $i< j$.   Given a vertex $v$, we label the corresponding lamination diagram as follows. 
Given an index $i \in \{1, \ldots, d-1\}$ and an integer $k \geq 0$, consider the point $f^k(c_i)$ and how this point is located relative to $X_v(f)$. 
\begin{itemize}
\item  If $f^k(c_i)$ lies in one of the bounded complementary components of $X_v(f)$, we label the corresponding gap in the lamination diagram for $X_v(f)$ with the symbol $k_i$. 

\item  If $f^k(c_i)$ lands on the central leaf of $X_v(f)$, we label the corresponding equivalence class in the lamination diagram for $X_v(f)$ by the symbol $k_i$.  When indicated by a drawing, we label only one representative point in the equivalence class.

\item Otherwise, the label $k_i$ does not appear in the lamination diagram for $X_v(f)$; note that the point $f^k(c_i)$ lies neither in the outer annulus nor in an inner annulus of $X_f(f)$. 
\end{itemize}

Thus, the data in the pictograph is the same as the static data of the collection, for the above vertices $v$ in the spine, of the local model surfaces $X_v(f)$ labelled in the above fashion, with the map forgotten. 

Suppose now $f$ and $g$ are two polynomials.  We say $f$ and $g$ have {\em equivalent} pictographs if, after applying some permutation of the set of indices $i=1, \ldots, d-1$ for the critical points, there exists 
\begin{itemize}
\item a simplicial isomorphism $S(f) \to S(g)$ between the subtrees of the tree spines for $f$ and $g$, sending a vertex $v$ for $f$ to a vertex $\phi(v)$ for $g$, 
\item for each such vertex $v$, a corresponding isometry $\phi_v: X_v(f) \to X_{\phi(v)}(g)$ which sends each label $k_i$ appearing in the diagram for $X_v(f)$ to the same label $k_i$ appearing in the diagram for $X_{\phi(v)}(g)$. 
\end{itemize}

An abstract tree of local models $(\cF, \cX)$ determines a pictograph as well.   Counting both critical points in $\cX$ and those critical ends in $J(F)$ with multiplicity, there are again $d-1$ critical points.  Given a vertex $v$ with local model surface $(X_v, \omega_v)$, we regard a critical point $c_i$ below $v$ (in the tree) as``lying in'' the gap of $X_v$ corresponding to the edge leading to $c_i$.  

\subsection{Examples of pictographs}
\begin{figure}
\includegraphics[width=1.0in]{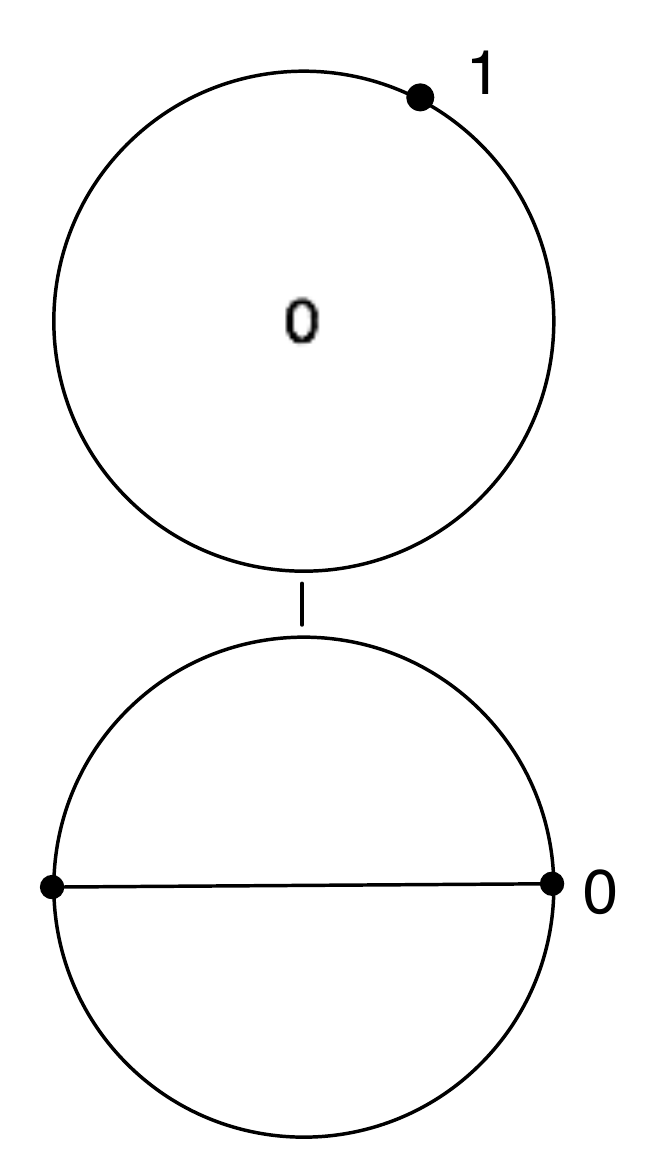}
\caption{The pictograph for every quadratic polynomial with disconnected Julia set.}
\label{quadratic spine}
\end{figure}

The degree 2 pictographs are the easiest to describe: in fact, there is only one possibility.  For quadratic polynomials $z^2+c$ with disconnected Julia set, the spine of the tree is the ray from the unique critical vertex $v_0$ heading to $\infty$.  The lamination diagram over the vertex $v_0$ is a circle cut by a diameter, representing the figure 8 level set $\{z\in\C: G_c(z) = G_c(0)\}$, with arclength measured by external angle.  The pictograph includes the data of this single diagram together with the image lamination (the trivial equivalence relation corresponding to level set $\{G_c = G_c(c)\}$), labelled by the symbol 0 to mark the critical point and 1 to mark the critical value.  See Figure \ref{quadratic spine}.  Because there is a unique critical point, we have dropped the subscript indexing.  Because angles are {\em not} marked on lamination diagrams, the pictographs are equivalent for all $c$ outside the Mandelbrot set.  (Recall that the tree of local models, and therefore the pictograph, is not defined for polynomials with connected Julia set.)  

\begin{figure} 
\begin{center}\includegraphics[width=1in]{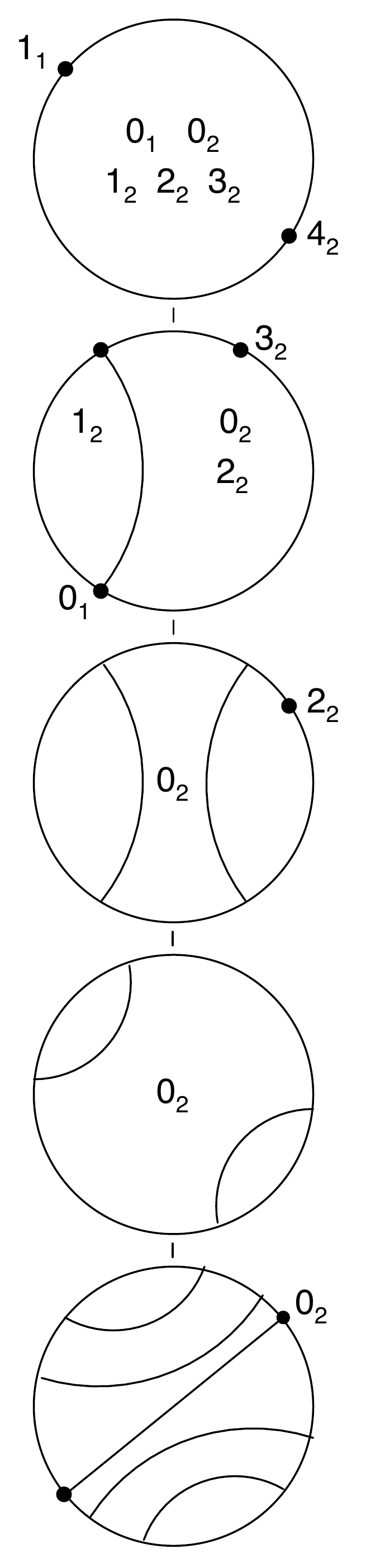}\end{center}
\caption{A cubic pictograph,  with critical escape rates $(M, M/3^3)$ for some $M>0$.}
\label{nongeneric cubic}
\end{figure}

For degree 3, Figure \ref{nongeneric cubic} shows a pictograph for a cubic polynomial with critical escape rates $G_f(c_1) = M$ and $G_f(c_2) = M/3^3$ for some $M>0$.  The spine of its tree is the linear subtree containing the four edges between critical point $c_2$ and critical value $f(c_1)$.  For the pictograph, we include five lamination diagrams at heights $M/3^i$, $i = -1, 0 , 1, 2, 3$.  The two critical points are labelled by $0_1$ and $0_2$.  Note that every spine in degree 3 will be a linear subtree of $T(f)$, because there are only two critical points.  

\begin{figure}
\includegraphics[width=2.5in]{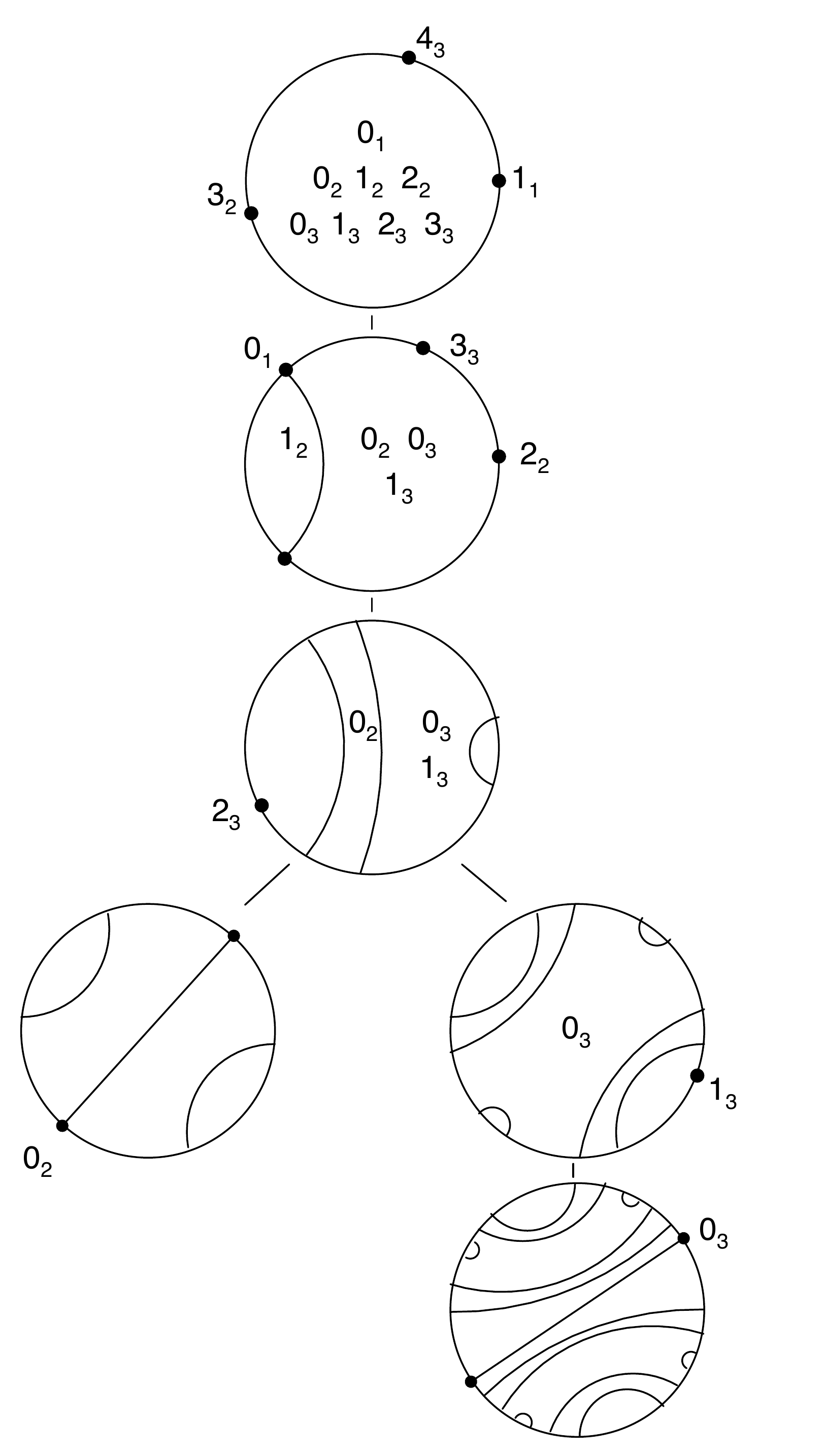}
\caption{A degree 4 pictograph, with critical escape rates $(M, M/4^2, M/4^3)$ for some $M>0$. } \label{degree 4 example}
\end{figure}

Figure \ref{degree 4 example} shows an example pictograph for a degree 4 polynomial with critical escape rates $G(c_1) = M$, $G(c_2) = M/4^2$, and $G(c_3) = M/4^3$ with a non-linear spine.  

In each of these examples, there is only one fundamental edge.  Figure \ref{0123} shows a cubic example with two fundamental edges.

\subsection{Proof of Theorem \ref{spine invariance}}

Suppose $f$ and $g$ are topologically conjugate.  Then there exists a quasiconformal conjugacy between basins $(f, X(f))$ and $(g, X(g))$.  By applying stretching deformations, we may assume the heights of the fundamental subannuli are the same, and that $(f, X(f))$ and $(g, X(g))$ are conjugate via a twisting deformation. By Theorem \ref{tlm invariance}, the trees of local models $(\cF_f, \cX(f))$ and $(\cF_g, \cX(g))$ are isomorphic via a holomorphic conjugacy $\phi$.   Choose arbitrarily an indexing of the critical orbits for $f$.  This indexing can be transported via $\phi$ to an indexing of those for $g$, so $f$ and $g$ will have equivalent pictographs.
\qed

\subsection{Reconstructing the tree of local models}  We can now prove that a pictograph plus the list of critical escape rates determines the full tree of local models over a metrized polynomial tree.  The strategy is the following.  The critical orbit labels allow us to first reconstruct the first-return map $(R, S(T))$ on the spine of the underlying tree $(F, T)$.  Then we use the lamination diagrams (and Theorem \ref{branched cover}) to reconstruct the local model maps and thus the first-return map $(\cR, \cS)$ on the tree of local models.  The heights of the local model surfaces and the metric on the underlying tree are determined by the critical heights.

\medskip\noindent
{\em Proof of Theorem \ref{spine}.}
Suppose we are given the pictograph $\cD$ for a polynomial $f$ of degree $d$ and the list of critical heights $h_1 \geq h_2 \geq \cdots \geq h_{d-1} \geq 0$.  By Theorem \ref{tlm spine}, it suffices to reconstruct the spine $\cS$ of the tree of local models and its first-return map.  We begin with the reconstruction of the first-return map $(R, S(T))$ on the spine of the underlying tree. 

Let $N$ be the number of independent critical heights; heights $h$ and $h'$ are independent if there is no integer $n$ such that $h = d^n h'$.   Denote by $v_0$ the vertex associated to the highest non-trivial lamination in $\cD$.  There are exactly $N$ trivial laminations above $v_0$ in the pictograph, each marked by points of the critical orbits.  Denote these vertices by $v_1, v_2, \ldots, v_N$, in ascending order.   The spine $S(T)$ is part of the data of the pictograph, after adjoining the ray from $v_0$ to $\infty$.  As usual, to reconstruct the action of $R$ we proceed inductively on descending height.  Above $v_0$, we have $R = F$, acting as translation by combinatorial distance $N$.  Each vertex of $S(T)$ below $v_0$ at combinatorial distance $j$ from $v_0$, with $0 < j \leq N$, is sent by $F$ to the vertex $v_{N-j}$.  

Now suppose we have computed the action of $R$ on $S(T)$ for all vertices at combinatorial distance $\leq n$ from $v_0$, and assume $n\geq N$.  Let $v$ be a vertex in $S(T)$ at combinatorial distance $n+1$ from $v_0$, and let $v'$ be the adjacent vertex above it.  Suppose $w' = R(v')$.  From the combinatorial distance between $w'$ and $v_0$, we can determine the iterate $k$ for which $w' = F^k(v')$.  Namely, if there are $n'$ edges on the path from $w'$ to $v_0$, then necessarily we have $n-n' = kN$ for some positive integer $k$.  Then $w' = F^k(v')$.  

Choose any index $j$ so that the symbol $0_j$ appears in the lamination diagram of $v$.  Then the symbol $0_j$ must also appear in the lamination diagram of $v'$, and the symbol $k_j$ must appear in the lamination diagram $L_{w'}$ of $w'$.  If $k_j$ lies in a gap of $L_{w'}$ together with a symbol $0_\ell$ for some index $\ell$, then necessarily $k_j$ must appear in the lamination diagram below $w'$ also containing $0_\ell$.  This vertex $v(\ell)$ is uniquely determined by $\ell$, and we may conclude that $R(v) = v(\ell)$.  

If $k_j$ lies in a gap of $L_{w'}$ containing no symbols of the form $0_\ell$, then $F^k(v)$ is not in the spine.  We must pass to a further iterate.  For each iterate $R^m$, define $k(m)$ by $R^m(v') = F^{k(m)}(v')$.   Choose the smallest positive integer $m$ so that the symbol $k(m)_j$ lies in a gap together with a symbol of the form $0_\ell$ for some index $\ell$ in the lamination over $R^m(v')$.  Such an integer always exists because some iterate of $R$ must send $v'$ to one of the vertices $\{v_1, \ldots, v_N\}$.  For this integer $m$, we choose the vertex $v(\ell)$ below $R^m(v')$ containing $0_\ell$, and we set $R(v) = v(\ell)$.  In this way, we reconstruct $(R, S(T))$ to all vertices at combinatorial distance $n+1$ from $v_0$, completing the induction argument.  

Our next step is to reconstruct the height function $h$ on the spine $S(T)$.  It suffices to determine the height of the fundamental vertices $v_i$, for $i = 0, 1, \ldots, N-1$, because of the relation $h(F(v)) = d\, h(v)$ on vertices.  We are given that $h(v_0) = h_1$, the height of the highest critical point.  For each $i>0$, the lamination diagram over $v_i$ must contain at least one marked point, labelled by the symbol $k_j$ for some positive integer $k$ and index $j \in \{2, \ldots, d-1\}$.  It follows that $h(v_i) = d^k \, h_j$.  

At this point, we observe that we could have taken our initial data to be the lengths of the $N$ fundamental edges, rather than the heights of the critical points.  Indeed, if $l_i$ is the length of fundamental edge $e_i$, then the height function $h$ is determined as follows.  Set $l = l_1 + \cdots + l_N$.  Then $h(v_0) = \sum_{k=1}^\infty l/d^k$, the distance from $v_0$ to the Julia set $J(F)$.  Then $h(v_i) = h(v_0) + l_1 + \cdots + l_i$ for each $i = 1, \cdots, N$.  The height of all other vertices is determined by the relation $h(F(v)) = d\, h(v)$.

We now apply Theorem \ref{branched cover}(1) to reconstruct the local model surfaces $X_v$ over each vertex $v$ in $S(T)$.  Setting the length of the central leaf to $2\pi$, the heights of the inner and outer annuli coincide with the length of the underlying edges of the trees, scaled by a certain factor $c_v>0$.  The factor $c_v$ is the reciprocal of the weight $\mu_F(J(F, v))$ defined in (\ref{weight}); the weight of $v$ is computable from the first-return map $(R, S(T))$ because all vertices with degree $>1$ are contained in the spine.  

By Theorem \ref{branched cover}(2), the local model maps over the vertices in $S(T)$ can be reconstructed from the lamination diagrams.  Indeed, the degree is obtained by counting the number of symbols of the form $0_j$ and adding 1. Recall, however, that we are able to so reconstruct a local model map only up to pre- and post-composition with rotational symmetries.  

Since such symmetries consist entirely of rotations and must preserve all labels, the only configurations of labelled lamination diagrams which are symmetric are those for which all labels lie in a central gap which is fixed by this rotational symmetry. 
So suppose $R(v)=w=F^m(w)$, and consider the labelled lamination diagrams for $X_v$ and $X_w$.  We now consider several cases.
\begin{enumerate}
\item Suppose neither $X_v$ nor $X_w$ admit label-preserving symmetries.  Then there is a unique map $X_v \to X_w$ sending a label $k_i$ in $X_v$ to the corresponding label $(k+m)_i$ in $X_w$. 

\item Suppose $X_v$ admits label-preserving symmetries but $X_w$ does not.   Then all labels for $X_v$ lie in a common central gap and the covering $X_v \to X_w$ is cyclic. Up to isomorphism there is a unique such covering; we choose a representative arbitrarily.   Note that any two such choices differ by precomposition by an isometry which is a symmetry of the labelled diagram for $X_v$, yielding a holomorphic conjugacy between the two different extensions of the dynamics to $X_v$ which is the identity on $X_w$. 

\item Suppose both $X_v$ and $X_w$ admit label-preserving symmetries.  Again, all labels for $X_v$ and for $X_w$ must lie in central gaps fixed by the rotational symmetries, and the covering $X_v \to X_w$ is cyclic.   By elementary covering space theory, given any fixed covering $\cF_v: X_v \to X_w$ and any rotation $\beta: X_w \to X_w$, there exists a lift $\alpha: X_v \to X_v$ of $\beta$ under $\cF_v$.  This lift $\alpha$ again yields a holomorphic conjugacy between the two different extensions $\cF_v$ and $\beta \circ \cF_v$ of the dynamics to $X_v$ which is the identity on $X_w$.   

Given any fixed covering $\cF_v: X_v \to X_w$ and any label-preserving symmetry $\alpha: X_v \to X_v$, rotation by $\alpha^{-1}: X_v \to X_v$ yields a holomorphic conjugacy between the two different extensions $\cF_v$ and $\cF_v \circ \alpha$ of the dynamics to $X_v$ which is the identity on $X_w$.   

The previous two paragraphs cover all sources of ambiguity in the extension of the dynamics. 
\end{enumerate}

Thus as the induction proceeds, we see that at the inductive stage, we make choices for the extension of the dynamics, but that different choices are holomorphically conjugate by a map which affects only the surface over which the extension is made.  

It follows that any two different collections of choices will yield holomorphically conjugate dynamics $\cR: \cS \to \cS$ on the spine of the tree of local models.

The dynamical system $(\cR, \cS)$ determines the full tree of local models by Theorem \ref{tlm spine}.
\qed

%%%%%%
%%%%%%
\bigskip\bigskip
\section{Algorithmic construction of an abstract pictograph}
\label{sec:construction}

The pictographs defined in the previous section satisfy a collection of simple combinatorial conditions.  In this section, we provide algorithmic rules for the construction of an abstract pictograph.  We prove:

\begin{theorem}  \label{realization}
Every abstract pictograph arises for a polynomial.
\end{theorem}

\begin{theorem}  \label{compatible heights}
For each abstract pictograph of degree $d$, the set of allowable critical heights $(h_1, \ldots, h_{d-1})\in \R^{d-1}$, with $h_1 \geq \cdots \geq h_{d-1}$, is homeomorphic to an open $N$-dimensional simplex, where $1 \leq N \leq d-1$ is the number of fundamental edges.  
\end{theorem}

\noindent
We will use Theorems \ref{realization} and \ref{compatible heights} to construct examples and classify topological conjugacy classes in the following sections.

\begin{figure}
\includegraphics[width=6.5in]{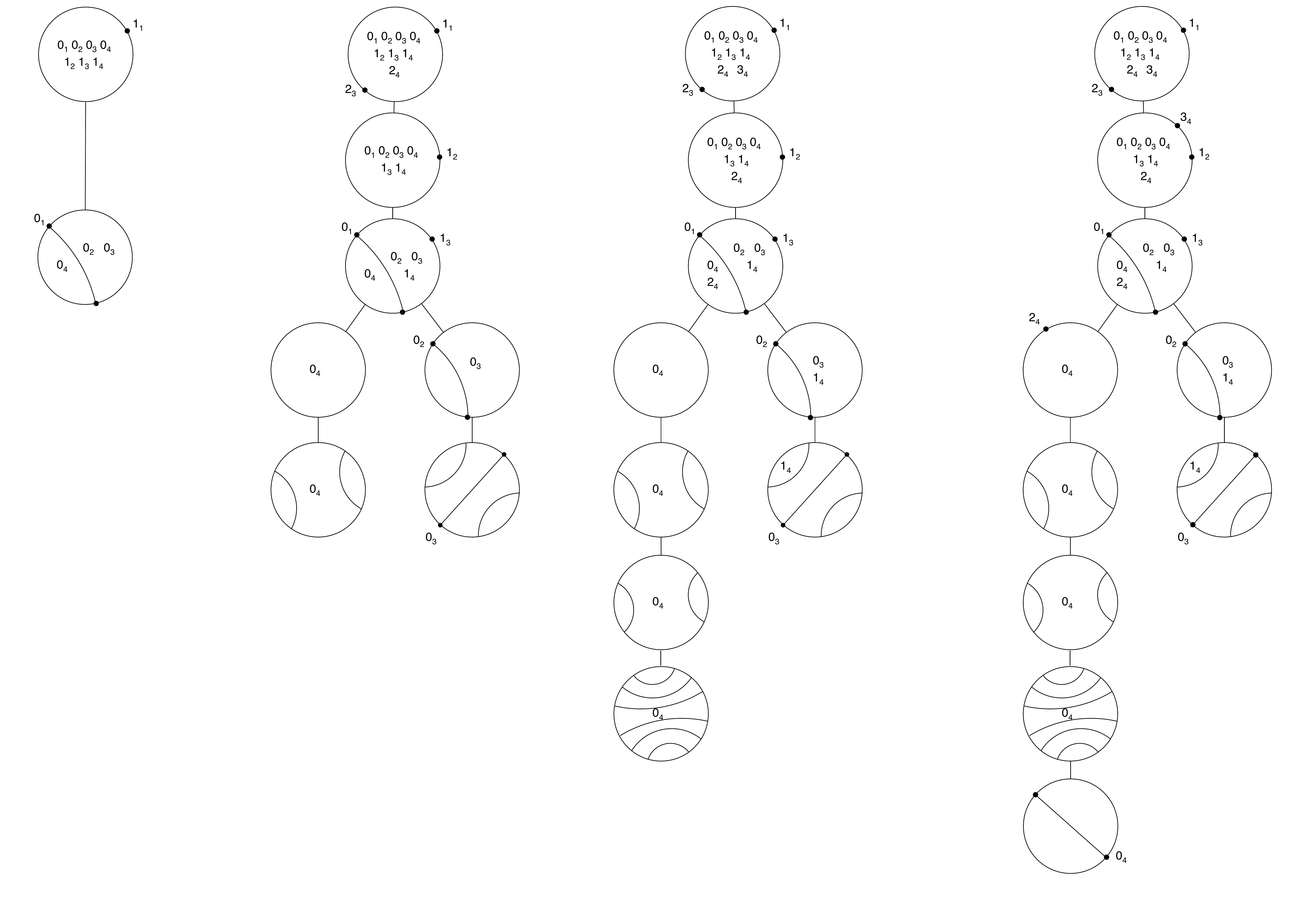}
\caption{From left to right, steps in the construction of an abstract pictograph of degree 5.}
\label{abstract spine}
\end{figure}

\subsection{Combinatorial pictograph construction}  \label{construction}
In this subsection, we present an algorithm for inductively constructing a pictograph.  The proof of Theorem \ref{spine}, which asserts that a pictograph determines uniquely a tree of local models, also applies here to show that the output of the algorithm does in fact yield a tree of local models.  This allows us to construct examples with desired properties.  

Though the spine is a static object, our algorithm produces branched coverings as it progresses, so that labels can be conveniently propagated forward under the dynamics. 

Fix a degree $d \geq 2$.  As the algorithm progresses, one makes choices.  The choices will influence the values of the following global variables:
\begin{itemize}
\item the number $0 < N < d$ of fundamental vertices $v_0, \ldots, v_{N-1}$
\item the geometry of the laminations at the fundamental vertices $v_i$
\item the presence or absence of specific critical orbit relations, including the multiplicity of critical points
\item the level of a critical point (defined as the number of iterations to reach a fundamental vertex), which may be infinite.
\end{itemize}

If desired, values of these global variables may be specified in advance; these values then constrain the possible choices as the algorithm progresses.   In particular, the geometry of the diagrams in the fundamental vertices can be specified ahead of time.

The algorithm terminates in finite time if and only if all critical orbits are eventually determined; e.g. the critical points have finite levels or lie in preperiodic gaps.  Otherwise, one must specify an infinite amount of data. 
 
 Upon working a few examples by hand, it is immediately evident that the set of choices at any given stage is rather large; the possibilities rapidly increase.    In particular, one of the steps of the algorithm amounts to the Hurwitz problem of choosing a planar covering, given the locations of branch values.  Thus use of the algorithm as a way to enumerate e.g. generic trees of local models in high degrees seems only theoretically and not practically possible.  
  
Before giving the algorithm, we introduce some terminology.
\medskip

{\bf Levels.}   Suppose $T$ is either the tree underlying a tree of local models, or a subtree thereof which contains the fundamental vertices.  Recall that a vertex $v$ at or below the fundamental vertices has {\em level $\ell$} if $\ell$ is the smallest nonnegative integer for which $F^l(v)$ is at or above $v_0$.  Our pictograph construction algorithm will start with a vertex $v_0$ at level $0$, its image, $v_1$ (thought of as lying above $v_0$) and will systematically add vertices in nondecreasing order of levels (that is, in nonincreasing order of its height).   
\medskip

{\bf Labels.}  Recall that a {\em label} is a symbol of the form  $k_i$ where $k\geq 0$ is an integer and $i \in \{1, \ldots, d-1\}$; it is a placeholder for the potential location of an iterate $f^k(c_i)$ of a critical point in the plane.  Given a label $k_i$, its {\em time} is the integer $k$, while its {\em critical point} is the integer $i$.  To {\em label} a lamination diagram is to assign a symbol of the form $k_i$ to an equivalence class (possibly trivial) or to a gap.  As we shall see, gaps will usually have many different labels.  Given such an equivalence class or gap, then, to say that it {\em contains the $i$th critical point at time $k$} means that it contains the label $k_i$.  

If a label $0_i$ ever appears on an equivalence class (rather than in a gap),  we say the orbit of $0_i$ is {\em determined}.  If a lowest vertex $v$ of the spine so far constructed has a gap with a label $0_i$, then the corresponding critical orbit is not determined, and we must continue to extend the spine downward from $v$ in order to complete the definition of the spine.  
\medskip

{\bf Propagation of labels.}  Suppose vertex $v$ has level $l\geq 0$, vertex $w$ has level $l-m \geq 1$, $m \geq 1$, and $S_v \to S_w$ is a branched cover of local model surfaces encoded by a branched covering map of laminations.  If $k_i$ is a label at $v$, then this label may be {\em propagated forward} via the map by labeling the image equivalence class or gap with the symbol $(k+m)_i$.  Given the branched cover of laminations, the image of a label under forward propagation is unique.

Suppose now vertex $v'$ is just above (in the tree) vertex $v$.   Since we are constructing the spine, we will have  $\deg(v)>1$, so there will be a distinguished critical gap of $v'$ corresponding to the edge joining $v$ up to $v'$.  If $k_i$ is a label in the diagram at $v$, then the gap of $v'$ corresponding to $v$ should have the label $k_i$ as well; we call this {\em upward propogation}.  

Suppose again that  vertex $v'$ is just above (in the tree) vertex $v$. If $k_i$ is a label in the gap of the lamination at $v'$ corresponding to $v$, then to {\em propagate $k_i$ downward} is to do one of the following:
\begin{itemize}
\item {\em bisect the edge} joining $v'$ and $v$ by adding a new vertex $v''$ between $v'$ and $v$, assigning it the trivial lamination, and labelling a point (trivial equivalence class) with the symbol $k_i$, or 
\item {\em drop the symbol $k_i$ down} by labeling an equivalence class or gap in the lamination  at $v$ with the symbol $k_i$.
\end{itemize}
Unlike forward and upward propagation, downward propagation will always involve choices.  
\medskip

{\bf Extension step.}  The following step incrementally grows the diagram for the pictograph.  See Figure \ref{abstract spine}.  Suppose 
\begin{itemize}
\item $k>0$, $l \leq 1$
\item $I \subset \{1, \ldots, d-1\}$ and  $J \subset \{1, \ldots, d-1\}$ are nonempty subsets
\item $v'$ is an existing vertex which is currently a lower end of the spine 
\item there is a gap $V$ in the lamination diagram of $v'$ whose label set consists exclusively of the symbols $0_i, i \in I$
\item $w'$ is an existing vertex which is not a lower end of the spine, and is at level $l$
\item the first-return map $\cR$ has already been defined at $v$, so that the data of the branched covering of labelled laminations encoding $\cR_{v'}: L_{v'} \to L_{w'}$ is already known
\item there is a gap $W$ in the lamination diagram for $w'$ containing the symbols $k_i, i \in I$ (and possibly other symbols as well) and the symbols $0_j, j \in J$; we assume that these labels $k_i$ all have a common time equal to $k$
\item there is a vertex $w$ just below $w'$, joined to $w'$ by an edge corresponding to the gap $W$, in whose lamination diagram the labels $k_i, i \in I$, have been placed.
\end{itemize}
The {\em extension step} is the creation of 
\begin{itemize}
\item a vertex $v$ below and joined to $v'$ by an edge corresponding to $V$
\item a lamination diagram $L_v$
\item a propagation of the labels $0_i, i \in I$, from $V$ down to the diagram for $v$
\item a choice of branched covering of laminations of degree $\#I + 1$ corresponding to a local model map $L_v \to L_w$ sending $0_i \mapsto k_i$; this is the extension of the return map over $v$.
\end{itemize}

Thinking of propagation and extension as subroutines, we now give the algorithm.  

\begin{enumerate}
\setcounter{enumi}{-1}
\item {\bf Initialization.} Fix a degree $d \geq 2$.  Set $N:=1$.   
Let $v_0, v_N=v_1$ be abstract vertices.  Define the set of fundamental vertices to be the set $\{v_0\}$; its level, by definition, is $0$.  
Join $v_0$ to $v_1$ by a vertical edge.

Choose any branched covering of laminations $\cR_{v_0}: L_{v_0} \to L_{v_1}$ of degree $d$ where the lamination for $L_{v_1}$ has no nontrivial equivalence classes, and $L_{v_0}$ has at least one nontrivial equivalence class.    

Draw the lamination diagram for $L_{v_1}$, and {\em below} it, draw the diagram for $L_{v_0}$. 
Join these diagrams by a vertical edge. 

Using the symbol set $\{0_1, 0_2, \ldots, 0_{d-1}\}$, label (arbitrarily) the critical equivalence classes and critical gaps of the lamination diagram $L_{v_0}$.  
Propagate these symbols forward via the above local model map.  To level $0$, we have now partially defined the pictograph $\cD$ and an associated first-return map $\cR$.   
\end{enumerate}

The algorithm loops the following steps.
\begin{enumerate}
\item {\bf Search for first return time.}   Let $v'$ be a vertex of the spine so far constructed which is a lower end; note that it might not be a vertex of minimal height.  Suppose there is a gap $V$ in the diagram for $v'$ which contains a collection of labels $0_i, i \in I$; these correspond to undetermined critical orbits.  Suppose $v'$ has level $l$.  There exists a unique triple $(k, w', W)$ consisting of a positive integer time $k$, a vertex $w'$, and a gap $W$ in the diagram of $w'$ such that 
\begin{enumerate}
\item $W$ contains some critical label $0_j$, 
\item for each $i \in I$, the label $k_i$ belongs to $W$,  
\item the label $k_i$ does not appear in the subtree below $w'$, 
\item the time $|k_i|=k$ is as small as possible with (a) and (b) above still holding.
%\marginremark{FIX:  should be lowest appearance of such a $k_i$, so change in level = time}
\end{enumerate}

If there is no nonempty subtree (of the spine already constructed) below $w'$, set $v':=w'$, and repeat the search for the first return time.

Otherwise, let $0_j$, $j \in J$, be the set of critical labels in $W$.

\item {\bf Propagate down.}   Propagate the labels $k_i, i \in I$, down the subtree strictly below $w'$.  This means that we are free to bisect edges, etc. as described above.  

\item {\bf Extend.} Inductively, one vertex at a time, extend the spine below the vertex $v'$ using the extension step; note that the first-return map is also extended.  

\item {\bf Propagate forward.}  The previous extension and downward propagation steps have created new labels in the tree; propagate them forward.  

\item {\bf Propagate upward.} If needed, propagate labels upward.

\end{enumerate}

\subsection{Proof of Theorems \ref{realization} and \ref{compatible heights}}
Let $\cD$ be a pictograph of degree $d\geq 2$, constructed according to the rules of \S\ref{construction}.  We first need to choose a compatible list of critical heights 
	$$\{h(c_1), h(c_2), \ldots, h(c_{d-1})\} \subset [0, \infty)$$
for the polynomial tree which will be determined by $\cD$.  For convenience, we reindex the critical points of $\cD$ so that we may have $h(c_1) \geq h(c_2) \geq \cdots \geq h(c_{d-1})$.  Fix any $M > 0$ to be the height of the highest critical vertex $v_0$ in $\cD$, and set $h(c_1) = M$.  Denote the vertices above $v_0$ in ascending order by $v_1, v_2, \ldots, v_N$ so the lamination of $v_N$ contains the symbol $1_1$ on its central leaf.  Choose any sequence of heights $M< h(v_1) < h(v_2) < \cdots < h(v_N) = dM$ for these vertices.  For each $i=2, \ldots, d-1$, there is at most one $j \in \{1, \ldots, N\}$ for which a symbol $k_i$ lies on the central leaf over $v_j$.  If such a $j$ exists, then we set
	$$h(c_i) = h(v_j)/d^k.$$
If no such $j$ exists, then we set $h(c_i) = 0$.  

The labelling rules in the spine construction guarantee that $h(c_1) \geq \cdots \geq h(c_{d-1})$, and $h(c_i) = 0$ if and only if the symbol $0_i$ appears in an infinite sequence of laminations descending from the vertex $v_0$.

As in the proof of Theorem \ref{spine}, the spine with the data of critical escape rates uniquely determines a full tree of local models.  Finally, we apply the tree-of-local-model realization Theorem \ref{tlm realization} to conclude that there exists a polynomial $f$ of degree $d$ with $\cD$ as its pictograph with its critical escape rates satisfying $G_f(c_i) = h(c_i)$.  
\qed

%%%%%%
%%%%%%

\bigskip\bigskip
\section{Combinatorics in degree 3}  \label{sec:degree3}

In this section, we provide the algorithm for computing $\Top(\cD)$, the number of topological conjugacy classes of basins $(f, X(f))$ with a given pictograph $\cD$, in degree 3.  The basin combinatorics of cubic polynomials are greatly simplified by the small number of critical points.   See for example the expansive treatment of cubic polynomials in \cite{Branner:Hubbard:1}, \cite{Branner:Hubbard:2} and the encoding of the cubic trees in \cite{DM:trees}.  We relate the combinatorics of pictographs with some results already known about the tableaux and turning curves of \cite{Branner:Hubbard:2}.  The proofs of the theorems below follow the arguments of Branner in \cite[Theorem 9.1]{Branner:cubics} and Harris in \cite{Harris:cubics}.  

To state the two theorems, we need the notion of a marked level.  Suppose $f$ is a cubic polynomial with a disconnected Julia set, so at least one critical point lies in the basin $X(f)$.  Suppose the critical points of $f$ are labelled as $c_1$ and $c_2$, with $G_f(c_2) \leq G_f(c_1)$.  For each integer $l > 0$, if $G_f(c_2) < G_f(c_1)/3^{l-1}$, there is a unique connected component $P_l$ of $\{z\in X(f): G_f(z) < G_f(c_1)/3^{l-1} \}$ containing $c_2$.  Let $B_l \subset P_l$ be the closed subset where $G_f(z) \leq G_f(c_1)/3^l$.  A {\em marked level} is an integer $l>0$ where the orbit of $c_2$ intersects $B_l\setminus P_{l+1}$.  The marked levels are called {\em semi-critical} in \cite{Milnor:lc} when $G_f(c_2) = 0$.

\begin{theorem} \label{deg3twist}
Suppose $\cD$ is a degree 3 pictograph with finitely many marked levels.  The number of topological conjugacy classes of basins $(f, X(f))$ with pictograph $\cD$ is 
		$$\Top(\cD) = \max_j \, \frac{2^j}{\max\{t_i: i\leq j\}},$$
where
\begin{enumerate}
\item the marked levels are indexed as $\{l_j\}_{j\geq 1}$;
\item for each $j$, we let $m_j$ be the sum of the relative moduli to level $l_j$; and
\item	 $t_j$ is the smallest positive integer for which $t_j m_j$ is integral.
\end{enumerate}
 If there are no marked levels, then $\Top(\cD) = 1$.  The computation of $\Top(\cD)$  depends only on the tableau of $\cD$. 
\end{theorem}

\noindent
The hypothesis of Theorem \ref{deg3twist} is clearly satisfied for all polynomials in the shift locus, because their pictographs have finite length. In that case, the conjugacy class of a basin $(f, X(f))$ determines the conjugacy class of the polynomial itself.    

In \cite{DS:count}, an explicit algorithm is developed and implemented using Theorem \ref{deg3twist} to enumerate all conjugacy classes of a given length (the greatest integer $l$ so that $G_f(c_2) < G_f(c_1)/3^{l-1}$) in the shift locus of cubic polynomials.  In particular, the algorithm includes an emumeration of all tableaux of a given finite length, followed by the enumeration of all possible pictographs (for structurally stable polynomials) associated to a given tableau.

For infinitely many marked levels, we find (in the notation of Theorem \ref{deg3twist}):

\begin{theorem}  \label{deg3infinite}
Suppose $\cD$ is a degree 3 pictograph with infinitely many marked levels.  Then there are infinitely many topological conjugacy classes of cubic polynomials with pictograph $\cD$.  Fixing the maximal critical escape rate $M>0$, either
\begin{itemize}
\item[(a)]  $\sup_j t_j = \infty$ and the conjugacy classes form 
	$$\Sol(\cD) = \sup_j \, \frac{2^j}{\max\{t_i: i\leq j\}}$$
solenoids in the moduli space $\cM_3$; or 
\item[(b)]  $\sup_j t_j < \infty$ and each conjugacy class is homeomorphic to a circle.
\end{itemize}
The computation of $\sup_j t_j$ and $\Sol(\cD)$ depends only on the tableau of $\cD$.
\end{theorem}

In this section, we give the proofs of Theorems \ref{deg3twist} and \ref{deg3infinite}.   We work with a simplified combinatorial object, the truncated spine.  We show how to compute the tableau and $\tau$-sequence and also the tree code from the truncated spine.  We conclude the section by comparing our constructions to those appearing in \cite{Branner:Hubbard:1}, \cite{Branner:Hubbard:2}, and \cite{Branner:cubics}.    

\subsection{The space of cubic polynomials}
Let $\cP_3\iso \C^2$ denote the space of monic and centered cubic polynomials.  It is a degree 2 branched cover of $\cM_3$.  Explicitly, a polynomial $f(z) = z^3 + a z + b$ is conformally conjugate to $g(z) = z^3 + a'z + b'$ if and only if they are conjugate by $z\mapsto -z$; consequently $a=a'$ and $b=-b'$.  Therefore, $\cM_3$ has the structure of a complex orbifold with underlying manifold $\C^2$; the projection $\cP_3\to \cM_3$ is given by $(a,b) \mapsto (a, b^2)$, and its branch locus $\{b=0\}$ is precisely the set of polynomials with a nontrivial automorphism (necessarily of the form $z\mapsto -z$).  Observe that the critical points for cubic polynomials with a nontrivial automorphism are interchanged by the automorphism; they therefore escape at the same rate.  

\subsection{The length of a cubic polynomial}
Fix a cubic polynomial $f$ with disconnected Julia set.  Denote its critical points by $c_1$ and $c_2$, so that $G_f(c_1) \geq G_f(c_2)$.  The {\em length} of $f$ is the least integer $L = L(f)$ such that $G_f(c_2) \geq G_f(c_1)/3^L$.  If no such integer exists, we set $L(f) = \infty$.  Thus, $L(f)=\infty$ if and only if $c_2$ lies in the filled Julia set of $f$; and $L(f) = 0$ if and only if $G_f(c_1) = G_f(c_2) > 0$.

\begin{lemma} \label{level 0 unique}
Let $\cD$ be a pictograph for a cubic polynomial of length $L(f)=0$.  There is a unique topological conjugacy class of polynomials in $\cM_3$ with pictograph $\cD$.  
\end{lemma}

\proof
Any length 0 cubic polynomial $f$ has $G_f(c_1) = G_f(c_2) = M$ for some $M>0$.  Fix $M>0$ and let $(\cF, \cX)$ be the unique tree of local models with pictograph $\cD$ and critical height $M$.  The underlying tree $(F, T)$ has a unique fundamental edge.  Let $v_0$ denote the vertex at height $M$, and set $v_1 = F(v_0)$.  To construct any polynomial $f$ with tree of local models $(\cF, \cX)$, we first glue the outer annulus of the local model surface $(X_{v_0}, \omega_{v_0})$ to the unique inner annulus of the local model surface $(X_{v_1}, \omega_{v_1})$.  The choice of gluing along the fundamental edge uniquely determines the gluing choices of all local models above $v_0$, because the local model maps must extend holomorphically.  Because $L(f)=0$, the local degree at all vertices below $v_0$ is 1, and therefore the choice of gluing along the fundamental edge also determines the gluing along every edge below $v_0$.  In other words, the gluing at the fundamental edge determines the conformal conjugacy class of an entire basin.  By uniformization, we may conclude that this gluing choice determines a unique point in $\cM_3$; see e.g. \cite[Lemma 3.4, Proposition 5.1]{DP:basins}.  

Finally, it is easy to see from the definition of the twisting deformation that all choices of gluing $X_{v_0}$ to $X_{v_1}$ can be obtained by twisting.  Therefore, all polynomials in $\cM_3$ with tree of local models $(\cF, \cX)$ are twist-conjugate.  Combined with Theorem \ref{compatible heights}, it follows that all polynomials in $\cM_3$ with pictograph $\cD$ are topologically conjugate.  
\qed

\subsection{Reducing to the case of 1 fundamental edge}  \label{reduction}
The main idea of the proofs of Theorems \ref{deg3twist} and \ref{deg3infinite} is the following.  We begin with a tree of local models $(\cF, \cX)$ with the given pictograph, and we fix a point in the base torus of the bundle of gluing configurations $\cB_3(\cF, \cX)$.  In the absence of symmetry, the basepoint corresponds to a unique choice of gluing along the fundamental edges of $(\cF, \cX)$.  That is, the conformal structure of the basin at and above the highest critical point is fixed.  

As usual, we let $v_0$ denote the highest branching vertex of the tree, so the local model surface at $v_0$ contains a critical point.  Inductively on descending height, we glue the local models along the spine of $(\cF, \cX)$ below the vertex $v_0$.  Each local model map below $v_0$ has degree 2, and there are exactly two ways to glue the local model surface so that the map extends holomorphically.  At each stage we determine which gluing choices are conformally conjugate and which gluing choices are topologically (twist) conjugate.  

It turns out that the only choices that contribute to our count of topological conjugacy classes are those made at the vertices of the spine that map to $v_0$.  Suppose the given tree has two fundamental edges, and let $v$ be a vertex in the spine that does not lie in the grand orbit of $v_0$.  Let $e$ be the edge above $v$.  If the forward orbit of the lower critical point does not contain $v$, then the two gluing choices along $e$ are easily seen to be conformally conjugate.  If the orbit of the lower critical point does contain $v$, then the distinct gluing choices for the local model at $v$ will always be topologically conjugate, as the following lemma shows.

\begin{lemma} \label{intermediate levels}
Fix a tree of local models for a cubic polynomial with two fundamental edges.  Let $v$ be a vertex in the spine, below $v_0$ and in the forward orbit of the lower critical point.  Suppose we have glued all local models in the spine above vertex $v$.  Then the two gluing choices of $(X_v, \omega_v)$ are topologically conjugate, via a conjugacy that preserves the orbit of the lower critical point and leaves the conformal structure above $v$ unchanged.  
\end{lemma}

\proof 
Let $f$ be any cubic polynomial with the given tree of local models $(\cF, \cX)$.  As $(\cF, \cX)$ has two fundamental edges, $f$ lies in the shift locus and the two critical points have distinct escape rates; let $c_1, c_2$ denote the critical points so that $G_f(c_1) > G_f(c_2)$.  Let $L=L(f)$ be the length of $f$.  The fundamental annulus is decomposed into two subannuli,
	$$A_0^2 = \{G_f(f^{L}(c_2)) < |z| < 3 \, G_f(c_1) \}$$
and
	$$A_0^1 = \{G_f(c_1) < |z| < G_f(f^{L}(c_2))\},$$
which can be twisted independently.  

For each $0 < n < L(f)$, denote by $A_n$ the annular component of $\{G_f(c_1)/3^n < |z| < G_f(c_1)/3^{n-1}\}$ separating the two critical points.  Let $0 < n_1 < \cdots < n_m < L(f)$ index the values of $n$ for which the orbit of $c_2$ intersects $A_n$.   The vertex $v$ corresponds to a level curve of $G_f$ in the annulus $A_{n_i}$ for some $i \in \{1, \ldots, m\}$.  We proceed inductively on $i$.  

For $i=1$, a full twist in $A_0^1$ followed by a full twist in the opposite direction in $A_0^2$ induces the opposite gluing choice for the local model surface $(X_v, \omega_v)$ without affecting the conformal class of the basin above $v$.  This is because the annulus $A_{n_1}$ maps with degree 2 by $f^{n_1}$ to the annulus $A_0$; a full twist in $A_0$ induces a half-twist in $A_{n_1}$.  

Similarly for each $i$: the annulus $A_{n_i}$ is mapped with degree $2^i$ by $f^{n_i}$ to $A_0$.   Therefore, $2^{i-1}$ full twists in $A_0^1$ followed by $-2^{i-1}$ full twists in $A_0^2$ will induce a half twist at the level of $v$ while preserving the conformal structure of the basin above $v$.   
\qed

\subsection{Truncated spines}  \label{truncated spine}
\begin{figure} 
\includegraphics[width=0.9in]{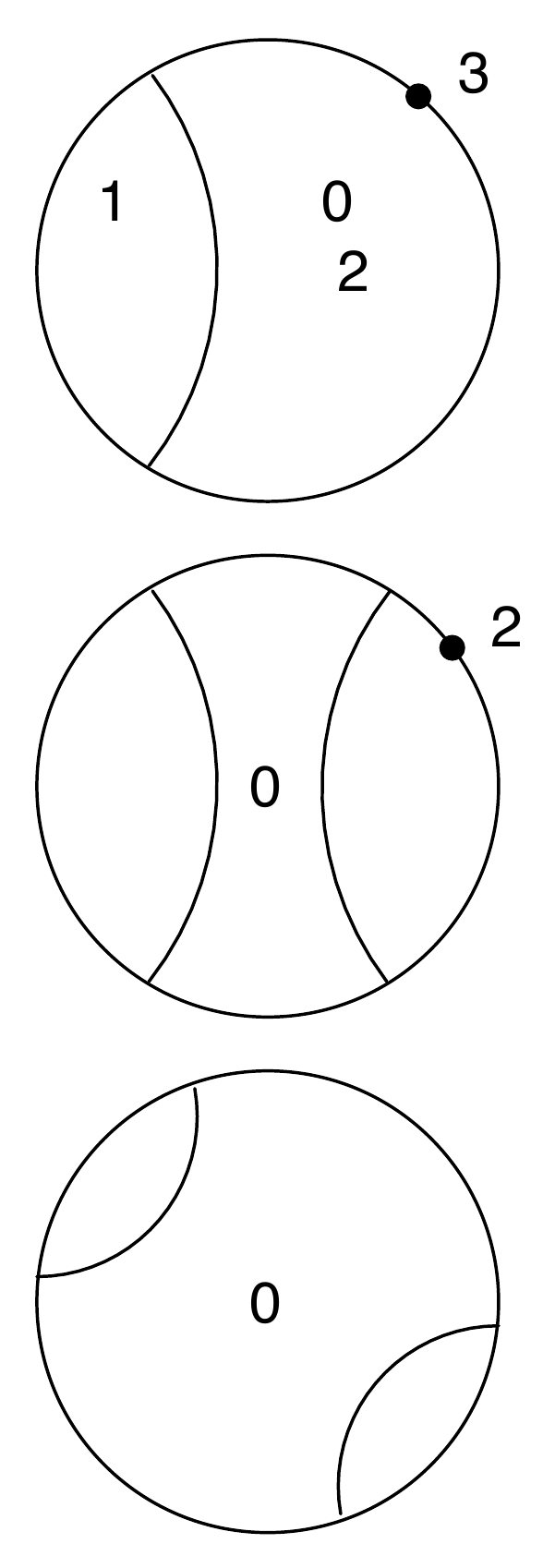}
\caption{The truncated spine for the pictograph in Figure \ref{nongeneric cubic}.} \label{truncated}
\end{figure}

\begin{figure} 
\includegraphics[width=4.75in]{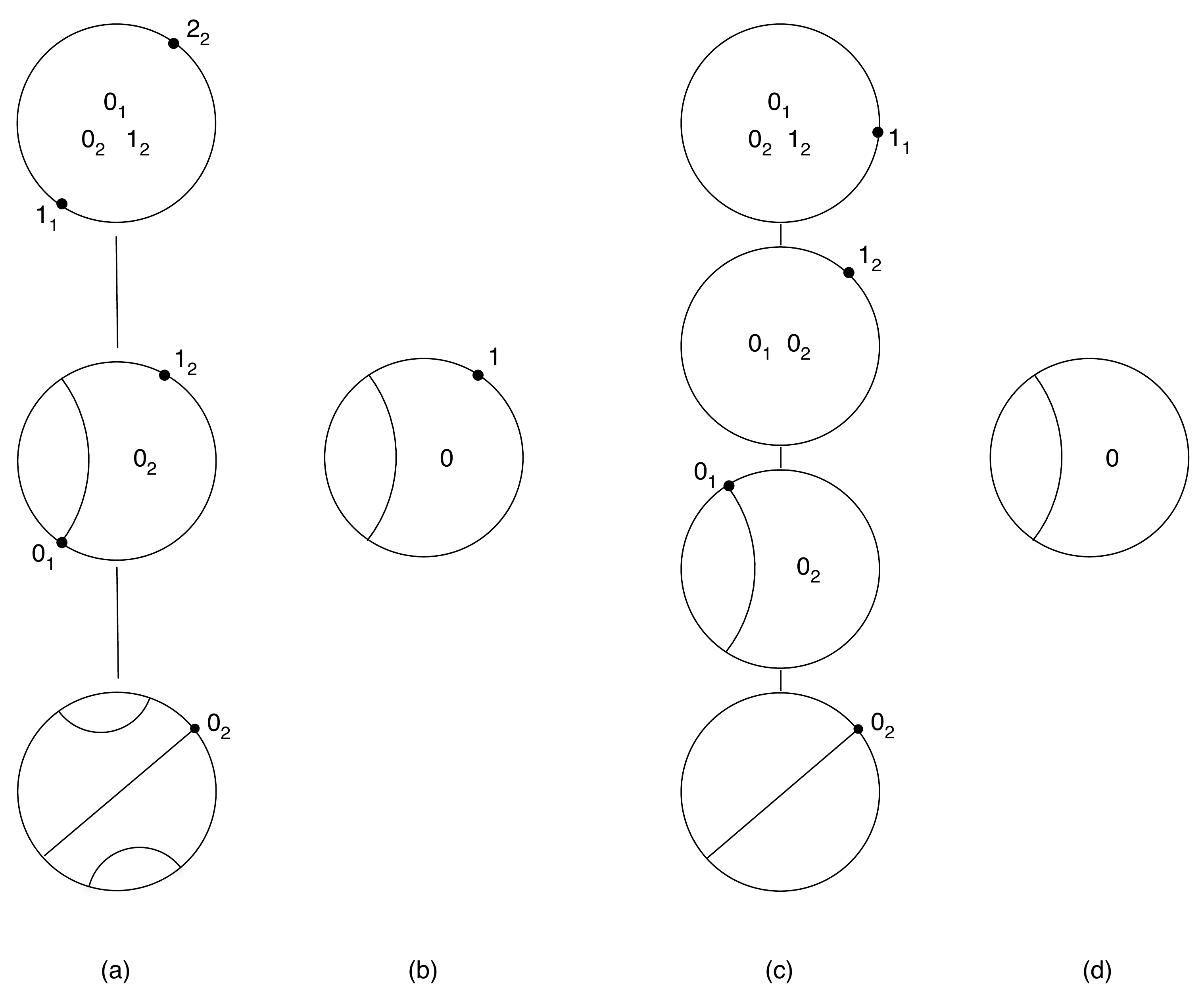}
\caption{Length 1 cubics: (a) a pictograph, with critical heights $(M, M/3)$; (b) the truncated spine for the pictograph of (a); (c) the generic pictograph at level 1, with critical heights $(M, M')$ where $M/3 < M' < M$; (d) the truncated spine for the pictograph of (c).} \label{length 1}
\end{figure}

The {\em truncated spine} of a cubic polynomial is a subset of its pictograph:  it consists of the laminations at vertices of height $G_f(c_1)/3^n$ for all $0 \leq n < L(f)$, where $L(f)$ is the length of $f$.  The lamination at level $n$ (height $G(c_1)/3^n$) is labelled by integer $0 \leq k \leq L(f) - n$, if $f^k(c_2)$ lies in one of its gaps or on the level curve.    If the two critical points have the same level, the truncated spine is an empty diagram.  Figure \ref{truncated} shows the truncated spine for the example of Figure \ref{nongeneric cubic}.  Figure \ref{length 1} contains pictographs for two cubic polynomials with length $L(f) =1$ and the truncated spine for each of them.

For drawing diagrams, it is more convenient to use the truncated spine rather than the full pictograph.  The next lemma shows that we do not lose any information by doing so.  

\begin{lemma} \label{spine from truncated}  
Except for length 0 cubics, the pictograph can be recovered from the truncated spine.  
\end{lemma}

\proof
For any polynomial $f$ with length $L(f)>0$, the truncated spine is non-empty.  The lamination at height $3 \, G_f(c_1)$ with its marked points is uniquely determined by the marked lamination at height $G_f(c_1)$ (see Lemma \ref{equivalence}).  For $L(f) = \infty$, the truncated spine together with the lamination at height $3 \, G_f(c_1)$ is the complete pictograph.

Assume $f$ is a cubic polynomial with finite length $L(f)>0$.  Observe that there are marked points on the lamination at height $G_f(c_1)$ (not just in the gaps) if and only if $f$ has one fundamental edge, meaning that $G_f(c_2) = G_f(c_1)/3^{L(f)}$.  In this case, the truncated spine is almost the complete pictograph.   The lamination diagram containing the critical point $c_2$ is also uniquely determined:  it is a degree 2 branched cover of the lamination at level $n$, where $n$ is the greatest integer such that the symbol $L(f) - n$ marks the lamination at height $G_f(c_1)/3^n$, branched over the marked point.  

In the case of two fundamental edges, it is easy to see how to fill in the pictograph. We first add the subscript ``2" to teach of the labels in the truncated spine, and we mark the unique non-trivial equivalence class at the height of $c_1$ with the symbol $0_1$.  We next include the trivial lamination (a circle) to the column of lamination diagrams above each lamination of the truncated spine.  These laminations correspond to vertices in the spine intersecting the grand orbit of the second critical point $c_2$, except at the height of $c_2$ itself.  The lamination at the height $G_f(c_2)$ will be the ``figure 8": it is a circle cut by a diameter; the diameter is marked with the symbol $0_2$.  Finally, at each height $G_f(c_1)/3^n$ where the symbol $(L(f) - n - 1)_2$ appears in the same component as the symbol $0_2$, we mark the trivial lamination below it with the symbol $(L(f) - n - 1)_2$, in one of two ways:  (1) if $(L(f) - n - 1)_2$ also appears at height $G_f(c_1)/3^{n+1}$ then we place $(L(f) - n - 1)_2$ in the gap, and (2) otherwise, we place the symbol on the lamination circle.

From the definition of the pictograph, we see that this is the complete diagram.
\qed

\subsection{Marked levels and the $\tau$ sequence}  \label{marked levels}
Fix a cubic polynomial $f$ with a disconnected Julia set, and recall that the length $L(f)$ is the least integer $L$ such that $G(c_2) \geq G(c_1)/3^L$; we set $L(f) =\infty$ if $G(c_2)=0$.  For each $0\leq n\leq L(f)$, the {\em level $n$ puzzle piece} $P_n$ is the connected component of $\{G < G(c_1)/3^{n-1}\}$ containing $c_2$.  Note that $P_0$ contains both critical points and $f^n(P_n) = P_0$ for all $n$.  The Yoccoz $\tau$-function associated to $f$,
	$$\tau_f: \{1, 2, 3, \ldots, L(f)\} \to \mathbb{N},$$
is defined by the following:  let $k(n) = \min\{k>0:  c_2 \in f^k(P_n)\}$; then $f^{k(n)}(P_n) = P_{\tau_f(n)}$.  In particular, $\tau_f(1)=0$ and $\tau_f(i+1) \leq \tau(i)+1$ for all $i$ and every $f$ with $L(f)>0$.  The data of $\tau_f$ is equivalent to the information in the tableau (or marked grid) of $f$ defined in \cite{Branner:Hubbard:2}.

We defined the notion of marked levels in the introduction, but we recall the definition here.  For each $n\geq 0$, let $B_n$ be the closed subset of the level $n$ puzzle piece $P_n$ where $G_f(z) \leq G_f(c_1)/3^n$.  A {\em marked level} is an integer $n>0$ where the orbit of $c_2$ intersects $B_n\setminus P_{n+1}$.

\begin{lemma}  \label{tau marked}
A level $n > 0$ is marked if and only if at least one of the following holds: 
\begin{enumerate}
\item	there exists $i<L(f)$ so that $\tau(i) = n$ and $\tau(i+1) \leq n$;  
\item	$G_f(c_2) = G_f(c_1)/3^{L(f)}$ and $n = \tau^k(L(f))$ for some $k>0$.  
\end{enumerate}
\end{lemma}

\noindent
The proof is immediate from the definitions.  From this lemma, we see that a marked level coincides with the ``semi-critical" levels of \cite{Milnor:lc}.   These marked levels also coincide with the ``off-center" levels $n_0, n_1, \ldots, n_j$ of $f$ of \cite[Theorem 9.1]{Branner:cubics}.  

In terms of the truncated spine, we may characterize the marked levels as follows.  

\begin{lemma}
The marked levels coincide with the degree 2 vertices of the truncated spine where the order of rotational symmetry is 1.
\end{lemma}

\noindent 
Again the proof is immediate from the definitions.  The marking of a level $n$ means that the symmetry is broken at that level.

\subsection{Counting topological conjugacy classes} \label{counting cubics}
We begin by defining relative moduli and twist periods, quantities involved in the computations of Theorems \ref{deg3twist} and \ref{deg3infinite}.  

Let $(\cF, \cX)$ be a tree of local models with a given pictograph $\cD$.  Let $\cS$ be the spine of $(\cF, \cX)$ and suppose we have chosen a gluing of all the local models along $\cS$.  From Proposition \ref{gluing the spine}, the gluing choices along $\cS$, together with its first-return map, uniquely determine a complete basin of infinity $(f, X(f))$.  Let $c_1$ and $c_2$ denote the critical points of $f$ so that $G_f(c_1)\geq G_f(c_2)$.  Note that $c_2$ may be an end of $X(f)$.  Let $L(f)$ be its length.  Let 
	$$A_0 = \{z: G_f(c_1) < |z| < 3 \, G_f(c_1)\}$$
denote the fundamental annulus.  For each $0 < n < L(f)$, denote by $A_n$ the annular component of $\{G_f(c_1)/3^n < |z| < G_f(c_1)/3^{n-1}\}$ separating the two critical points.

For each $0 \leq n < L(f)$, the {\em relative modulus} at level $n$ is the ratio
	$$m(n) =  \mod(A_n)/\mod(A_0).$$
Note that $m(n)$ is completely determined by the $\tau$-sequence $\tau_f$.  In fact, $m(n) = 2^{-k(n)}$, where $k(n)$ the least integer such that $\tau^{k(n)}(n) = 0$.  That is, $k(n)$ counts the number of times the orbit of $A_n$ returns to the critical nest.  

Note also that a full twist in $A_0$ induces a $m(n)$-twist in the annulus $A_n$.  For each $n < L(f)$, the {\em twist period} of the basin $(f, X(f))$ at level $n$ is the minimum number of twists $T_n>0$ in the fundamental annulus $A_0$ that returns all marked levels $\leq n$ to their original gluing configuration.  This means that the induced twist, summing along the annuli of the critical nest down to each marked level $j\leq n$, must be integral.   

We shall see in the proof of Theorem \ref{deg3twist} that these twist periods can be computed from the $\tau$-sequence and are independent of the choice of gluing configuration.  

\medskip
We are now ready to prove Theorems \ref{deg3twist} and \ref{deg3infinite}.  Our proofs follow the same reasoning as the proofs of \cite[Theorem 9.1]{Branner:cubics} and \cite{Harris:cubics}, but we need also to treat the case of finite length critical nests.

\bigskip\noindent{\em Proof of Theorem \ref{deg3twist}.}  
Let $\cD$ be a pictograph with finitely many marked levels.  Fix $M>0$, and let $(\cF, \cX)$ be any tree of local models with pictograph $\cD$ and maximal critical height $M$.  Let $L(f)$ denote the length of any cubic polynomial basin with tree of local models $(\cF, \cX)$.  If both critical points have the same height, the length is $L(f)=0$, there are no marked levels, and we are done by Lemma \ref{level 0 unique}.  We may assume that the critical points have distinct heights.  

Fix a point in the base torus for the bundle $\cB_3(\cF, \cX)$.  By Lemma \ref{lemma:cubic tlm symmetries}, there are no symmetries at the fundamental vertices; it follows that the base torus parametrizes the gluing choices in the fundamental edges.  From the structure of the bundle of gluing configurations, it follows that each point in the fiber corresponds to a unique conformal conjugacy class of basins (with the chosen gluing configuration above height $M$); the topological conjugacy classes are in one-to-one correspondence with their orbits under twisting.  As described in \S\ref{reduction}, we will proceed inductively, on descending height, to glue the local models along the spine.  By Lemma \ref{intermediate levels}, we may disregard the gluing choices at the ``intermediate levels", corresponding to vertices in the grand orbit of the lower critical point, when there are two fundamental edges.  At each vertex of height $M/3^n$, for integers $0 < n < L(f)$, we will compute the number of distinct topological conjugacy classes arising from the two gluing choices.

First assume that $(\cF, \cX)$ has no marked levels.    Then for every $n$, the two gluing choices of the local model at height $M/3^n$ are conformally equivalent.  If there are two fundamental edges, we conclude from Lemma \ref{intermediate levels} that there is a unique topological conjugacy class.  If there is only one fundamental edge, then in fact there is a unique point in any fiber of the bundle $\cB_3(\cF, \cX)$, so clearly $\Top(\cD) = 1$.  Note that the absence of marked levels can be discerned from the $\tau$-sequence of $\cD$, by Lemma \ref{tau marked}. 

Suppose now that $(\cF, \cX)$ has $k>0$ marked levels, and let 
	$$0 < l_1 < l_2 < \cdots < l_k < L(f)$$ 
denote the marked levels.  Let $m_j$ be the sum of relative moduli to level $l_j$:
	$$m_j = \sum_{n=1}^{l_j} m(n).$$
Let $t_j$ be the smallest positive integer so that $t_jm_j$ is integral; it is always a power of 2, because the vertices of the critical nest are mapped with degree 2 to their images.  Then the twist period at level $l_j$ is easily seen to be the maximum 
	$$T_{l_j} = \max\{t_i: i\leq j\}.$$

\begin{lemma} \label{times 2}
We have $T_{l_j} \in \{T_{l_{j-1}}, 2T_{l_{j-1}}\}$ for all $j$.
\end{lemma}

\proof
Suppose it takes $t$ full twists to return to a given configuration at marked level $l$, and suppose the next marked level is $l'>l$.  The level $\tau(l')$ must also be a marked, so $\tau(l')\leq l$.  Thus, after $t$ twists, $\tau(l')$ is in its original configuration; therefore $l'$ must be either in its original configuration or twisted halfway around.  Therefore, at most $2t$ twists are needed to return to the original configuration at level $l'$.   
\qed

\medskip
Lemma \ref{times 2} says that the twist period between marked levels can increase at most by a factor of two.  It follows that twisting reaches all possible gluing configurations if and only if $T_{l_k} = 2^k$.  Or, more precisely, the number of distinct topological conjugacy classes associated to the given tree of local models is the ratio $2^k/T_{l_k}$.  This completes the proof of the theorem.
\qed

\bigskip\noindent{\em Proof of Theorem \ref{deg3infinite}.}  
Suppose $\cD$ is a pictograph wth infinitely many marked levels.  For any fixed value $M>0$, there is a unique tree of local models $(\cF, \cX)$ with the given pictograph $\cD$ and critical point at height $M$.  Note that the infinitely many marked levels implies, in particular, that one critical point lies in the filled Julia set.  

As in the proof of Theorem \ref{deg3twist}, we begin by fixing a point in the base torus of $\cB_3(\cF, \cX)$.   Note that the base torus is a circle in this case.  At each marked level of the pictograph, the two gluing choices produce conformally inequivalent basins in $\cB_3(\cF, \cX)$; it follows that there are infinitely many points in any fiber of $\cB_3(\cF, \cX)$.  From Lemma \ref{Cantor fiber}, the fiber is then homeomorphic to a Cantor set;  in particular, there are uncountably many points in the fiber.  On the other hand, two basins in a fiber are topologically conjugate if and only if they lie in the same twist orbit; thus, there can only be countably many topologically conjugate points in a fiber.  Therefore, there are infinitely many topological conjugacy classes within $\cB_3(\cF, \cX)$.  

As proved in \cite{Branner:Hubbard:2} and discussed further in \cite{Branner:cubics} and \cite{Harris:cubics}, the topological conjugacy classes in the cubic moduli space $\cM_3$ organize themselves into solenoids or a union of circles, depending on the twist periods.  In the notation of the proof of Theorem \ref{deg3twist}, the twist period of a complete basin with tree $(\cF, \cX)$ is the value 
	$$T = \lim_{j\to\infty} T_{l_j}.$$
If $T$ is infinite, then necessarily the bundle of gluing configurations is a union of 
	$$\Sol(\cD) := \lim_{j\to\infty}  \frac{2^j}{T_{l_j}}$$
solenoids.  Note that the limit $\Sol(\cD)$ exists and lies in the set $\{1, 2, 2^2, 2^3, \ldots, \infty\}$, because $2^j/T_{l_j}$ defines a non-decreasing sequence of powers of 2.  For $T<\infty$, the bundle of gluing configurations forms an infinite union of closed loops.
\qed

\subsection{Combinatorial relations: from truncated spines to trees and tableaux}  \label{treecode}
It is easy to read the tableau and tree code from the truncated spine of a cubic polynomial.  Suppose first that the polynomial lies in the shift locus, so the length $L(f)$ is finite.  Then its tableau (and $\tau$-sequence) is finite.  As defined in \cite{Branner:Hubbard:2}, the tableau is a subset of the 4th quadrant of the $\mathbb{Z}^2$ lattice:  a lattice-point $(i,j)$ with $i\geq 0$ and $j\leq 0$ is marked if the lamination at level $-j$ contains the integer symbol $i$.  Furthermore, when $i-j = L(f)$, we mark $(i,j)$ if the lamination at level $-j-1$ contains the integer symbol $i$ in its {\em central} gap (the gap containing the symbol 0).  

In terms of the $\tau$-sequence, we have $\tau(1) = 0$ for every polynomial with length $L(f)>0$.  For $L(f)>1$ and each $0 < n < L(f)$, we set 
	$$\tau(n) = \max\{j: \mbox{lamination at level } j \mbox{ is labelled by } (n-j)\}.$$
To compute $\tau(L(f))$, we consider the set
	$$\cL = \{j : \mbox{the central gap in lamination at level } j \mbox{ is labelled by } (L(f) - j - 1)\}.$$
We then have 
	$$\tau(L(f)) = \left\{ \begin{array}{ll}
		1 + \max\{j: j\in \cL \} & \mbox{ if } \cL \not= \emptyset \\
		0	& \mbox{ if } \cL = \emptyset  \end{array} \right. $$
As an example, the $\tau$-sequence for the truncated spine of Figure \ref{truncated} is $(0, 0, 1)$.

The tree code for a cubic polynomial of length $L(f)$ is a sequence of pairs $(k(i), t(i))$, where $i$ ranges from 1 to $L(f)$, defined in \cite[\S11]{DM:trees}.  Fix a truncated spine for a polynomial of length $L(f)$.  A {\em minimal symbol} in a gap of a lamination is the smallest integer in a labelled gap.  The lifetime $k(i)$ is equal to the number of times the symbol $j$ appears as a minimal symbol in a gap at level $i-j-1$, as $j$ ranges from 0 to $i-1$.  In particular, $k(1)=1$.  The terminus $t(i)$ is computed as follows:
\begin{enumerate}
\item	Let $j(i)$ be the smallest $j$ which appears at level $i-j-1$ but is {\em not} a minimal symbol.  If such a $j$ does not exist, then let $j(i)=i$.
\item	Let $m(i)$ be the minimal symbol at level $i-j(i)-1$ in the gap containing $j(i)$.  When $j(i)=i$, set $m(i) = 0$.  
\item	Let $t(i) = i - j(i) + m(i)$.
\end{enumerate}
As an example, the tree code for the truncated spine of Figure \ref{truncated} is $(1,0), (2,0), (1, 1)$.

\subsection{The Branner-Hubbard description}
Branner and Hubbard showed that there are two important dynamically-defined fibrations in the space of monic and centered cubic polynomials $\cP_3$.  Let $\cC_3$ denote the connectedness locus, the set of polynomials in $\cP_3$ with connected Julia set.  First, the maximal critical escape rate 
	$$M: \cP_3\setminus \cC_3 \to (0,\infty)$$
defined by
	$$M(f) = \max\{G_f(c): f'(c)=0\}$$
is a trivial fibration with fibers homeomorphic to the 3-sphere \cite[Theorem 6.1]{Branner:cubics} (which follows from \cite[Theorem 11.1, Corollary 14.6]{Branner:Hubbard:1}).  Branner and Hubbard analyzed the quotient of a fiber of $M$ in $\cM_3$; it follows from \cite[Cor 14.9]{Branner:Hubbard:1} that the induced map 
	$$M: \cM_3\setminus \cC_3 \to (0,\infty)$$
is also trivial fibration with fibers homeomorphic to the 3-sphere.  The trivialization is given by the stretching deformation; see \S\ref{qc}.

For each $r>0$, let $H_r\subset M^{-1}(r)\subset \cP_3$ be the locus of polynomials with $G(c_2)< G(c_1) = r$.  Let $c'_1$ denote the {\em cocritical point} of $c_1$, so that $f^{-1}(f(c_1)) = \{c_1, c_1'\}$, and let $\theta(c'_1) \in \mathbb{R} / 2\pi \mathbb{Z}$ be its external angle.  Then
	$$\Phi_r: H_r \to S^1$$
defined by $\Phi_r(f) = \theta(c'_1)$ is a trivial fibration with fibers homeomorphic to the unit disk $\D$ \cite[Theorem 6.2]{Branner:cubics}.  The fiber of $\Phi_r$ over $\theta$ will be denoted $F_r(\theta)$.  Note that every polynomial in $F_r(\theta)$ is conjugate by $z\mapsto -z$ to a unique polynomial in $F_r(\theta + \pi)$.  (It is worth observing that the polynomials with nontrivial automorphism cannot be in $H_r$; either $f(z) = z^3$ or $z\mapsto -z$ interchanges the two critical points, and they therefore escape at the same rate.)

The turning deformation of \cite{Branner:Hubbard:1} induces a monodromy action on a fiber $F_r(\theta)$; its first entry into $F_r(\theta + \pi)$ determines the {\em hemidromy action}.  Alternatively, the hemidromy action is the monodromy of the induced fibration on the quotient of $H_r$ in $\cM_3$:
	$$\Phi_r: [H_r] \to S^1$$
given by $\Phi_r(f) = 2\theta(c_1') \modspace 2\pi$ which is well-defined on the conjugacy class of $f$. The fibers of $\Phi_r$ in $\cM_3$ are again topological disks.  

\begin{figure}
\includegraphics[width=2.3in]{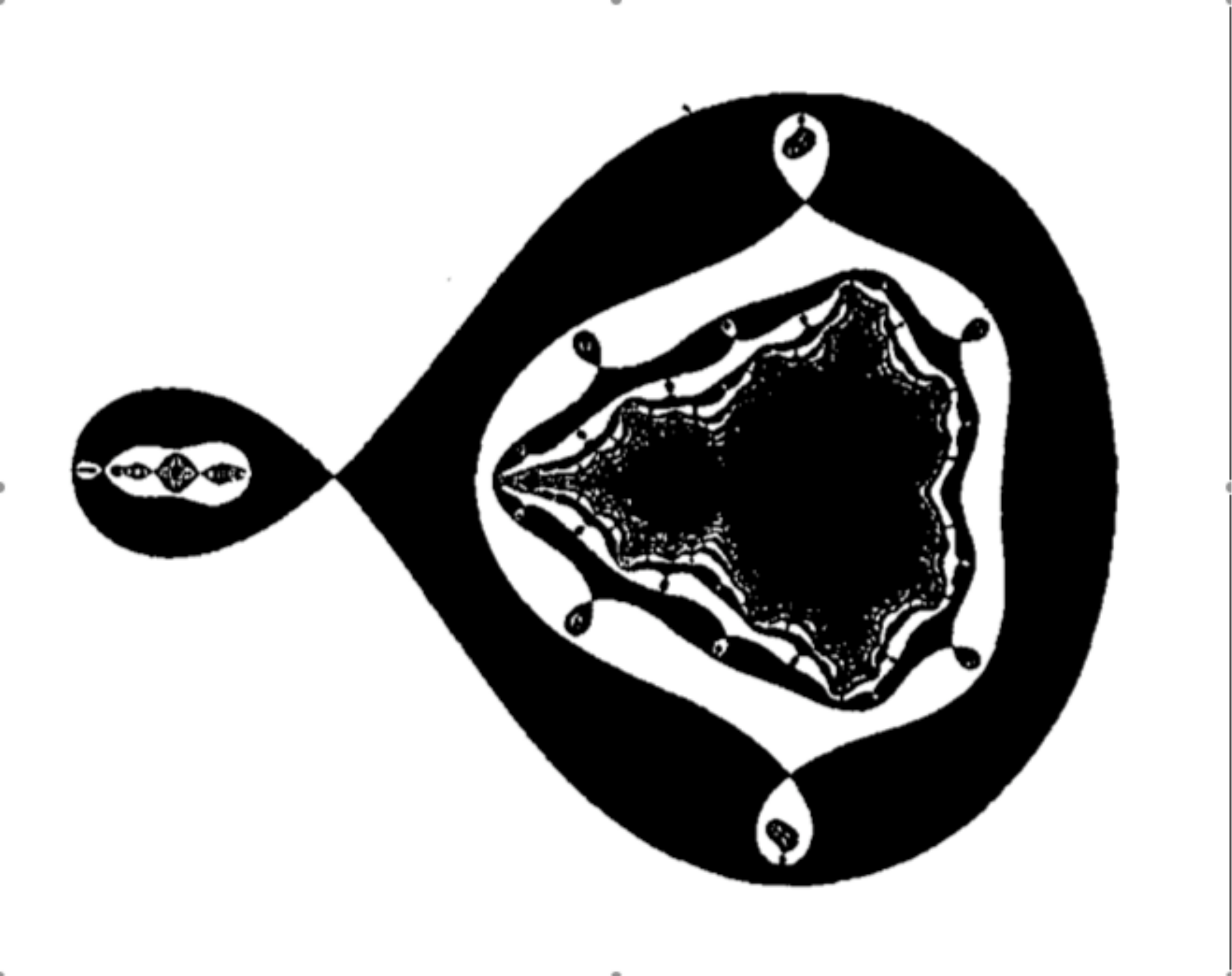}
\caption{A fiber of $\Phi_r$ showing critical level sets of $f\mapsto G_f(c_2)$, from \cite[Figure 9.3]{Branner:Hubbard:2}.}
\end{figure}

Let us fix $\theta = 0$ and consider the fiber $[F_r(0)]$ in $\cM_3$.  The hemidromy action on this fiber will be denoted 
	$$T_r: [F_r(0)] \to [F_r(0)].$$
It corresponds to a full twist in the fundamental annulus; compare \S\ref{qc}.  The escape rate of the second critical point further decomposes $[F_r(0)]$.  The critical level sets of $G(c_2)$ are precisely the levels $r/3^n$ for integers $n>0$.  The connected components of $\{G(c_2) < r/3^n\}$ are called the {\em level n disks}.  The hemidromy action permutes these disks.  The {\em period} of a level $n$ disk $D$ is the least number of iterates $p>0$ such that $T_r^p(D) = D$.  Branner-Hubbard showed that these periods are always powers of 2.  The period of the level $n$ disk coincides with the twist period $T_n$ (defined in \S\ref{counting cubics}) for any cubic polynomial in that disk.  

Suppose $(\cF, \cX)$ is a cubic tree of local models with both critical heights positive.  If there is only one fundamental edge, then the Branner-Hubbard turning curves through any polynomial $f$ with tree $(\cF, \cX)$ constitute the connected components of the bundle $\cB_3(\cF, \cX)$.  In fact, twisting coincides with the turning deformation, up to the normalization of the parametrization.  For a cubic tree $(\cF, \cX)$ with two fundamental edges, the base torus of $\cB_3(\cF, \cX)$ is two-dimensional.  Intersecting with a fiber of the Branner-Hubbard bundle $\Phi_r$, the bundle $\cB_3(\cF, \cX)$ consists of finitely many connected components of a level set of $f\mapsto G_f(c_2)$.  In this case, a full turn around the base of $\Phi_r: [H_r]\to S^1$ corresponds to a twist by 
	$$\frac{m_1}{m_1+m_2} \, {\bf e}_1 + \frac{m_2}{m_1 + m_2} \, {\bf e}_2$$
where $m_i$ is the modulus of the fundamental annulus $A_i$, so that $m_1 + m_2 = 2r$.

Among the polynomials $f$ with length $L(f) = \infty$, there are two types of connected components in the Branner-Hubbard slice $F_r(0)$.  They showed that the Mandelbrot sets in their picture correspond to cubic polynomials where the connected component of the filled Julia set containing the critical point is periodic.  Equivalently, the tableau is periodic.  The twist periods of these Mandelbrot sets are always finite; there are only finitely many marked levels in the corresponding pictograph.  In this case, the bundle $\cB_3(\cF, \cX)$ is a finite union of circles.  Branner and Hubbard proved that all other cubic polynomials with infinite length correspond to points in their slice.  For these polynomials, the bundle $\cB_3(\cF, \cX)$ may be a union of circles or a union of solenoids.

%%%%%%
%%%%%%

\bigskip\bigskip\section{Examples in degree 3} \label{sec:examples}

In this section, we give examples of truncated spines for cubic polynomials; the definition of the truncated spine can be found in \S\ref{truncated spine}.  These examples illustrate the existence of cubic polynomials with 
\begin{enumerate}
\item	the same tableau (or $\tau$-sequence) but different trees $(F, T)$,
\item	the same tree $(F, T)$ but different trees of local models $(\cF, \cX)$, and
\item	the same tree of local models but different topological conjugacy classes.
\end{enumerate}
The examples we provide are structurally stable in the shift locus; they are the shortest examples that exist.  We give a final example with an infinite $\tau$-sequence, where the bundle of gluing configurations for any associated tree of local models $(\cF, \cX)$ forms exactly two solenoids in the moduli space $\cM_3$.

\begin{figure} 
\includegraphics[width=2.5in]{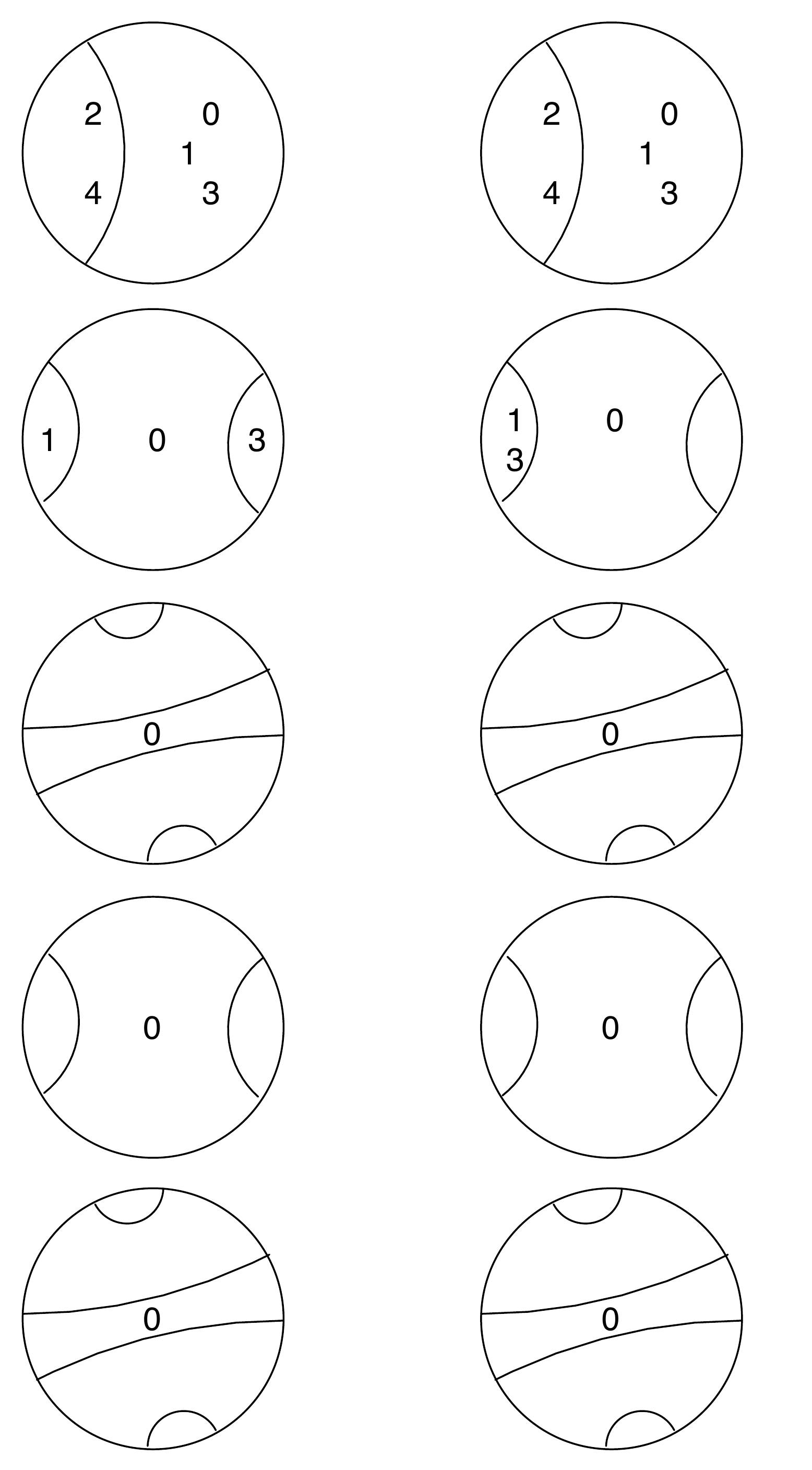}
\caption{Truncated spines associated to two different trees with the same $\tau$-sequence $0, 1, 0, 1, 0$.} \label{ex1}
\end{figure}

\begin{figure}
\includegraphics[width=2.8in]{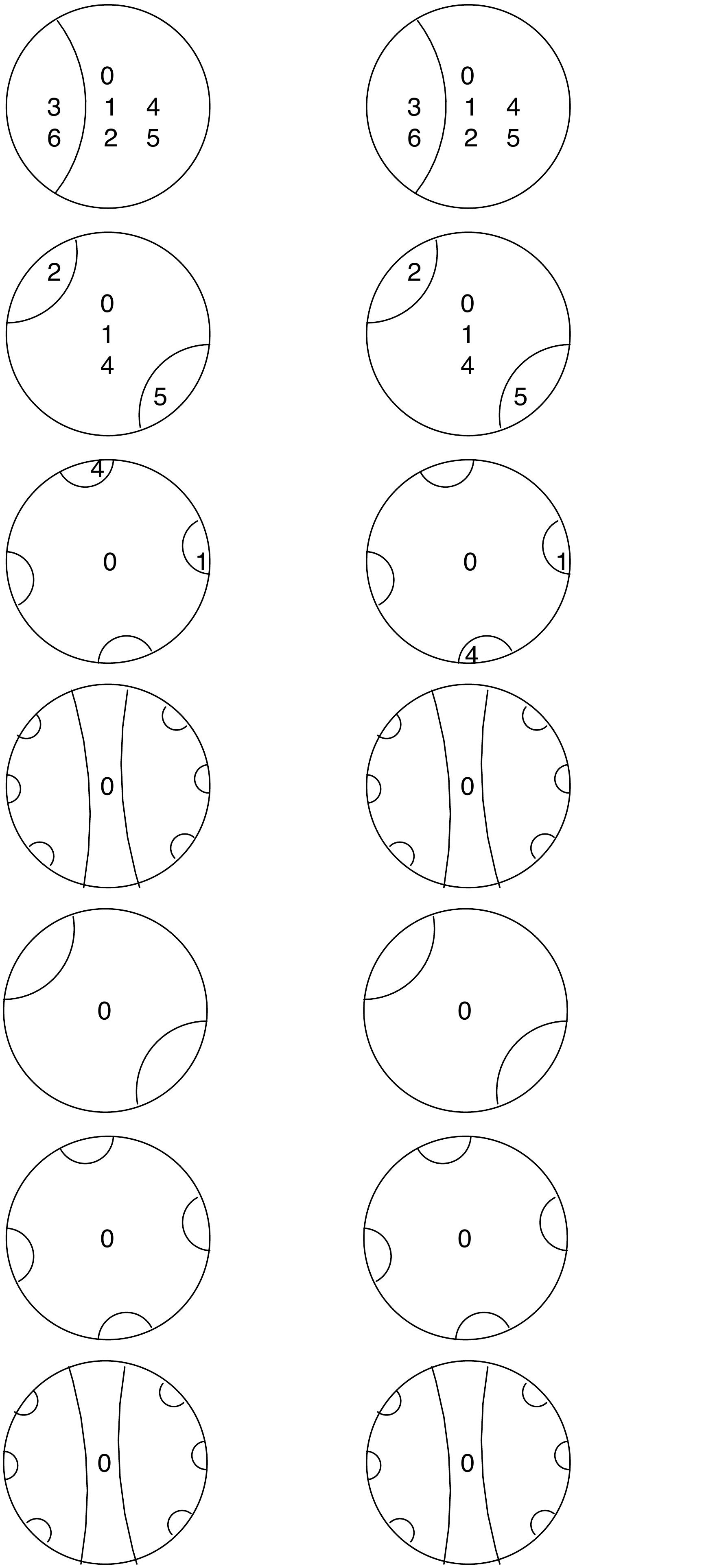}
\caption{Inequivalent truncated spines associated to the same tree.  The two spines differ in the cyclic ordering of the 1 and 4 at level 2 (the third lamination).}  \label{ex2}
\end{figure}

\begin{figure}
\includegraphics[width=2.2in]{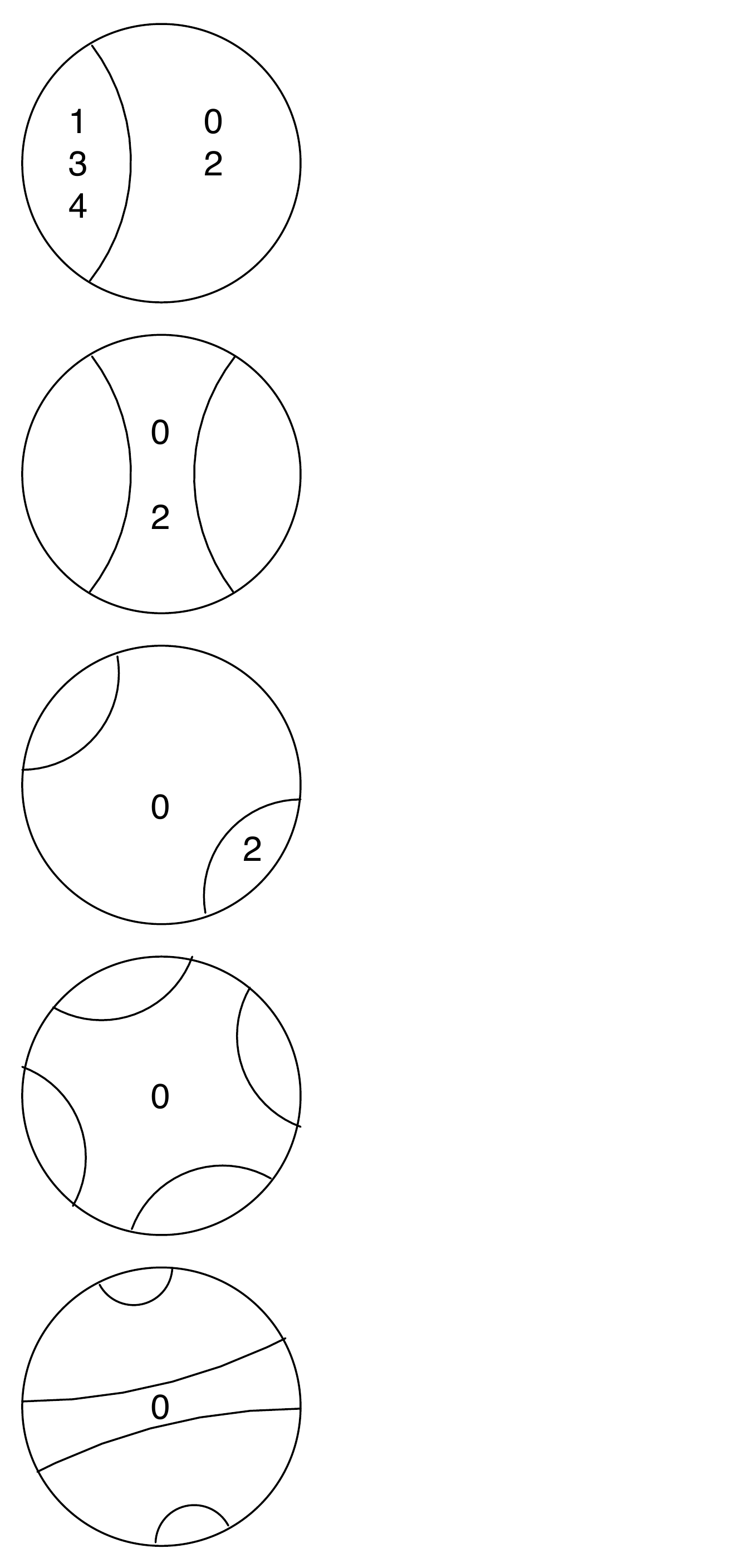}
\caption{Truncated spine of length 5 determining two topological conjugacy classes}  \label{ex3}
\end{figure}

\subsection{Examples from Figures \ref{ex1}, \ref{ex2}, and \ref{ex3}}
The tree codes for truncated spines in Figure \ref{ex1}, using the computation of \S\ref{treecode},  are 
	$$(1,0), (1,1), (3,0), (1,1), (2,3)$$
for the spine on the left, and
	$$(1,0), (1,1), (3,0), (1,1), (1,3)$$
for the spine on the right.  These examples are also presented in \cite[\S11]{DM:trees}.  As tree codes characterize the trees of a given level, this shows the trees are different.  The $\tau$-sequence for both truncated spines is $0,1,0,1,0$.

The second example (Figure \ref{ex2}) shows two truncated spines associated to the same tree.  The tree code for these spines is $(1,0), (1,1), (1,2), (4,0), (1,1), (1,2), (3,4)$.  The difference is in the relative locations of the first and fourth iterates of critical point $c_2$ in the puzzle piece at level 2.

In Figure \ref{ex3}, we give an example of a truncated spine determining two distinct conjugacy classes.  The $\tau$-sequence for this example is $(0,0, 1, 2, 0)$.  Following the algorithm of Theorem \ref{deg3twist}, we first enumerate the marked levels $l_0 = 0, l_1 = 2$; then compute the sums of relative moduli $m_0 = 0, m_1=1$; then compute $t_0 = 1, t_1=1$.  The number of conjugacy classes is the maximum of $2^0/1$ and $2^1/1$; therefore there are 2 conjugacy classes determined by this truncated spine.

It is not hard to show that these examples are the shortest of their type; that is, any $\tau$-sequence giving rise to more than one tree must have length $\geq 5$; any two truncated spines giving rise to the same tree must have length $\geq 7$; any truncated spine giving rise to more than one conjugacy class must have length $\geq 5$.  Indeed, one can easily compute by hand all combinatorial possibilities to length 6.   An enumeration of all admissible $\tau$-sequences, truncated spines, and stable topological conjugacy classes to length 21 is given in \cite{DS:count}, implementing an algorithm derived from Theorem \ref{deg3twist}, while an enumeration of all cubic trees to length 17 was given in \cite{DM:trees}.

\subsection{Solenoidal twist classes} 
This final example is similar to the Fibonacci solenoid, described in \cite{Branner:Hubbard:2}, but here we provide a $\tau$-sequence that determines exactly two solenoids in $\cM_3$.  Let $l_0=0$, $l_1 = 2$, $l_2=4$, and $l_j = 2l_{j-1}+1$ for all $j>2$.  Consider the $\tau$-sequence given by 
	$$0,0,1,2,0,1,2,3,4,0,1,2, \,\cdots,9,0,1,\,\cdots,l_4,0,1,\,\cdots,l_5,0,1,\,\cdots,l_6,0,1,\,\cdots$$
It can be proved inductively that this sequence determines a unique pictograph; the first five laminations of the truncated spine (for the $\tau$-sequence $(0,0,1,2,0)$) are shown in Figure \ref{ex3}.  Therefore, for any critical escape rate $M>0$, the $\tau$-sequence determines a unique tree of local models.  The marked levels are given by the sequence $\{l_j: j > 0\}$.  Computing inductively, the relative moduli sums are $m_1 = 1$, $m_2 = 3/2$, and $m_j = m_{j-1} + 1/2 + m_{j-1}/2$ for all $j>2$.  Therefore $t_1 = 1$, $t_2 = 2$, and $t_j = 2^{j-1}$ for all $j$.  Because the twist  periods $T_{l_j} = \max\{t_i: i\leq j\} = 2^{j-1}$ are unbounded, Theorem \ref{deg3infinite} implies that this $\tau$-sequence determines 
	$$\lim_{j\to\infty} 2^j/2^{j-1} = 2$$ 
solenoids in the moduli space $\cM_3$.

%%%%%%%
%%%%%%%

%%%%%%%
%%%%%%%

\bigskip\bigskip\section{Counting argument, all degrees}
\label{sec:counting}

In this final section, we give the proof of Theorem \ref{maintheorem1}.  We also provide an example of a pictograph associated to multiple topological conjugacy classes in any degree $d>2$.  Towards Theorem \ref{maintheorem1}, we have already established the topological-conjugacy invariance of the pictograph (Theorem \ref{spine invariance}).  Here, we show:

\begin{theorem} \label{conjugacy classes}
Let $\cD$ be a pictograph.  The number $\Top(\cD)$ of topological conjugacy classes of basins $(f, X(f))$ with pictograph $\cD$ is inductively computable from the discrete data of $\cD$.  Specifically, the computation depends only on the first-return map along the spine $(R, S(T))$ of the underlying tree and the automorphism group of the full tree of local models.  
\end{theorem}

\noindent
It is useful to compare this statement to those of Theorems \ref{deg3twist} and \ref{deg3infinite} containing the degree 3 computation.  In degree 3, the data of $(R, S(T))$ is equivalent to the Branner-Hubbard tableau and Yoccoz $\tau$-sequence.  Also in degree 3, the symmetries of the tree of local models are easy to describe.  Recall that the automorphism group is itself inductively computable from a pictograph in all degrees; see Section \ref{sec:symmetries}.

For the general degree case, we introduce the {\em restricted basin} of infinity for a polynomial, with a notion of equivalence that carries information from the full tree of local models.  This allows us to define an analog of the ``marked levels" in degree 3 (see \S\ref{marked levels}).  We introduce the {\em lattice of twist periods} to generalize the sequence of twist periods used to compute the number of conjugacy classes in degree 3.  In higher degrees, the markings and symmetries are significantly more complicated, so computing the number of twists needed to return to a given gluing configuration involves more ingredients.

\subsection{Restricted basins, conformal equivalence}  \label{restricted basins}
Let $(f, X(f))$ be a basin of infinity.  Fix any real number $t>0$.  Let $G_f$ denote the escape-rate function on $X(f)$, and set 
	$$X_t(f) = \{z\in X(f): G_f(z) > t\}.$$  
We refer to the pair $(f, X_t(f))$ as a {\em restricted basin}.  Here we introduce a special notion of equivalence of restricted basins that will be useful in the proof of Theorem \ref{conjugacy classes}.  

Let $(\cF_1, \cX(f_1))$ and $(\cF_2, \cX(f_2))$ denote the trees of local models for the basins $(f_1, X(f_1))$ and $(f_2, X(f_2))$, respectively.  By construction, there are gluing quotient maps 
	$$g_i: (\cF_i, \cX(f_i)) \to (f_i, X(f_i))$$
that are conformal isomorphisms from each local model surface to its image, inducing conjugacies between the restrictions of $\cF_i$ and $f_i$.  

We say the restricted basins $(f_1, X_t(f_1))$ and $(f_2, X_t(f_2))$ are {\em conformally equivalent over $(\cF, \cX)$} if their trees of local models are both isomorphic to $(\cF, \cX)$, and there exists a conformal isomorphism 	
	$$\phi: X_t(f_1) \to X_t(f_2)$$
inducing a conjugacy between the restrictions $f_1|X_t(f_1)$ and $f_2|X_t(f_2)$ that extends to the full tree of local models.  Specifically, there is an isomorphism between trees of local models
	$$\Phi:  (\cF_1, \cX(f_1)) \to (\cF_2, \cX(f_2))$$
which restricts to the induced isomorphism $\tilde{\phi}$ on the truncated trees, at heights $>t$, defined by lifting $\phi$ via the gluing quotient maps $g_i$:
$$\xymatrix{ (\cF_1, \cX_t(f_1)) \ar[r]^{\tilde{\phi}}    \ar[d]_{g_1} & (\cF_2, \cX_t(f_2))\ar[d]^{g_2} \\
			(f_1, X_t(f_1)) \ar[r]_{\phi}  &  (f_1, X_t(f_2))  }$$

Similarly, we define $\Aut_{(\cF, \cX)}(f, X_t(f))$ to be the group of conformal automorphisms of the restricted basin $(f, X_t(f))$ that extend to automorphisms of the tree $(\cF, \cX)$.  Denoting by $\Aut(f, X(f))$ and $\Aut(f, X_t(f))$ the groups of conformal isomorphisms (of $X(f)$ and $X_t(f)$, respectively) commuting with $f$, we find:

\begin{lemma}  \label{aut inclusions}
For any basin of infinity $(f, X(f))$ of degree $d\geq 2$, and any $t>0$, we have a chain of subgroups
	$$\Aut(f, X(f)) \subset \Aut_{(\cF, \cX)}(f, X_t(f)) \subset \Aut(f, X_t(f)) \subset C_{d-1},$$
where $C_{d-1}$ is the cyclic group of order $d-1$, acting by rotation in the uniformizing coordinates near $\infty$.
\end{lemma}

\proof
The first two inclusions follow easily from the definitions.  Indeed, any automorphism of a basin $(f, X(f))$ induces an automorphism of the tree of local models and of any restricted basin.  The last inclusion follows because an automorphism of $(f, X_t(f))$ must commute with $f$ near infinity, where it is conformally conjugate to $z^d$.  
\qed

\subsection{Restricted basins, topological equivalence}  \label{restricted, topological}
As for conformal equivalence of restricted basins, defined in \S\ref{restricted basins}, we say restricted basins $(f_1, X_t(f_1))$ and $(f_2,X_t(f_2))$ are {\em topologically equivalent over $(\cF, \cX)$} if there exists a topological conjugacy 
	$$\psi: X_t(f_1) \to X_t(f_2)$$
that extends to an isomorphism of the full tree of local models.  

It is important to observe that topologically conjugate restricted basins are also quasiconformally conjugate; the proof is identical to the one for full basins of infinity.  Further, if the restricted basins come from basins with the same critical escape rates, the quasiconformal conjugacy can be taken to be a twist deformation.  On each level set of $G_f$, the twist deformation acts by isometries (in the $|\del G_f|$ metric), and therefore it preserves the conformal structure of the local models in the tree of local models (compare the proof of Theorem \ref{tlm invariance}).  This proves:

\begin{lemma}
A topological conjugacy between restricted basins $(f_1, X_t(f_1))$ and $(f_2, X_t(f_2))$ with the same critical escape rates induces, via the gluing quotient maps, an isomorphism of truncated trees of local models $(\cF_1, \cX_t(f_1))$ and $(\cF_2, \cX_t(f_2))$.  
\end{lemma}

\noindent
In the definition of topological equivalence of restricted basins, then, we are requiring that this induced isomorphism of truncated trees of local models can be extended to an isomorphism on the full tree of local models.

\subsection{Twist periods}  \label{periods}
As described in \S\ref{qc}, a quasiconformal deformation of a basin of infinity has a canonical decomposition into its twisting and stretching factors.  Fix $f\in \cM_d$ and consider the analytic map of \S\ref{twist normalization},
	$$\Tw_f: \R^N \to \cM_d,$$
which parametrizes the twisting deformations in the $N$ fundamental subannuli of $f$, sending the origin to $f$.   Recall that the basis vector
	$${\bf e}_j = (0, \ldots, 0, 1, 0, \ldots, 0)\in \R^N$$
induces a full twist in the $j$-th fundamental subannulus.  

A {\em twist period} of $f$ is any vector $\tau\in\R^N$ which preserves the conformal conjugacy class of $(f, X(f))$; that is, $\Tw_f(\tau) = \Tw_f(0)$.  When $f$ is in the shift locus, the set of twist periods forms a lattice in $\R^N$ \cite[Lemma 5.2]{DP:heights}.  In general, the set of twist periods forms a discrete subgroup 
	$$\TP(f) \subset \R^N.$$ 
As we shall see in the proof of Theorem \ref{conjugacy classes}, polynomials with equivalent pictographs can have different lattices of twist periods; this can happen when one gluing configuration has automorphisms while another does not.  Nevertheless, we will see that the possibilities for $\TP(f)$ can still be computed from the data of the pictograph.

For each $t>0$, we define $\TP_t(f) \supset \TP(f)$ to be the set of vectors $\tau\in \R^N$ that preserve the conformal equivalence class of the restricted basin $(f, X_t(f))$ over $(\cF, \cX(f))$; the equivalence of restricted basins was defined in \S\ref{restricted basins}.

\begin{lemma}  \label{twist period lattice}
For any $t_1 > t_2 >0$ and any $f$ with $N$ fundamental subannuli, each group $\TP_{t_i}(f)$ forms a lattice in $\R^N$, with index $[\TP_{t_1}(f): \TP_{t_2}(f)] < \infty$ and 
	$$\TP(f) = \bigcap_{t>0} \TP_t(f).$$  
\end{lemma}
	
\proof
The argument is similar to the proof of \cite[Lemma 5.2]{DP:heights}.  Let $X_t(f) = \{G_f > t\} \subset X(f)$.  For each fundamental subannulus $A_j$, there are only finitely many connected components $B_j$ of preimages of $A_j$ inside $X_t(f)$ under any iterate of $f$.  A full twist in the annulus $A_j$ induces a $1/k$-twist in a preimage $B_j$, where $k$ is the local degree of the iterate $f^n$ sending $B_j$ to $A_j$.  Let 
	$$d_j = \mathrm{lcm} \{k : k = \deg (f^n|B_j \to A_j )\}$$
over all such components $B_j \subset X_t(f)$.  Then the subgroup $\TP_t(f)$ of $\R^N$ must contain the vector
	$$d_j \, {\bf e}_j = (0, \ldots, 0, d_j, 0, \ldots, 0)$$
for each $j$; indeed, this vector induces an automorphism of $(f, X_t(f))$ that extends to the identity automorphism on the full tree of local models.  Because $\TP_t(f)$ is a discrete subgroup of $\R^N$, we see that it must be a lattice.  The same argument also shows that $[\TP_{t_1}(f): \TP_{t_2}(f)] < \infty$.

Finally, if $\tau\in \TP(f)$, then $\tau$ must induce an equivalence of the restricted basin $(f, X_t(f))$ over $(\cF, \cX(f))$ for every $t>0$, because a conjugacy between basins induces an equivalence on trees of local models.  Therefore, 
	$$\TP(f) \subset \bigcap_{t>0} \TP_t(f).$$
Conversely, we observe that if $\tau\cdot (f, X_t(f))$ is conformally equivalent to $(f, X_t(f))$ for all $t>0$, then there is a conformal conjugacy between the basins $\tau\cdot (f, X(f))$ and $(f, X(f))$, so $\tau\in \TP(f)$.  We conclude that $\TP(f) = \bigcap_{t>0} \TP_t(f)$.
\qed

\subsection{Twist periods in degree 3}Ê \label{degree 3 twist periods}
We remark that the definition of twist period given here differs from that given in 
\S\ref{counting cubics} for degree 3 maps.Ê They coincide in the case of one fundamental 
annulus, in the sense that the twist periods $\{T_n\}$ form a sequence of generators for 
the one-dimensional twist lattices $\TP_{t_n}(f) \subset \R$ for heights $t_n$ just below 
level $n$.

A cubic polynomial $f$ with two fundamental subannuli is necessarily in the shift locus 
(and structurally stable) with finite length $L(f) > 0$. We describe here how to recover 
the sequence of twist periods $T_n$ at levels $0 \leq n < L(f)$ from the lattices of 
twist periods $\TP_t(f)$.

For a cubic polynomial $f$ with two fundamental subannuli, let $(F, T(f))$ denote its 
tree, let $v_0$ be the highest branching vertex of $T(f)$, and let $h: T(f) \to \R$ be 
the height function.Ê Choose a sequence of descending heights
\begin{equation} \label{cubic t_n}
Ê Ê Ê Ê t_0 > t_1 > t_2 > \cdots
\end{equation}
so that $t_n$ is a height ``just below" a vertex of combinatorial distance $n$ from 
$v_0$.Ê That is, $t_0 = h(v_0) - \eps$ for any sufficiently small $\eps>0$, and there is 
a unique vertex in each connected component of $h^{-1}(t_{n+1}, t_n)$.

From the definitions, and the absence of symmetries at $v_0$, we have $T_0 = 1$ and
Ê Ê Ê Ê $$\TP_{t_0} = \; \< {\bf e}_1, \, {\bf e}_2 \> \; = \Z^2 \subset \R^2$$
for every such polynomial.

Let $S(T)$ the spine of the tree $(F, T(f))$.Ê Denote by $w$ the lowest vertex in the 
spine (the lower critical point).Ê For each positive integer $n < L(f)$, set
Ê Ê Ê Ê $$J(n) = \# \{L(f)-n \leq j < L(f) : F^j(w) \in S(T) \}, $$
the number of times the critical orbit intersects the spine, above level $n$ and below 
$v_0$.Ê From the proof of Theorem \ref{conjugacy classes} given below, an inductive 
argument shows that
\begin{equation} \label{cubic twists}
Ê Ê Ê Ê T_n = \frac{[\TP_{t_0}: \TP_{t_{2n}}]}{2^{J(n)}}.
\end{equation}

\subsection{A cubic example}Ê \label{cubic twist example}
Because the computation of the twist period lattices is crucial in the proof of Theorem 
\ref{conjugacy classes}, we illustrate with an example in degree 3.

\begin{figure}
\includegraphics[width=1.7in]{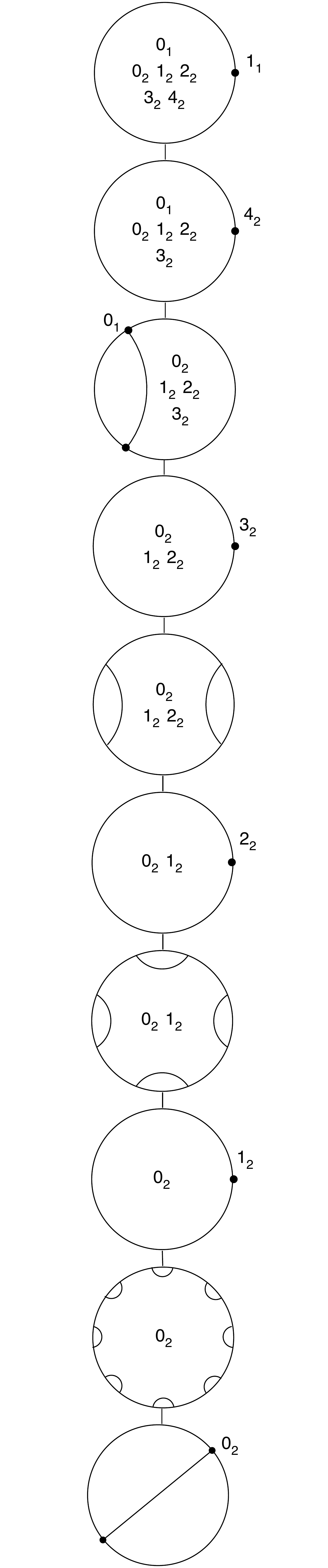}
\caption{A cubic pictograph with $\tau$ sequence $(0,1,2,3)$.Ê Its twist periods are 
computed in \S\ref{cubic twist example}.} \label{0123}
\end{figure}

Consider the cubic pictograph shown in Figure \ref{0123}.Ê It is the unique pictograph 
associated to the $\tau$-sequence $(0,1,2,3)$ with two fundamental edges.Ê Any associated 
polynomial has length $L(f) = 4$. There are no marked levels, in the sense defined in 
\S\ref{marked levels}. Consequently, $T_0 = T_1 = T_2 = T_3 = 1$.Ê From Theorem 
\ref{deg3twist}, there is a unique topological conjugacy class of cubic polynomials with 
this pictograph.

Now let $\{t_n\}$ be a descending sequence of escape rates, as defined in (\ref{cubic 
t_n}).Ê Every polynomial $f$ with this pictograph has the following lattices of twist 
periods:
\begin{eqnarray*}
Ê Ê Ê Ê \TP_{t_0} &=&Ê \; \< {\bf e}_1, \, {\bf e}_2 \>Ê \\
Ê Ê Ê Ê \TP_{t_1}\; = \; \TP_{t_2} &= & \; \< {\bf e}_1, \, 2{\bf e}_2 \> \\
Ê Ê Ê Ê \TP_{t_3}\; = \; \TP_{t_4} &=& \; \< 2{\bf e}_1, \, {\bf e}_1 + 2{\bf e}_2 \>Ê \\
Ê Ê Ê Ê \TP_{t_n}Ê &=& \; \< 4{\bf e}_1, \, 3{\bf e}_1 + 2{\bf e}_2 \> \mbox{ for all $n\geq 5$}
\end{eqnarray*}
These lattices are computed inductively with $n$, with the base case $\TP_{t_0} = \< {\bf 
e}_1, \, {\bf e}_2 \>$.Ê To determine $\TP_{t_n}$ from $\TP_{t_{n-1}}$, we compute the 
induced twisting at each vertex of the spine, down to height $t_n$, for elements of 
$\TP_{t_{n-1}}$. Specifically, let $v$ be the lowest vertex in the spine above height 
$t_n$, and fix $\tau = \tau_1{\bf e}_1 + \tau_2{\bf e}_2 \in\TP_{t_{n-1}}$. We compute 
(1) the relative modulus of $\mod(e)$ of each edge $e$ in the spine (as the reciprocal of 
the degree by which $e$ maps to a fundamental edge) down to the vertex $v$, and (2) the 
sum
Ê Ê Ê Ê $$R_\tau(v) = \sum_e \tau_{j(e)} \mod(e)$$
where $j(e)$ is the index of the fundamental edge in the orbit of $e$.Ê If the  
lamination diagram at vertex $v$ is invariant under rotation by the amount $R_\tau(v)$ (leaving also all labels invariant),
then $\tau\in \TP_{t_n}$.Ê If not, we consider integer multiples of $\tau$.

Finally, $\TP(f) = \< 4{\bf e}_1, \, 3{\bf e}_1 + 2{\bf e}_2 \>$ by Lemma \ref{twist 
period lattice}.Ê Observing that $J(n) = n$ for each $n = 1, 2, 3$, we also see that 
equation (\ref{cubic twists}) holds.
  
\subsection{Proof of Theorem \ref{conjugacy classes}}
Fix a pictograph $\cD$.  We aim to show that the number $\Top(\cD)$ of topological conjugacy classes of basins $(f, X(f))$ with pictograph $\cD$ can be inductively computed from the discrete data of $\cD$.  

For any $M>0$, there exists a vector of compatible critical heights with maximal critical height $M$.  Let  $(\cF, \cX)$ be the tree of local models with pictograph $\cD$ and the chosen vector of critical heights given by Proposition \ref{spine}.  Let $(F, T)$ be the underlying polynomial tree with height function $h: T\to (0,\infty)$, so that $h(v_0) = M$ for the highest branching vertex $v_0$ and $h(Fx) = d\cdot h(x)$ for all $x\in T$.  Choose a descending sequence of real numbers 
	$$t_{-1} > M > t_0 > t_1 > t_2 > \cdots > 0$$ 
so that 
\begin{itemize}
\item	the points of $h^{-1}(t_i)$ are not vertices of $T$ for any $i\geq 0$, and 
\item	there is a unique vertex of $T$ in each connected component of $h^{-1}[t_i, t_{i-1}]$ for all $i \geq 0$.
\end{itemize}

For each $i \geq 0$, we will inductively compute the number $\Top(\cD, i)$ of topological conjugacy classes of restricted basins 
	$$f : X_{t_i}(f) \to X_{t_i}(f)$$
over $(\cF, \cX)$; the equivalence was defined in \S\ref{restricted, topological}.

The number $\Top(\cD)$ is not simply the limit of these numbers $\Top(\cD,i)$ as $i\to \infty$, as we shall see, but it can be determined from the sequence $\{\Top(\cD,i)\}_i$ and the data used to compute it.  

Let $i=0$.  The conformal equivalence class of the restricted basin $(f, X_{t_0}(f))$ over $(\cF, \cX)$ depends only on the gluing along each of the fundamental subannuli.  It is easy to see, from the definition of twisting, that all gluing choices within the fundamental annulus are equivalent under twisting.  We conclude that 
	$$\Top(\cD, 0) = 1$$
for any pictograph $\cD$.  

For the induction argument, we need to compute the lattice of twist periods $\TP_{t_0}(f)$ and the automorphism group $\Aut_{(\cF, \cX)}(f, X_{t_0}(f))$ for any choice of restricted basin $(f, X_{t_0}(f))$ from the discrete data of $\cD$ at and above the vertex $v_0$ and automorphism group $\Aut(\cF, \cX)$.  Recall that the automorphism group of the full tree of local models is isomorphic to the automorphism group of the first-return map $(\cR, \cS)$ on the spine (Lemma \ref{spine automorphisms}), so any information we need about $\Aut(\cF, \cX)$ is determined by $\cD$.  

Suppose $\cD$ has $N$ fundamental  edges (so any basin with pictograph $\cD$ has $N$ fundamental subannuli).  As usual, we label the ascending consecutive vertices $v_0, v_1, \ldots, v_N = F(v_0), v_{N+1}, \ldots$ in the tree $(F, T)$, where $v_0$ is the highest branching vertex.  In \S\ref{local symmetry}, we defined the order of local symmetry of the tree $(\cF, \cX)$ at the vertex $v_j$; we denote this order by $k_j$.  

\begin{lemma}  \label{ascending symmetry}
For each $j\geq 0$, the orders of symmetry at $v_j$ and $v_{j+N}$ satisfy
	$$k_{j+N} = k_j/ \gcd(k_j, d).$$
\end{lemma}

\proof
The local degree of $F$ at each vertex $v_j$ is $d$.  From Lemma \ref{symmetry orders}, we know that $k_j/\gcd(k_j, d)$ must divide $k_{j+N}$.  On the other hand, by the definition of the local symmetry order (coming from an automorphism of $(\cF, \cX)$), any symmetry at $v_{j+N}$ must lift to the domain $v_j$.  Therefore, we have equality.  
\qed

\begin{lemma} \label{restricted aut group}
The automorphism group $\Aut_{(\cF, \cX)}(f, X_{t_0}(f))$ is cyclic of order equal to
	$$\alpha = \gcd\{ k_0, k_1, \ldots, k_{N-1}, d-1\}.$$
\end{lemma}

\proof
From Lemma \ref{ascending symmetry}, the value $\alpha$ will divide the orders of local symmetry at every vertex $v_j$, $j\geq 0$.  Further, a rotation by $2\pi/\alpha$ at any vertex $v_j$ in the tree $(\cF, \cX)$ will induce a rotation by $2\pi d/\alpha \equiv 2\pi/\alpha \mod 2\pi$ at its image, because $\alpha | (d-1)$.  Therefore, rotation by $2\pi/\alpha$ can act on any gluing of $(\cF, \cX)$ to form a restricted basin $(f, X_{t_0}(f))$.  The automorphism must extend to the full tree of local models, by the definition of the orders $k_j$.  It follows that the order of  $\Aut_{(\cF, \cX)}(f, X_{t_0}(f))$ is at least $\alpha$.  On the other hand, the order of any element in  $\Aut_{(\cF, \cX)}(f, X_{t_0}(f))$ must divide $\alpha$, combining the definitions with Lemma \ref{symmetry orders}.
\qed

\medskip
Fix a conformal equivalence class of restricted basin $(f, X_{t_0}(f))$ over $(\cF, \cX)$.  To compute $\TP_{t_0}(f)$, note first that each of the basis vectors 
	 $${\bf e}_j = (0, \ldots, 0, 1, 0, \ldots, 0)$$
for $j=1, \ldots, N$, are contained in $\TP_{t_0}(f)$, by construction.  Indeed, a full twist in any subannulus induces the identity automorphism on the tree of local models $(\cF, \cX)$.  

We claim that for each $0 < j < N$, the twist vector	$$\tau_j = \frac{1}{k_j} \, ({\bf e}_{j+1} - {\bf e}_j)$$
is also contained in $\TP_{t_0}(f)$.  This vector $\tau_j$ twists by $1/k_j$ in the fundamental subannulus $A_{j+1}$ and by $-1/k_j$ in $A_j$.  It therefore twists by the order of symmetry at the vertex $v_j$, and it induces a twist by $d/k_j$ in the image of $A_{j+1}$ and by $-d/k_j$ in the image of $A_j$.  By Lemma \ref{ascending symmetry}, these twists commute with the action of $f$.  The restricted basins $(f, X_{t_0}(f))$ and $\tau_j\cdot (f, X_{t_0}(f))$ are conformally conjugate; for $k_j>1$, the isomorphism extends to a non-trivial isomorphism of the underlying tree of local models where the action on $(\cF, \cX)$ rotates vertices in the grand orbit of $v_j$.

Finally, we treat the symmetry at $v_0$.  Set
	$$\tau_0 = \frac{1}{k_0} \, ({\bf e}_1 - d \, {\bf e}_N).$$
As for $\tau_j$, $j>0$, the twist vector $\tau_0$ induces the symmetry at $v_0$.  The term $(d/k_0) \, {\bf e}_N$ is subtracted off so that the correct order of symmetry is induced at the image $v_N$, as in Lemma \ref{ascending symmetry}.  Putting the pieces together, we find that $\TP_{t_0}(f)$ is generated by all twist vectors of the form 
	$$\tau = (a_1, \ldots, a_N),$$
where $0 \leq a_j \leq 1$, $\sum_{j = j_0}^N a_j$ is an integer multiple of $1/k_{j_0-1}$ for each $j_0>1$, and $\sum_{j=1}^N a_j$ is an integer multiple of $(1-d)/k_0 \mod 1$.  Observe that this computation is independent of the initial choice of restricted basin $(f, X_{t_0}(f))$.

Now fix $i \geq 0$ and a restricted basin $(f, X_{t_0}(f))$.  Let $\cB_i(\cD)$ denote the set of conformal equivalence classes of restricted basins $(f, X_{t_i}(f))$ over $(\cF, \cX)$ that extend the restricted basin $(f, X_{t_0}(f))$.  Suppose we have computed 
\begin{enumerate}
\item	the number of conformal equivalence classes $|\cB_i(\cD)|$; 
\item the order of the automorphism group $\Aut_{(\cF, \cX)}(f, X_{t_i}(f))$ for each element of $\cB_i(\cD)$; and
\item	the lattice of twist periods $\TP_{t_i}(f)$ for each element of $\cB_i(\cD)$.
\end{enumerate}
As explained above, the conformal classes in $\cB_i(\cD)$ are topologically equivalent if and only if they are equivalent by twisting, via a conjugacy that extends to an isomorphism of the full tree of local models, so we need only compute the number of classes in each twist orbit to obtain $\Top(\cD, i)$ from this data.  That is,
\begin{equation} \label{top computation} 
\Top(\cD, i) = \sum_{(f, X_{t_i}(f)) \, \in \, \cB_i(\cD)} \frac{1}{ [\TP_{t_0}(f): \TP_{t_i}(f)] } \, .
\end{equation}

Now we pass to $i+1$.  Let $\{(X_v, \omega_v)\}$ be the set of local models in the spine of $(\cF, \cX)$ with vertex $v$ in the height interval $(t_{i+1}, t_i)$.  Let $d_v$ be the degree of the local model map with domain $(X_v, \omega_v)$.  Let $k_v$ be the order of local symmetry of $(\cF, \cX)$ at $v$.

Fix a conformal class $(f, X_{t_i}(f)) \in \cB_i(\cD)$, and assume that $\Aut_{(\cF, \cX)}(f, X_{t_i}(f))$ is the trivial group.  Then the number of classes in $\cB_{i+1}(\cD)$ that extend $(f, X_{t_i}(f))$ is given by 
	$$\prod_v \frac{d_v}{\gcd( k_v, d_v) } ,$$
where the product is taken over all vertices $v$ of the spine in the height interval $(t_{i+1}, t_i)$.  Indeed, the extension to height $t_{i+1}$ along any edge of degree 1 is uniquely determined.  We need only compute how many distinct ways we may glue each local model $(X_v, \omega_v)$ of degree $d_v>1$ along the edge above $v$ so that $f$ extends holomorphically.  The absence of automorphisms shows that the local symmetry (fixing $v$) and local degree are the only contributing factors.  It is easy to see that each extension will also have a trivial automorphism group.  

Now suppose $(f, X_{t_i}(f)) \in \cB_i(\cD)$ has automorphism group $\Aut_{(\cF, \cX)}(f, X_{t_i}(f))$ of order $m>1$.  By Lemma \ref{aut inclusions}, the automorphism group is cyclic, acting by rotation in the uniformizing coordinates near infinity.  By construction, every such automorphism extends to the full tree of local models, so there is a certain amount of symmetry among the vertices in the height interval $(t_{i+1}, t_i)$.  First, there is at most one vertex $v'$ in the spine at this height fixed by the automorphism of order $m$, and $m$ must divide the local symmetry $k_{v'}$.  All other vertices $v$ of the spine have orbit of length $m$, and the order of local symmetry $k_v$ is constant along an orbit.  Choose a representative vertex $\hat{v}$ for each orbit.  The number of conformal classes in $\cB_{i+1}(\cD)$ extending $(f, X_{t_i}(f))$ are organized as follows.  There are
	$$N(m) = \frac{d_{v'}}{\gcd(k_{v'}, d_{v'})} \cdot \prod_{\tiny\mbox{orbits of length } m} \frac{d_{\hat{v}}}{\gcd(k_{\hat{v}}, d_{\hat{v}}) }$$
conformal conjugacy classes of extensions with an automorphism of order $m$; indeed, a choice of gluing at vertex $\hat{v}$ determines the choice (up to local symmetry) at each vertex in its orbit.  For each factor $l |m$, we can also compute the number of extensions of $(f, X_{t_i}(f))$ with automorphism of order $l$.  A simple inclusion-exclusion argument shows that there are 
	$$N(l) = \frac{d_{v'}}{\gcd(k_{v'}, d_{v'})} \cdot \prod_{\tiny{\mbox{orbits of length }} l} \frac{d_{\hat{v}}}{\gcd(k_{\hat{v}}, d_{\hat{v}}) } \; - \sum_{\{l'\, : \; l|l'|m,\, l'\not=l\}} N(l')$$ 
conformal equivalence classes of extensions with an automorphism of order $l$, under the extra restriction that we require the equivalence to act by the identity on the restriction $(f, X_{t_i}(f))$.  Consequently, there are 
	$$\frac{l}{m} \cdot N(l)$$
distinct conformal equivalence classes of extensions with automorphism of order $l$ (without the extra assumption): for each new basin in the initial count of $N(l)$, there is an isomorphism that acts as rotation by $l/m$ near infinity, producing another gluing configuration in the count of $N(l)$.  We observe that the computation depends only on the data $\{d_v, k_v\}$ at each of the vertices in the spine.

Now we compute the lattice of twist periods $\TP_{t_{i+1}}(f)$ for each class $(f, X_{t_{i+1}}(f)) \in \cB_{i+1}(\cD)$.  Note that any element $\tau\in \TP_{t_{i+1}}(f)$ also induces a conformal equivalence of the restricted basin $(f, X_{t_i}(f))$ over $(\cF, \cX)$, so $\TP_{t_{i+1}}(f)$ forms a subgroup of $\TP_{t_i}(f)$.   We will examine linear combinations of basis elements of $\TP_{t_i}(f)$ to determine which elements lie in $\TP_{t_{i+1}}(f)$.  We need to use the order of the automorphism group of the chosen basin $(f, X_{t_{i+1}})$ and the order of the automorphism group of the restriction $(f, X_{t_i}(f))$.  

First assume that both $(f, X_{t_{i+1}}(f))$ and the restriction $(f, X_{t_i}(f))$ have trivial automorphism group.  Fix any element $\tau\in \TP_{t_i}(f)$, so that $\tau\cdot(f, X_{t_i}(f))$ is conformally equivalent to $(f, X_{t_i}(f))$ over $(\cF, \cX)$.  One can check algorithmically whether a multiple $a\tau$ lies in $\TP_{t_{i+1}}(f)$, for each $a = 1, 2, 3, \ldots$, via the following steps:
\begin{enumerate}
\item		Compute all relative moduli, down to the vertex $v$, for each $v$ in the spine in the height interval $(t_{i+1}, t_i)$:  each edge $e$ between $v$ and $v_0$ is mapped by a degree $d(e) > 1$ to one of the fundamental edges, and its relative modulus is $1/d(e)$.
\item		Compute the rotation induced by $a\tau$ at each vertex $v$:  if $a\tau = (t_1, \ldots, t_N)$, then $v$ is rotated by 
		$$R_{a\tau}(v) = \sum_e \frac{t_{j(e)}}{d(e)},$$
where the sum is over all edges $e$ between $v$ and $v_0$, and $j(e)$ is the index of the unique fundamental edge in the orbit of $e$.
\item		If the rotation $R_{a\tau}(v)$ is an integer multiple of $1/k_v$ at each vertex $v$, then $a\tau \in \TP_{t_{i+1}}(f)$.
\end{enumerate}
Note that this process terminates at some finite value of $a$:  by Lemma \ref{twist period lattice}, we know that $\TP_{t_{i+1}}$ is a lattice of finite index within $\TP_{t_i}$.  With trivial automorphism group, the computation of $\TP_{t_{i+1}}(f)$ is independent of the choice of extension $(f, X_{t_{i+1}}(f))$.  

We remark that, even in the absence of global automorphisms, there can be local symmetries of $(\cF, \cX)$ that act nontrivially on the spine of the underlying tree.  A given twist vector $\tau$ may induce one of these nontrivial automorphisms of the spine; so the computation of steps (2) and (3) above requires that we compare action on $v$ to symmetry at a different vertex, say $v'$.  In that case, the local symmetries $k_v$ and $k_{v'}$ will coincide, so the computation is the same.  

Suppose now that our chosen $(f, X_{t_{i+1}}(f))$ has automorphism group of order 
	$$|\Aut_{(\cF, \cX)}(f, X_{t_{i+1}}(f))| = l$$ 
and the restriction $(f, X_{t_i}(f))$ has 
	$$|\Aut_{(\cF, \cX)}(f, X_{t_i}(f))| = m \geq l.$$
Note that $l| m$.  Fix any element $\tau\in \TP_{t_i}(f)$.  For $l=m$, we may proceed as above:  we check integer multiples of $\tau$ and compute the rotation induced at each of the lowest vertices $v$.  For $l < m$, we have an additional possibility.  It can happen that the twisted basin $a\tau \cdot (f, X_{t_{i+1}})$ is conformally equivalent to the basin $(f, X_{t_{i+1}})$ via an isomorphism that acts as rotation by $k/m$ near infinity, for some integer $1 \leq k < m/l$.  Thus our algorithmic procedure involves an extra computation.  The three steps above become:
\begin{enumerate}
\item		Compute all relative moduli, down to the vertex $v$, for each $v$ in the spine in the height interval $(t_{i+1}, t_i)$, as before.
\item		Compute the rotation induced by $\phi^k\circ (a\tau)$ at each vertex $v$, where $\phi^k$ acts as rotation by $k/m$ near infinity, for each $k = 1, \ldots, m/l$:  it is given by the simple relation $R_{\phi^k\circ (a\tau)}(v) = R_{a\tau}(v) + (k/m)$.
\item		If for any $k$, the rotation $R_{\phi^k\circ (a\tau)}(v)$ is an integer multiple of $1/k_v$ at each vertex $v$, then $a\tau \in \TP_{t_{i+1}}(f)$.
\end{enumerate}
We illustrate with one example in degree 5 the delicacy of computing twist periods in the presence of automorphisms; see \S\ref{symmetry example} and Figure \ref{degree 5 aut}.  

To make the above algorithmic process implementable, it is useful to compute an explicit basis for $\TP_{t_{i+1}}(f)$.  For example, let $\tau_1, \ldots, \tau_N$ be a set of basis vectors for the lattice $\TP_{t_i}(f)$.  We can apply the above steps to each basis vector.  Let $a_j$ is the smallest positive integer so that $a_j \tau_j\in\TP_{t_{i+1}}(f)$.  We next compute the rotation effect of each vector of the form
	$$n_1 \tau_1 + \cdots + n_N \tau_N$$
for all tuples of non-negative integers $\{n_i\}$ with $n_i \leq a_i$.  This is a finite process and will produce a basis for $\TP_{t_{i+1}}(f)$.   

Once we have computed the twist periods for each class $(f, X_{t_{i+1}}(f)) \in \cB_{i+1}(\cD)$, the number $\Top(\cD, i+1)$ is computed by equation (\ref{top computation}).

Finally, we need to compute $\Top(\cD)$.  Note that the number of conformal classes extending a given $(f, X_{t_0}(f))$ is non-decreasing with $i$; that is, 
	$$|\cB_i(\cD)| \leq |\cB_{i+1}(\cD)|.$$
We claim
\begin{enumerate}
\item	$\Top(\cD) = \infty$ if and only if $\lim_{i\to\infty} |\cB_i(\cD)| = \infty$; and 
\item	if $\lim_{i\to\infty} |\cB_i(\cD)| = |\cB_n(\cD)|$ for some $n$, then $\Top(\cD) = \Top(\cD, n)$.
\end{enumerate}
From Theorem \ref{thm:bundle_of_gluing}, the number $\Top(\cD)$ is bounded above by the number of points in the fiber of the bundle of gluing configurations.  By the construction of the bundle (from the proof of Theorem \ref{thm:bundle_of_gluing}) the number of points in a fiber is equal to an integer multiple of $\lim_{i\to\infty} |\cB_i(\cD)|$.  Therefore, if $\Top(\cD) = \infty$, then it must be that $\lim_{i\to\infty} |\cB_i(\cD)| = \infty$.  On the other hand, when the fiber of the bundle of gluing configurations has infinite cardinality, Lemma \ref{Cantor fiber} states that the fibers are Cantor sets.  In particular, the fibers are uncountable.  A topological conjugacy class within the fiber contains at most countably many elements, as the image of a lattice in $\R^N$.  Therefore, there are infinitely many topological conjugacy classes.  This proves statement (1).  The second statement is immediate from the arguments and definitions above;  once the number of classes $|\cB_i(\cD)|$ has stabilized, the lattice of twist periods $\TP_{t_i}(f)$ also stabilizes.
\qed

\subsection{Example:  $\Top(\cD) = 1$ for the pictograph in Figure \ref{abstract spine}}  
Following the steps in the proof of Theorem \ref{conjugacy classes}, we compute the number of topological conjugacy classes associated to the pictograph constructed in Figure \ref{abstract spine}, with two fundamental subannuli.  Fix any restricted basin $(f, X_{t_0}(f))$ with the given pictograph.  It is easy to see from the diagram that the automorphism group is trivial, as are the local symmetries at the two fundamental vertices.  The lattice of twist periods $\TP_{t_0}$ is generated by the standard basis vectors, 
	$$\TP_{t_0} = \; \< {\bf e}_1, \, {\bf e}_2 \> \;\subset \R^2.$$   
There are two vertices $v_1$ and $w_1$ in the height interval $(t_1, t_0)$, with local degrees $d_{v_1} = 2$ and $d_{w_1} = 3$.  Each has trivial local symmetry, so we compute that 
	$$|\cB_1(\cD)| = 2\cdot 3 = 6,$$
and each class has trivial automorphism group.  A full twist in fundamental subannulus $A_1$ leaves $v_1$ and $w_1$ invariant, but a full twist in subannulus $A_2$ induces a $1/2$ twist at $v_1$ and a $1/3$ twist at $w_1$.  It requires 6 twists in $A_2$ to return to the given gluing configuration at level 1.  We compute,
	$$\TP_{t_1} = \; \< {\bf e}_1, \, 6{\bf e}_2 \>.$$
The computation is independent of the conformal class in $\cB_1(\cD)$.   We find that 
	$$\Top(\cD, 1) = 6\cdot \frac{1}{6} = 1.$$
In the height interval $(t_2, t_1)$, there are again two vertices, say $v_2$ below $v_1$ and $w_2$ below $w_1$.  We have local degrees $d_{v_2} = 2$, $d_{w_2} = 2$ and local symmetries $k_{v_2} = 2$, $k_{w_2} = 1$.  Therefore, 
	$$|\cB_2(\cD)| = |\cB_1(\cD)| \cdot 1\cdot 2 = 12,$$
and each class has trivial automorphism group.  A full twist in fundamental subannulus $A_1$ induces a $1/2$ twist at $v_2$, and $1/2$ is an integer multiple of $1/k_{v_2}$.  On the other hand, a $1/2$ twist is also induced at $w_2$ with $k_{w_2} = 1$, so we find that ${\bf e}_1 \not\in \TP_{t_2}$ but $2 {\bf e}_2 \in \TP_{t_2}$.  For the subannulus $A_2$, the action of $6{\bf e}_2$ induces full rotations of both $v_2$ and $w_2$, so $6{\bf e}_2 \in \TP_{t_2}$.    We find that
	$$\TP_{t_2} =  \; \< 2{\bf e}_1, \, 6{\bf e}_2 \>$$
and 
	$$\Top(\cD, 2) = 12 \cdot \frac{1}{12} = 1.$$
For all vertices in the spine below $v_2$ and $w_2$, the local symmetry at a vertex coincides with the local degree.  Therefore 
	$$|\cB_i(\cD)| = |\cB_2(\cD)| = 12$$
for all $i\geq 2$.  As explained at the end of the proof of Theorem \ref{conjugacy classes}, the number of topological conjugacy classes also stabilizes, so we may conclude that the pictograph determines
	$$\Top(\cD) = \Top(\cD, 2) = 1$$
topological conjugacy class of basins $(f, X(f)) \in \cB_5$.  In fact, because these polynomials are in the shift locus, this pictograph determines a unique topological conjugacy class of polynomials in $\cM_5$.  
	
\subsection{A degree 5 example with symmetry}  \label{symmetry example}
\begin{figure} 
\includegraphics[width=5in]{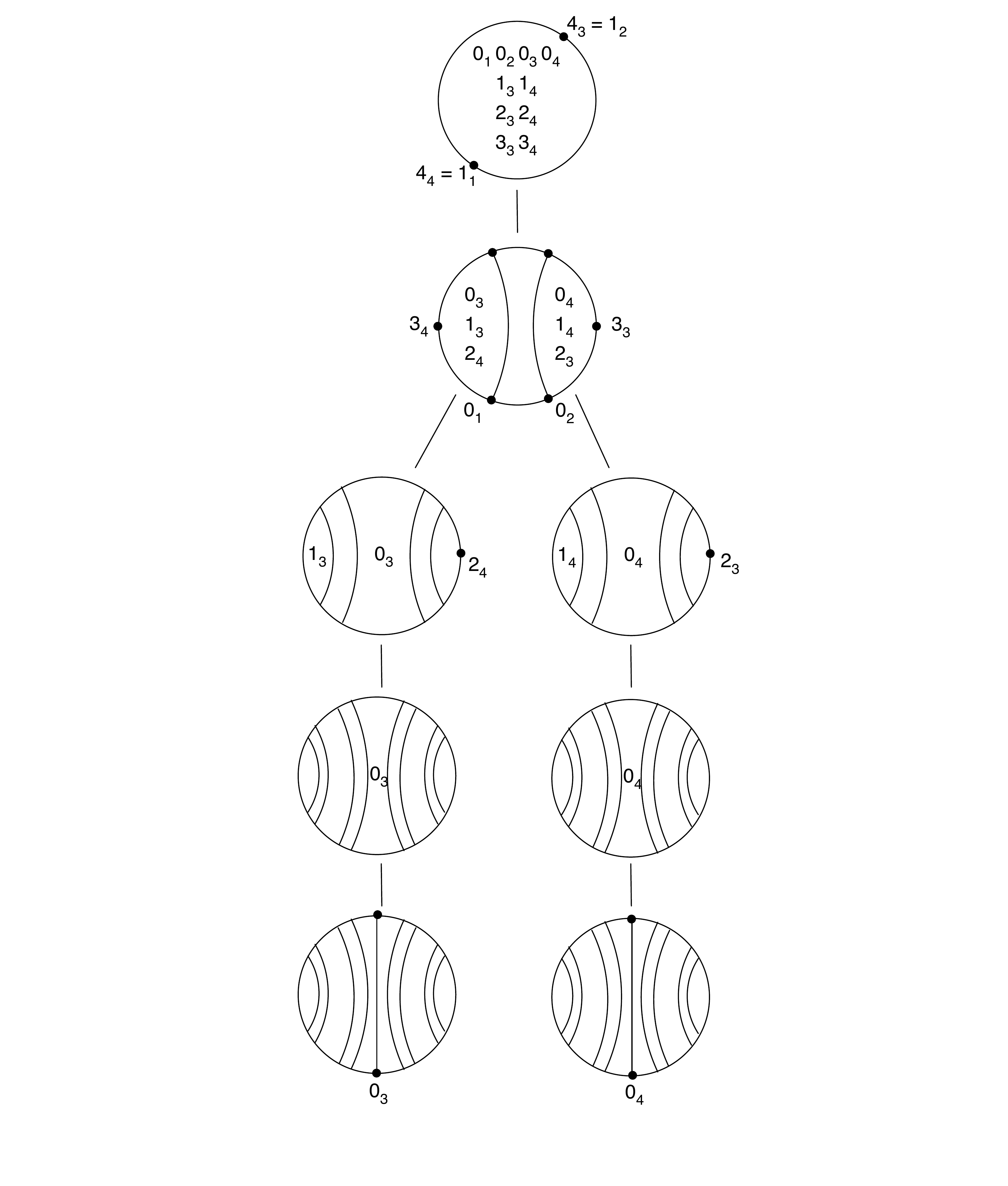}
\caption{A degree 5 pictograph with symmetry.  See \S\ref{symmetry example}.} \label{degree 5 aut}
\end{figure}

Consider the pictograph of Figure \ref{degree 5 aut}.  It has one fundamental edge.  Fix the critical escape rate $M>0$ of the highest critical points $c_1$ and $c_2$, and choose heights 
	$$M> t_0 > M/5 > t_1 > M/25$$ 
as in the proof of Theorem \ref{conjugacy classes}.  

Any choice of restricted basin $(f, X_{t_0}(f))$ has  $\Aut_{(\cF, \cX)}(f, X_{t_0}(f))$ of order 2, interchanging the critical points labelled by $0_1$ and $0_2$.  The extension to the full tree $(\cF, \cX)$ also interchanges the critical points labelled by $0_3$ and $0_4$.  We have 
	$$\TP_{t_0} = \Z.$$
Fixing $(f, X_{t_0}(f))$, there are three conformal equivalence classes of extensions over $(\cF, \cX)$ to level $t_1$.  Two of these extensions, say $(f_1, X_{t_1})$ and $(f_2, X_{t_1})$, will have an automorphism of order 2.  The third $(f_3, X_{t_1})$ has trivial automorphism group.  The restricted basins $(f_1, X_{t_1})$ and $(f_2, X_{t_1})$ are in the same topological conjugacy class over $(\cF, \cX)$, as one full twist in the fundamental annulus interchanges them; we have 
	$$\TP_{t_1}(f_1) = \TP_{t_1}(f_2) = 2\Z.$$
In the conformal class without automorphisms, one full twist arrives at a basin that is conformally equivalent via an isomorphism that rotates the basin by 180 degrees, and 
	$$\TP_{t_1}(f_3) = \Z.$$  
These restricted basins have unique conformal extensions to basins $(f_i, X(f_i))$, $i = 1, 2, 3$.  This pictograph determines exactly two topological conjugacy classes of polynomials, one with automorphisms and one without.

\subsection{Multiple topological conjugacy classes in arbitrary degree $>2$}  \label{multiple}
\begin{figure}
\includegraphics[width=2.5in]{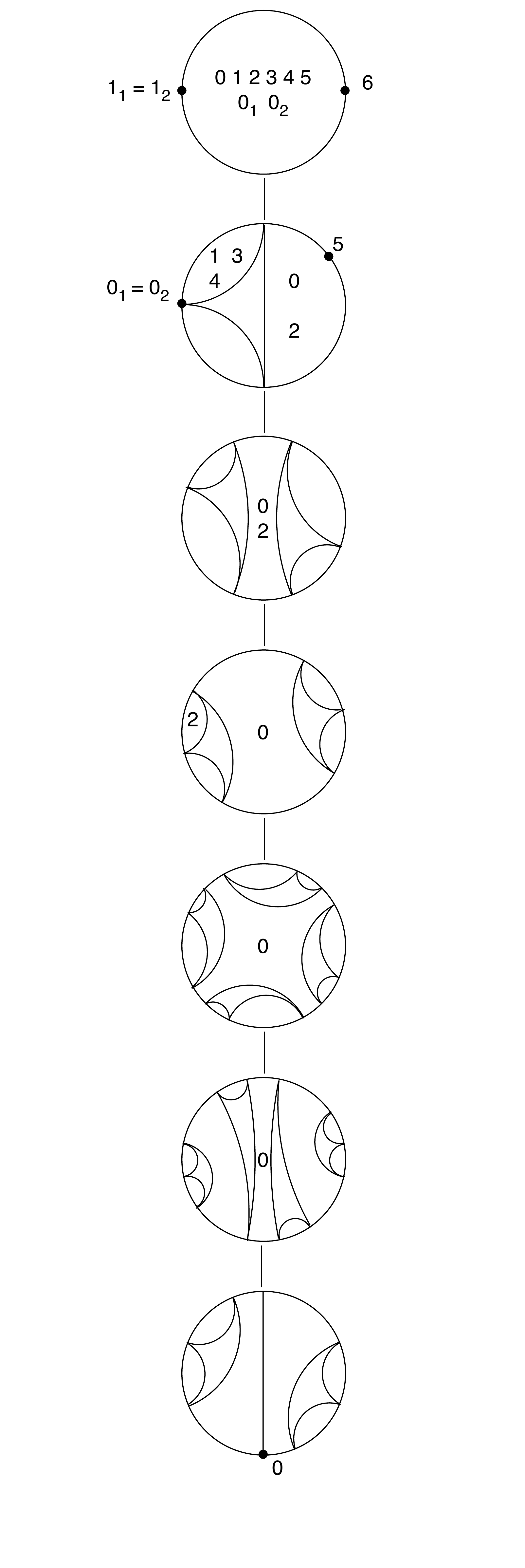}
\caption{A degree 4 pictograph determining two topological conjugacy classes; see \S\ref{multiple}.  The labels without subscript mark the orbit of the lowest critical point.}
\label{multiple twist class}
\end{figure}

Figure \ref{multiple twist class} shows a pictograph $\cD$ in degree 4 that determines two topological conjugacy classes.  This example can easily be generalized to any degree $d \geq 3$ by replacing the critical point at the highest branching vertex $v_0$ with one of multiplicity $d-2$.  It has one fundamental subannulus.

Let $v_0, v_{-1}, v_{-2}, \ldots$ denote the consecutive vertices in descending order.  To compute the number of topological conjugacy classes, we evaluate the twist periods at each level.  First, choose any restricted basin $(f, X_{t_0}(f))$.  The automorphism group $\Aut_{(\cF, \cX)}(f, X_{t_0}(f))$ is trivial.   As with every pictograph, $\Top(\cD,0) = 1$.  Because of the local symmetry at $v_{-1}$, there is only one conformal equivalence class of extension to $(f, X_{t_1}(f))$, so we also have $\Top(\cD, 1)$=1.  At $v_{-2}$, however, the symmetry is broken by the location of the second iterate of the lower critical point, so the two gluing choices determine distinct conformal equivalence classes.  The sum of relative moduli at $v_{-2}$ is $1/2 + 1/2 = 1$, so a full twist in the fundamental annulus induces a full twist at $v_{-2}$.  Consequently, the two conformal classes lie in two distinct topological conjugacy classes and $\Top(\cD, 2) = 2$.   

For each vertex below $v_{-2}$, there is a local symmetry of order at least 2, so the two gluing choices are conformally equivalent.  We conclude that there are exactly two conformal equivalence classes of basins extending the given $(f, X_{t_0}(f))$, and these lie in exactly two topological conjugacy classes.  Because this is the pictograph for polynomials in the shift locus, there are exactly two topological conjugacy classes of polynomials in $\cM_d$ with the given pictograph.

%%%%%%
%%%%%%

\bigskip\bigskip
%\bibliographystyle{../tex/bib/math}
%\bibliography{../tex/bib/math}

\def\cprime{$'$}

 \end{document}